\documentclass[11pt,leqno]{article}

\usepackage[english]{babel}

\usepackage{amsmath}

\usepackage{amssymb}

\usepackage{xypic}

\usepackage[all]{xy}

\usepackage{latexsym}

\usepackage{amsfonts}

\usepackage[dvips]{graphicx}

\usepackage[latin1]{inputenc}

\usepackage[T1]{fontenc}

\numberwithin{equation}{section}

\newcommand{\A}{\mathcal{A}}

\newcommand{\C}{\mathcal{C}}

\newcommand{\D}{\mathcal{D}}

\renewcommand{\mod}{\mathrm{Mod}}

\newcommand{\op}{\mathrm{Op}}

\newcommand{\rc}{\mathbb{R}\textrm{-}\mathrm{c}}

\newcommand{\Powl}{\hat{\mathcal{P}}(\Z_\ell)}

\newcommand{\CC}{\mathbb{C}}

\newcommand{\R}{\mathbb{R}}

\newcommand{\Z}{\mathbb{Z}}

\newcommand{\N}{\mathbb{N}}

\renewcommand{\P}{\mathcal{P}}

\newcommand{\OO}{\mathcal{O}}

\newcommand{\iso}{\stackrel{\sim}{\to}}

\newcommand{\OW}{\OO^\mathrm{w}}
\newcommand{\OWX}{\OO^\mathrm{w}_X}

\newcommand{\CW}{\C^{{\infty ,\mathrm{w}}}}
\newcommand{\wtens}{\overset{\mathrm{w}}{\otimes}}

\newcommand{\rh}{\mathit{R}\mathcal{H}\mathit{om}}

\newcommand{\ho}{\mathcal{H}\mathit{om}}

\newcommand{\id}{\mathrm{id}}

\renewcommand{\dim}{\textbf{Proof.}}

\newcommand{\qed}{\nopagebreak \phantom{} \hfill $\Box$ \\}

\newcommand{\supp}{\mathrm{supp}}

\newcommand{\RP}{\mathbb{R}^{{\scriptscriptstyle{+}}}}

\newcommand{\imin}[1]{#1^{-1}}

\newcommand{\lind}[1]{\underset{#1}{\underrightarrow{\lim}}}

\newcommand{\Lind}{\underrightarrow{\lim}}  

\newcommand{\lpro}[1]{\underset{#1}{\underleftarrow{\lim}}}

\newcommand{\exs}[3]{0 \to {#1} \to {#2} \to {#3} \to 0}

\newcommand{\lexs}[3]{0 \to {#1} \to {#2} \to {#3}}

\newcommand{\dt}[3]{{#1} \to {#2} \to {#3} \stackrel{+}{\to}}

\newtheorem{teo}{Theorem}[section]

\newtheorem{df}[teo]{Definition}

\newtheorem{cor}[teo]{Corollary}

\newtheorem{oss}[teo]{Remark}

\newtheorem{prop}[teo]{Proposition}

\newtheorem{lem}[teo]{Lemma}

\newtheorem{es}[teo]{Example}

\newtheorem{nt}[teo]{Notations}

\usepackage{latexsym}

\author{Naofumi Honda\ \ \ \ Luca Prelli}

\title{\bf{\sc{Multi-specialization and multi-asymptotic expansions}}}

\date{}

\begin{document}

\maketitle

\thispagestyle{empty}
\begin{abstract}
In this paper we extend the notion of specialization functor to the case of several closed submanifolds satisfying some suitable conditions. Applying this functor to the sheaf of Whitney holomorphic functions
we construct different kinds
 of sheaves of multi-asymptotically developable functions,
whose definitions are natural extensions of the definition of strongly asymptotically developable functions introduced by Majima.
\end{abstract}

\tableofcontents

\section*{Introduction}

Asymptotically developable expansions of  holomorphic functions on a sector are an
important tool to study ordinary differential equations
with irregular singularities.
When we study, in higher dimension,
a completely integrable connection with irregular singularities
along a normal crossing divisor $H = H_1 \cup \dots \cup H_\ell \subset X$,
it is known that these asymptotic expansions are too weak for this purpose.
Hence H.~Majima, in \cite{Ma84}, introduced the notion of {\it{strongly asymptotically
developable expansion}}
along $H$ for a holomorphic function defined on a poly-sector $S$, and the one of
{\it{consistent family of coefficients}} to which $f$ is
strongly asymptotically developable.

We can understand these notions
from a view point of a locally defined multi-action. For each smooth submanifold
$H_k$ ($k=1,2,\dots,\ell$), we can locally identify $X$ with the normal bundle
$T_{H_k} X$ of $H_k$ near $H_k$. A conic action on $T_{H_k} X$ by $\R^+$
induces a local action $\mu_k$ on $X$ near $H_k$. Then a poly-sector $S$
on which $f$ is defined can be regarded as {\it{a multi-cone}}
with respect to a multi-action $\mu_1$, \dots, $\mu_\ell$
in the sense that
it is an intersection of open sets $V_k$ ($k=1,2,\dots,\ell$) where each $V_k$
is a (locally) conic subset with respect to the action $\mu_k$ and its
edge is contained in $H_k$. A strongly asymptotically
developable expansion of $f$ is, roughly speaking, a formal Taylor
expansion with respect to an orbit
generated by these actions $\mu_1$, \dots, $\mu_\ell$.

Hence, from this point of view, one can expect that
strongly asymptotically developability extends to that along a more general
$H$. As a matter of fact, we have succeeded to construct
the sheaves of {\it{multi-asymptotically developable functions}}
along several kinds of $H$ by the aid of a multi-action,
whose definitions are natural extensions of the one introduced by H.~Majima in \cite{Ma84}.
These sheaves contain, as important cases, not only that
of strongly asymptotically developable functions
but also that associated with a multi-cone which appears in
{\it{a bi-normal deformation}} introduced by P.~Schapira and K.~Takeuchi in \cite{ST94}.

Now an important problem is to establish relations between these sheaves
on different spaces along different kinds of $H$,
to be more precise, we need to construct operations such as inverse images including
restrictions and direct images for these sheaves.
For that purpose, we need a uniform machinery allowing us to treat geometries
associated with these multi-actions. Namely, we need the notion of {\it{multi-normal deformation}} and the one of {\it{multi-specialization}} introduced in this paper,
which are our main subjects.

Let us briefly explain these new notions.

Let $X$ be a $n$-dimensional real analytic manifold with $\operatorname{dim}{X} = n$, and let
$\chi = \{M_1,\dots,M_\ell\}$ be
a family of connected closed submanifolds satisfying some suitable conditions. The multi-normal deformation of $X$ with respect to $\chi$ is constructed as follows. We first construct the normal deformation $\widetilde{X}_{M_1}$ of $X$ along $M_1$ defined by M.~Kashiwara and P.~Schapira in \cite{KS90}. Then, taking the pull-back of $M_2$ in $\widetilde{X}_{M_1}$, we can obtain the normal deformation $\widetilde{X}_{M_1,M_2}$ of $\widetilde{X}_{M_1}$ along  the pull-back of $M_2$. Then we can define recursively the normal deformation along $\chi$ as
$\widetilde{X}=\widetilde{X}_{M_1,\ldots,M_\ell}:=(\widetilde{X}_{M_1,\ldots,M_{\ell-1}})^\sim_{\widetilde{M}_\ell}.$ This manifold is of dimension $n+\ell$, it is locally isomorphic to $X \times \R^\ell$ and in the zero section $X \times \{0\}$ it is isomorphic to
$$
\underset{X,1\leq j\leq \ell}{\times}T_{M_j}\iota(M_j):=T_{M_1}\iota(M_1) \underset{X}{\times} T_{M_2}\iota(M_2) \underset{X}{\times} \cdots \underset{X}{\times} T_{M_\ell}\iota(M_\ell),
$$
where $\iota(M_j)$ denotes the intersection of the $M_k$'s strictly  containing $M_j$ (or $X$ if $M_j$ is maximal). There is also an action of $(\RP)^\ell$ on $\widetilde{X}$, which is obtained as a natural extension of the $\RP$-action of the normal deformation with respect to one submanifold. Its restriction to the zero section will be crucial for the definition of multi-asymptotic functions.
There are also natural notions of multi-cone and multi-normal cone extending the one of P.~Schapira and K.~Takeuchi, which will be the key to understand the geometry of the sections of the multi-specialization functor.
Given a morphism of real analytic manifolds $f:X \to Y$, we are also able to extend $f$ to a morphism $\widetilde{f}:\widetilde{X}\to\widetilde{Y}$. This is done by repeatedly employing the usual construction
of a morphism between normal deformations,
i.e. we extend $f$ to  $\widetilde{f}_1:\widetilde{X}_{M_1} \to \widetilde{Y}_{N_1}$,
then we extend $\widetilde{f}_1$ to $\widetilde{f}_{1,2}:\widetilde{X}_{M_1,M_2} \to \widetilde{Y}_{N_1,N_2}$.
Then we can define recursively $\widetilde{f}:\widetilde{X} \to \widetilde{Y}$
by extending the morphism $\widetilde{f}_{1,\ldots,\ell-1}:\widetilde{X}_{M_1,\ldots,M_{\ell-1}} \to \widetilde{Y}_{N_1,\ldots,N_{\ell-1}}$
to the normal deformations with respect to $M_\ell$ and $N_\ell$ respectively. This morphism enables us to make a link between different kinds of multi-normal deformations. As a kind of example, the desingularization map makes a link between the normal deformation with respect to a normal crossing divisor and the binormal deformation of P.~Schapira and K.~Takeuchi.

Once we have constructed the multi-normal deformation, we are able to extend the definition of the specialization functor to the case of several submanifolds. Roughly speaking, it is a functor associating to a subanalytic sheaf on $X$ a subanalytic sheaf on $\underset{X,1\leq j\leq \ell}{\times}T_{M_j}\iota(M_j)$. We perform all these constructions in the subanalytic setting in order to treat sheaves of functions with growth conditions. Given a morphism of real analytic manifolds $f:X \to Y$ we give conditions for the commutation of the direct and inverse images with respect to the multi-specialization functor.

When we apply the multi-specialization functor to the subanalytic sheaf of Whitney holomorphic functions we obtain the sheaf of multi-asymptotically developable functions along $\chi$ and, outside the zero section in the normal crossing case, the sheaf of strongly asymptotically developable functions of H.~Majima. When we apply the multi-specialization functor to the subanalytic sheaf of Whitney holomorphic functions vanishing up to infinity on $M_1 \cup \cdots \cup M_\ell$ (resp. Whitney holomorphic functions on $M_1 \cup \cdots \cup M_\ell$) we obtain the sheaf of flat multi-asymptotically developable functions (resp. consistent families of coefficients) along $\chi$. The vanishing of the $H^1$ of flat multi-asymptotically developable functions allows us to prove a Borel-Ritt exact sequence for multi-asymptotic functions.

The paper is organized in the following way. In Section \ref{1} we introduce the notion of multi-normal deformation. In Section \ref{2} we define the multi-action of $(\RP)^\ell$ on the zero section of $\widetilde{X}$. In Sections \ref{3} and \ref{4} we give the definitions of multi-cone and multi-normal cone which are essential to understand the sections of the multi-specialization. Morphisms between multi-normal deformations and their restrictions to the zero sections are studied in Section \ref{5}. The functorial construction is performed in Section \ref{6} with the definition of the multi-specialization functor and its relations with the functors of direct and inverse image. In Section \ref{7} we define
the sheaves of multi-asymptotically developable functions along $\chi$ whose functorial nature is proved in Section \ref{8}, where we apply the multi-specialization functor to the subanalytic sheaf of holomorphic functions with Whitney growth conditions.

We end this work with an Appendix in which we introduce the category of multi-conic sheaves. Using o-minimal geometry we construct suitable coverings of subanalytic open subsets which are helpful in order to study the sections of multi-conic sheaves.
\par For the reader convenience, we add some references. For basic notions on
sheaves and the definition of specialization we refer to \cite{KS90} and to \cite{ST94,Ta96} for the theory of bispecialization.  For subanalytic sheaf
theory and specialization in this setting we refer to \cite{KS01,Pr08,Pr07}. A
complete introduction  to subanalytic sets can be found in \cite{BM88} (see
also \cite{Co00,VD98} for the more general theory of o-minimal structures).
Strongly asymptotic expansions were introduced in \cite{Ma84}, we also recommend \cite{GS99,HS00,HS01} for further developments and \cite{Mo99,Mo99x} for Gevrey asymptotics.

\begin{section}{Multi-normal deformation}\label{1}

In this paper, all the manifolds are assumed to be countable at infinity.
Let $X$ be a real analytic manifold with $\operatorname{dim}{X} = n$, and let
$\chi = \{M_1,\dots,M_\ell\}$ be
a family of closed submanifolds in $X$ ($\ell \ge 1$).
We set, for $N \in \chi$ and $p \in N$,
$$
	\operatorname{NR}_p(N) := \{M_j \in \chi;\,
p \in M_j,\,N \nsubseteq M_j \text{ and } M_j \nsubseteq N\}.
$$
Let us consider the following conditions for $\chi$.
\begin{itemize}
\item[H1] Each $M_j \in \chi$ is connected and the
submanifolds are mutually distinct, i.e. $M_j \neq M_{j'}$ for $j \neq {j'}$.
\item[H2] For any $N \in \chi$ and $p \in N$ with $\operatorname{NR}_p(N) \ne \emptyset$,
we have
\begin{equation}{\label{eq:geometric-condition}}
\left(
\underset{
M_j \in \operatorname{NR}_p(N)
}{\bigcap} T_pM_j\right) + T_pN = T_pX.
\end{equation}
\end{itemize}
Note that, if $\chi$ satisfies the condition H2, 
the configuration of two submanifolds
must be either 1.~or 2.~below.
\begin{enumerate}
\item Both submanifolds intersect transversely.
\item One of them contains the other.
\end{enumerate}
We set, for $N \in \chi$,
$$
\iota_\chi(N) := \left\{
\begin{matrix}
&X \qquad & \text{There exists no $M_j\in \chi$ with $N \varsubsetneqq M_j$}, \\
\\
&\underset{N \varsubsetneqq M_j}{\bigcap} M_j  \qquad &\text{Otherwise}.
\end{matrix}
\right.
$$
When there is no risk of confusion, we write for short $\iota(N)$ instead of $\iota_\chi(N)$.

Since $T_{M_j}\iota (M_j)$ is contained in the zero section $T_XX$ if $M_j$ satisfies
$M_j = \iota(M_j)$, we assume the condition H3 below for simplicity.

\begin{itemize}
\item[H3] $M_j \ne \iota(M_j)$ for any $j \in \{1,2,\dots, \ell\}$.
\end{itemize}

\begin{es} Let $X=\R^3$ with coordinates $(x_1,x_2,x_3)$.
\begin{itemize}
\item[(i)] Let $\chi=\{M_1,M_2,M_3\}$ with $M_i=\{x_i=0\}$, $i=1,2,3$. Then clearly $\chi$ satisfies H1. We have $\iota(M_i)=X$, $i=1,2,3$ and $T_0M_i + T_0M_j = T_0X$, $i,j \in \{1,2,3\}$, $i \neq j$. Hence $\chi$ satisfies H2,H3.
\item[(ii)] Let $\chi=\{M_1,M_2,M_3\}$ with $M_1=\{0\}$, $M_2=\{x_2=x_3=0\}$, $M_3=\{x_3=0\}$. Then clearly $\chi$ satisfies H1. We have $\iota(M_1)=M_2$, $\iota(M_2)=M_3$, $\iota(M_3)=X$ and $NR_0(M_i)=\emptyset$ for $i=1,2,3$. Hence $\chi$ satisfies H2,H3.
\item[(iii)] Let $\chi=\{M_1,M_2,M_3\}$ with $M_1=\{0\}$, $M_2=\{x_1=0\}$, $M_3=\{x_2=0\}$. Then clearly $\chi$ satisfies H1. We have $NR_0(M_1)=\emptyset$ and $T_0M_2 + T_0M_3 = T_0X$. We have $\iota(M_1)=M_2 \cap M_3 \supsetneqq M_1$ and $\iota(M_2)=\iota(M_3)=X$. Hence $\chi$ satisfies H2,H3.
\item[(iv)] Let $\chi=\{M_1,M_2,M_3\}$ with $M_1=\{0\}$, $M_2=\{x_1=x_2=0\}$, $M_3=\{x_3=0\}$. Then clearly $\chi$ satisfies H1. We have $NR_0(M_1)=\emptyset$ and $T_0M_2 + T_0M_3 = T_0X$. We have $\iota(M_2)=\iota(M_3)=X$ and $\iota(M_1)=M_2 \cap M_3=M_1$. Hence $\chi$ satisfies H2 but not H3.
\item[(v)] Let $\chi=\{M_1,M_2\}$ with $M_1=\{x_1=x_2=0\}$, $M_2=\{x_2=x_3=0\}$. Then clearly $\chi$ satisfies H1. We have $T_0M_1+T_0M_2 \subsetneqq T_0X$. Then $\chi$ does not satisfy H2.
\item[(vi)]  Let $\chi=\{M_1,M_2,M_3\}$ with $M_1=\{x_1=x_2\}$, $M_2=\{x_1=0\}$, $M_3=\{x_2=0\}$. Then clearly $\chi$ satisfies H1. We have $T_0M_i+\bigcap_{i \neq j}T_0M_j=T_0M_i \subsetneqq T_0X$ for $i=1,2,3$. Then $\chi$ does not satisfy H2.
\item[(vii)] Let $\chi=\{M_1,M_2,M_3\}$ with $M_1=\{x_1=x_2^2\}$, $M_2=\{x_1=0\}$, $M_3=\{x_2=0\}$. Then clearly $\chi$ satisfies H1. We have $T_0M_1+\bigcap_{1 \neq j}T_0M_j=T_0M_1 \subsetneqq T_0X$. Then $\chi$ does not satisfy H2 (even if $\bigcap_{3\neq j}T_0M_j+T_0M_3=T_0X$).
\end{itemize}

\end{es}

\begin{prop}{\label{prop:canonical-coordinate}}
The following conditions are equivalent.
\begin{enumerate}
\item The family $\chi$ satisfies the condition H2. 
\item For any $p \in X$, there exist a neighborhood $V$ of $p$ in $X$, a system
of local coordinates $(x_1,\dots,x_n)$ in $V$ and a family of subsets $\{I_j\}_{j=1}^\ell$
of the set $\{1,2,\dots,n\}$ for which the following conditions hold.
\begin{enumerate}
\item Either $I_k \subset I_j$, $I_j \subset I_k$ or $I_k \cap I_j = \emptyset$ holds
$($$k,j \in \{1,2,\dots,\ell\}$$)$.
\item A submanifold $M_j \in \chi$ with $p \in M_j$ $($$j = 1,2,\dots,\ell$$)$
is defined by $\{x_i = 0; i \in I_j\}$ in $V$.
\end{enumerate}
\end{enumerate}
\end{prop}
\dim\ \
Clearly 2.~ implies 1. We will show the converse implication.
We may assume, by taking $V$ sufficiently small, that $p \in M_j$ or $V \cap M_j = \emptyset$ for any $j$.
We set, for $N \in \chi$,
$$
d(N) := \max\, \{k;\,N = M_{j_1} \varsubsetneqq M_{j_2} \varsubsetneqq \dots \varsubsetneqq M_{j_k},\, M_{j_1},\dots,M_{j_k} \in \chi\},
$$
and
$$
d(\chi) := \max \{d(N);\, N \in \chi\}.
$$
We show the proposition by induction with respect to $d(\chi)$.
Clearly the equivalence holds for a family $\chi$ with $d(\chi) = 1$.
Let us consider $\chi$ with $d(\chi) = \kappa > 1$, and
suppose that $1. \Rightarrow 2.$ were true for any $\chi$ with $d(\chi) < \kappa$.
Let $L_1$, $\dots$, $L_m$ denote the least elements in $\chi$ with respect
to the partial order $\subset$ of sets.

For each $k=1,2,\dots,m$,  we will determine defining functions $f_{k,1}$, $\dots$,
$f_{k, i_k}$ of the submanifold $L_k$ in the following way.
Set, for any $N \in \chi$,
$$
\chi_N = \{L \in \chi;\, N \varsubsetneqq L\}.
$$
Since we have $d(\chi_{L_k}) < \kappa$ ($k=1,2,\dots,m$),
and since the family $\chi_{L_k}$
also satisfies the condition 1.~ of the proposition,
by the induction hypothesis,
there exist local coordinate functions ($\varphi_{k,1}$, $\dots$, $\varphi_{k,n}$)
that satisfy the condition 2.~for the family $\chi_{L_k}$.
Then, for $k=1,2,\dots,m$, we take defining functions
$\{f_{k,1},\dots, f_{k,i_k}\}$
of $L_k$ so that
they contain all the coordinate functions
$\varphi_{k, i}$ which vanish on $L_k$.

As $L_{j'} \in \operatorname{NR}_p(L_{j})$ holds for $1 \le j \ne j' \le m$,
it follows from the condition \eqref{eq:geometric-condition}
that we have
$$
\underset{1\le k \le m,\, 1 \le i \le i_k}{\wedge} df_{k,i} \ne 0
\qquad \text{near $p$}.
$$
Therefore the family of these functions $\{f_{k,i}\}$
can be extended to a system of local coordinates near $p$, for which
we can easily verify 2.~of the proposition. This completes the proof.
\qed

\begin{oss}
As a consequence of the proposition, for the family $\chi$ satisfying
the conditions H1, H2 and H3,
each $M_j \in \chi$ is defined by linear equations $x_i = 0$
with $i \in I_j$ near $p \in M_1 \cap \dots \cap M_\ell$
where $I_j$'s satisfy
\begin{equation}\label{eq:conditions-indices-set}
\begin{aligned}
	&\text{(i) either }  I_j \subsetneqq I_k,\, I_k \subsetneqq I_j \text{ or }
	I_j \cap I_k = \emptyset \text{ holds for any $j \ne k$}, \\
	&\text{(ii) }  \left(\underset{I_k \subsetneqq I_j}{\bigcup}I_k \right)
	\subsetneqq I_j
	 \text{ for any $j$}.
\end{aligned}
\end{equation}
For these $I_j$'s, we set
\begin{equation}\label{eq:def-hat-I}
\hat{I}_j := I_j \setminus \left(\underset{I_k \subsetneqq I_j}{\bigcup}I_k \right)
\end{equation}
which is not empty by (ii) of (\ref{eq:conditions-indices-set}).
\end{oss}


Let $X$ be a $n$-dimensional real analytic manifold and let $\chi=\{M_1,\ldots,M_\ell\}$ be a family of closed smooth submanifolds of $X$ satisfying
H2. First recall the construction of the normal deformation of $X$ along $M_1$.
We denote it by $\widetilde{X}_{M_1}$ and we denote by $t_1\in\R$
the deformation parameter.
Let $\Omega_{M_1}=\{t_1>0\}$ and let us identify $\{t_1 = 0\}$
with $T_{M_1}X$. We have the commutative diagram
\begin{equation}\label{normal deformation}
\xymatrix{T_{M_1}X \ar[r]^{s_{M_1}} \ar[d]^{\tau_{M_1}} & \widetilde{X}_{M_1} \ar[d]^{p_{M_1}} &
\Omega_{M_1} \ar[l]_{\ \ i_{\Omega_{M_1}}} \ar[dl]^{\widetilde{p}_{M_1}} \\
M \ar[r]^{i_{M_1}} & X. & }
\end{equation}
Set $\widetilde{\Omega}_{M_1}=\{(x,t_1)\;;\;t_1\neq 0\}$ and define
$$
\widetilde{M}_2:=\overline{(p_{M_1}|_{\widetilde{\Omega}_{M_1}})^{-1}M_2}.
$$
Then $\widetilde{M}_2$ is a closed smooth submanifold of $\widetilde{X}_{M_1}$. Now we can define the normal deformation along $M_1,M_2$ as
$$
\widetilde{X}_{M_1,M_2}:=(\widetilde{X}_{M_1})^\sim_{\widetilde{M}_2}.
$$
Then we can define recursively the normal deformation along
$\chi$ as
$$
\widetilde{X}=\widetilde{X}_{M_1,\ldots,M_\ell}:=(\widetilde{X}_{M_1,\ldots,M_{\ell-1}})^\sim_{\widetilde{M}_\ell}.
$$
Set $S=\{t_1,\ldots,t_\ell=0\}$, $M=\bigcap_{j=1}^\ell M_j$ and $\Omega=\{t_1,\ldots,t_\ell>0\}$. Then we have the commutative diagram
\begin{equation}\label{multi normal deformation}
\xymatrix{S \ar[r]^s \ar[d]^\tau & \widetilde{X} \ar[d]^p &
\Omega \ar[l]_{\ \ i_\Omega} \ar[dl]^{\widetilde{p}} \\
M \ar[r]^{i_M} & X. & }
\end{equation}
\begin{oss}{\label{rem:coordinate-systems}}
By Proposition \ref{prop:canonical-coordinate}, there exists
a family $\{X^{(\alpha)}\}$ of subanalytic open subsets in $X$
satisfying
\begin{enumerate}
\item $\{X^{(\alpha)}\}$ covers $M$ and it is locally finite in $X$.
\item Each $X^{(\alpha)}$ is isomorphic to the whole $\mathbb{R}^n$
as a real analytic manifold. On the coordinates
$(x_1, \dots, x_n)$ of $\mathbb{R}^n$,
the submanifold $M_j$  is defined by $x_i = 0$
with $i \in I_j$ where $I_j$'s satisfy
the condition (a) of Proposition {\ref{prop:canonical-coordinate}}.
\end{enumerate}
In what follows, we replace $X$ with $\cup_{\alpha} X^{(\alpha)}$.
Therefore $\widetilde{X}$ has a locally finite subanalytic open covering
$\{\widetilde{X}^{(\alpha)} := p^{-1}(X^{(\alpha)})\}_\alpha$ and each local model
$\widetilde{X}^{(\alpha)}$ is isomorphic to the multi-normal deformation
$\widetilde{\mathbb{R}^n}$ along $M_j := \{x_i = 0,\, i\in I_j\}$'s.
For further details about gluing of these local models,
see the proof of Proposition \ref{prop:zero section}, in particular,
the equations (\ref{eq:gluing-rule}).
\end{oss}

Let $(x_1, \dots, x_n)$ (resp. $I_j$'s) be the coordinates
(resp. subsets of $\{1,\dots,n\}$)
of $X^{(\alpha)}$ given in Remark \ref{rem:coordinate-systems}.
Set, for $i \in \{1,\dots,n\}$,
\begin{equation}{\label{eq:def-J_i}}
J_i=\{j\in\{1,\ldots,\ell\}\;;\;i\in I_j\}, \ \ \ \ t_{J_i}=\prod_{j\in J_i}t_j,
\end{equation}
where $t_1,\ldots,t_\ell\in\R$ and $t_{J_i}=1$ if $J_i=\emptyset$.
Then $p:\widetilde{X}\to X$ is locally defined by
$$
(x_1,\ldots,x_n,t_1,\ldots,t_\ell) \mapsto (t_{J_1}x_1,\ldots,t_{J_n}x_n).
$$
We are interested in the bundle structure of the zero section
$$
S :=\{t_1=\ldots=t_\ell=0\} \subset \widetilde{X}.
$$
Let us consider the canonical map $T_{M_j}\iota(M_j) \overset{}{\to} M_j \hookrightarrow X$, $j=1,\dots,\ell$,
we write for short $$\underset{X,1\leq j\leq \ell}{\times}T_{M_j}\iota(M_j)$$ instead of
$$
T_{M_1}\iota(M_1) \underset{X}{\times} T_{M_2}\iota(M_2) \underset{X}{\times} \cdots \underset{X}{\times} T_{M_\ell}\iota(M_\ell).
$$
\begin{prop}\label{prop:zero section} Assume that $\chi$ satisfies the conditions H1, H2 and H3.
Then we have
\begin{equation}
S = \underset{X,1\leq j\leq \ell}{\times}T_{M_j}\iota(M_j).
\end{equation}
\end{prop}
\dim\ \
Let $\hat{X}$ be a copy of $X$, and $\varphi: X \to \hat{X}$
a local coordinate transformation near $p \in X$. We set $\hat{M}_j = \varphi(M_j)$.
We may assume that $X = {\mathbb R}^n$
(resp. $\hat{X} = {\mathbb R}^n$)
with coordinates  $(x_1,\dots,x_n)$
(resp. $(\hat{x}_1, \dots, \hat{x}_n)$), and $\varphi$ is given by
$\hat{x}_i = \varphi_i(x_1,\dots,x_n)$ ($i=1,2,\dots,n$).
Moreover, by Proposition \ref{prop:canonical-coordinate}, we may also suppose
that there exists a family of subsets $\{I_j\}_{j=1}^\ell$ of $\{1,2,\dots,n\}$
satisfying the condition (\ref{eq:conditions-indices-set})
for the both coordinate systems.

Recall that the set $J_k \subset \{1,\dots,\ell\}$ was defined
in (\ref{eq:def-J_i}).
Then, outside of $\{t_1t_2\ldots t_\ell = 0\}$, the coordinate transformation
between multi-normal deformations
$$
(x_1,\dots,x_n,t_1,\dots,t_\ell) \to
(\hat{x}_1,\dots,\hat{x}_n,\hat{t}_1,\dots,\hat{t}_\ell)
$$
is given by
\begin{equation}{\label{eq:gluing-rule}}
\left\{
\begin{aligned}
&\hat{t}_{J_k}\hat{x}_k = \varphi_k(t_{J_1}x_1,\, t_{J_2}x_2,\,\dots,\, t_{J_n}x_n)
\qquad (k = 1,2,\dots,n),\\
&\hat{t}_j = t_j \qquad (j = 1,2,\dots,\ell),
\end{aligned}
\right.
\end{equation}
which naturally extends to $\{t_1t_2\ldots t_\ell=0\}$ (see below for the reason).
For any $k\in \underset{1 \le j \le \ell}{\bigcup} I_j$, we set
$$
I(k) := \underset{k \in I_j}{\bigcap} I_j = \underset{j \in J_k}{\bigcap} I_j.
$$
As the condition 2.~(a) of Proposition \ref{prop:canonical-coordinate} is equivalently
saying that
$$
I_{j} \cap I_{j'} \ne \emptyset \Rightarrow
I_{j} \subset I_{j'} \text{ or } I_{j'} \subset I_{j},
$$
the set $\{I_j;\, k \in I_j\}$
is totally ordered with respect to ``$\subset$``,
and $I(k)$ is its minimal element.
Hence, for any $k \in \underset{1 \le j \le \ell}{\bigcup} I_j$,
there exists $j(k) \in \{1,2,\dots,\ell\}$ such that $I(k) = I_{j(k)}$.
By expanding $\varphi_k(x_1,\dots,x_n)$ along the submanifold $M_{j(k)}$,
we obtain
$$
\varphi_k(x_1,\dots,x_n) =
\sum_{i \in I_{j(k)}}\left.\frac{\partial \varphi_k}{\partial x_i}\right|_{M_{j(k)}} x_i
+ {1\over 2}\sum_{i_1,i_2 \in I_{j(k)}}
\left.
\frac{\partial^2 \varphi_k}{\partial x_{i_1}\partial x_{i_2}}\right|_{M_{j(k)}}
x_{i_1}x_{i_2}
+ \dots,
$$
as $\varphi_k\vert_{M_{j(k)}} = 0$ holds.
Then we get
$$
t_{J_k}\hat{x}_k =
\sum_{i \in I_{j(k)}}\left.\frac{\partial \varphi_k}{\partial x_i}\right|_{M_{j(k)}}
t_{J_i}x_i
+ {1\over 2}\sum_{i_1,i_2 \in I_{j(k)}}
\left.
\frac{\partial^2 \varphi_k}{\partial x_{i_1}\partial x_{i_2}}\right|_{M_{j(k)}}
t_{J_{i_1}}t_{J_{i_2}}x_{i_1}x_{i_2}
+ \dots.
$$
We can easily see $J_k \subseteq J_i$ for any $i \in I_{j(k)}$
(in particular, the $\hat{x}_k$ extends to $t_1t_2 \ldots t_\ell = 0$).
Indeed,
this follows from the facts
$$
\begin{aligned}
&\text{$j \in J_k$ and $i \in I_{j(k)}$} \Longrightarrow
\text{$k \in I_j$ and $i \in I(k) = \underset{k \in I_\beta}{\bigcap}I_\beta$} \\
&\Longrightarrow i \in I_j \Longrightarrow j \in J_i.
\end{aligned}
$$
Therefore, by letting $t \to 0$,  we have
$$
\hat{x}_k =
\sum_{i \in I_{j(k)},\, J_k = J_i}
\left.\frac{\partial \varphi_k}{\partial x_i}\right|_M
x_i.
$$
where $M := \underset{1 \le j \le l}{\bigcap} M_j$.
Now we claim, for $k \in \underset{1 \le j \le \ell}{\bigcup} I_j$,
$$
\left\{i \in I_{j(k)};\, J_k = J_i\right\} =
\left\{i \in \{1,2,\dots,n\};\, j(k) = j(i)\right\},
$$
which is proved in the following way:
We first prove the implication ($\Leftarrow$).
Assume that $i$ satisfies $j(k) = j(i)$, which implies
$i \in I_{j(i)} = I_{j(k)}$. We have already proved the fact
$J_k \subseteq J_i$ for $i \in I_{j(k)}$.
Therefore it suffices to show $J_i \subseteq J_k$.
Let $\beta \in J_i$. Then, as $i \in I_\beta$, we have
$k \in I_{j(k)} = I_{j(i)} \subset I_\beta$, from which $\beta \in J_k$ follows.

The converse implication comes from
$$
I_{j(k)} = \underset{\beta \in J_k}{\bigcap}I_\beta
= \underset{\beta \in J_i}{\bigcap}I_\beta = I_{j(i)}
$$
as $j(k) = j(i)$ is equivalent to saying $I_{j(k)} = I_{j(i)}$
by the condition H1.

We divide the set $\underset{1 \le j \le \ell}{\bigcup} I_j \subset \{1,2,\dots,n\}$
into equivalent classes $\{B_\alpha\}$
by the equivalence relation ``$i \sim k \iff j(i) = j(k)$''.
Then, for an equivalent class $B$, we obtain
$$
\hat{x}_i =
\sum_{k \in B}
\left.\frac{\partial \varphi_i}{\partial x_k}\right|_M
x_k \qquad \text{for $i \in B$}.
$$
We denote by $E_B$ the vector bundle over $M$ defined by
the above equations where its fiber coordinates are given by $x_k$'s ($k \in B$).
Then, by the above observation,
$S$ is a direct sum of bundles $E_B$'s over $M$.
Note that, since each equivalent class can be written in the form
$$
I_{j(k)} \setminus \left(\underset{I_{j(i)} \varsubsetneqq I_{j(k)}}{\bigcup} I_{j(i)}\right)
=
I_{j(k)} \setminus \left(\underset{I_j \varsubsetneqq I_{j(k)}}{\bigcup} I_j\right)
$$
for some $j(k) \in \{1,2,\dots,l\}$ with $k \in \underset{1 \le j \le \ell}{\bigcup} I_j$, the bundle $E_B$ is isomorphic to
$T_{M_{j(k)}} \iota(M_{j(k)}) \underset{X}{\times} M$.
For any $j \in \{1,2,\dots,\ell\}$, the set
$
\hat{I}_j
$
defined by (\ref{eq:def-hat-I})
is not empty and it gives an equivalent class
of $\underset{1 \le j \le \ell}{\bigcup} I_j / \sim$. Further it follows from the
condition (\ref{eq:conditions-indices-set})
that we get
$$
\underset{1 \le j \le \ell}{\bigcup} I_j = \hat{I}_1 \sqcup \dots \sqcup \hat{I}_l.
$$
Hence we have obtained that $S$ is a direct sum of bundles
$T_{M_j} \iota(M_j) \underset{X}{\times} M$
($j=1,2,\dots.\ell$). This completes the proof.
\qed

\begin{es} Let us see three typical examples of multi-normal deformations.
\begin{enumerate}
\item (Majima) Let $X={\mathbb C}^2$ ($\simeq\R^4$ as a real manifold)
with coordinates $(z_1, z_2)$ and let $\chi=\{M_1,M_2\}$ with $M_1 = \{z_1 = 0\}$ and $M_2 = \{z_2 = 0\}$. Then $\chi$ satisfies H1, H2 and H3.
We have $I_1=\{1\}$, $I_2=\{2\}$, $J_1=\{1\}$, $J_2=\{2\}$ (in $\R^4$, if $z_1=(x_1,x_2)$ and $z_2=(x_3,x_4)$ we have $I_1=\{1,2\}$, $I_2=\{3,4\}$, $J_1=J_2=\{1\}$, $J_3=J_4=\{2\}$).
The map $p:\widetilde{X} \to X$ is defined by
$$
(z_1,z_2,t_1,t_2) \mapsto (t_1z_1,t_2z_2).
$$
Remark that the deformation is real though $X$ is complex. In particular $t_1, t_2 \in \R$.
We have $\iota(M_1)=\iota(M_2)=X$ and then the zero section $S$ of $\widetilde{X}$ is isomorphic to $T_{M_1}X \underset{X}{\times} T_{M_2}X$.

\item (Takeuchi)
Let $X={\mathbb R}^3$ with coordinates $(x_1, x_2, x_3)$ and let $\chi=\{M_1,M_2,M_3\}$ with
$M_1 = \{0\}$, $M_2 = \{x_2 = x_3 = 0\}$ and $M_3 = \{x_3 = 0\}$. Then $\chi$ satisfies H1, H2 and H3. We have $I_1=\{1,2,3\}$, $I_2=\{2,3\}$, $I_3=\{3\}$, $J_1=\{1\}$, $J_2=\{1,2\}$, $J_3=\{1,2,3\}$.
The map $p:\widetilde{X} \to X$ is defined by
$$
(x_1,x_2,x_3,t_1,t_2,t_3) \mapsto (t_1x_1,t_1t_2x_2,t_1t_2t_3x_3).
$$
We have $\iota(M_1)=M_2$, $\iota(M_2)=M_3$, $\iota(M_3)=X$ and then the zero section $S$ of $\widetilde{X}$ is isomorphic to $T_{M_1}M_2 \underset{X}{\times} T_{M_2}M_3 \underset{X}{\times} T_{M_3}X$.

\item (Mixed)
Let $X={\mathbb R}^3$ with coordinates $(x_1, x_2, x_3)$ and let $\chi=\{M_1,M_2,M_3\}$ with
$M_1 = \{0\}$, $M_2 = \{x_2 = 0\}$ and $M_3 = \{x_3 = 0\}$. Then $\chi$ satisfies H1, H2 and H3. We have $I_1=\{1,2,3\}$, $I_2=\{2\}$, $I_3=\{3\}$, $J_1=\{1\}$, $J_2=\{1,2\}$, $J_3=\{1,3\}$.
The map $p:\widetilde{X} \to X$ is defined by
$$
(x_1,x_2,x_3,t_1,t_2,t_3) \mapsto (t_1x_1,t_1t_2x_2,t_1t_3x_3).
$$
We have $\iota(M_1)=M_2 \cap M_3$, $\iota(M_2)=\iota(M_3)=X$ and then the zero section $S$ is isomorphic to $T_{M_1}(M_2 \cap M_3) \underset{X}{\times} T_{M_2}X \underset{X}{\times} T_{M_3}X$.
\end{enumerate}
\end{es}

The following easy lemma is needed later as we work on the subanalytic site.
Recall the definitions of $X^{(\alpha)} \subset X$ and the local model $\widetilde{X}^{(\alpha)} \subset \widetilde{X}$
given in Remark \ref{rem:coordinate-systems}.
\begin{lem}\label{lem:subanalytic-local-model}
	Let $U$ be a subanalytic subset in the local model $\widetilde{X}^{(\alpha)}$.
	If $p(U)$ is relatively compact in $X^{(\alpha)}$, then $U$ is subanalytic in $\widetilde{X}$.
\end{lem}
\dim \ It suffices to show that $U$ is subanalytic near $\tilde{x}$ for any $\tilde{x} \in \overline{U}$. Here the closure is taken in $\widetilde{X}$.
As $p(U)$ is relatively compact in $X^{(\alpha)}$,
we have $p(\tilde{x}) \in \overline{p(U)} \subset X^{(\alpha)}$.
This implies $\tilde{x} \in p^{-1}(X^{(\alpha)}) = \widetilde{X}^{(\alpha)}$.
Hence $U$ is subanalytic near $\tilde{x}$.
\qed
\end{section}

\begin{section}{Multi-actions}{\label{2}}



Let $X$ be a real analytic manifold and let $\chi = \{M_1, M_2, \dots, M_\ell\}$ be closed submanifolds. In what follows, we always assume that $\chi$ satisfies the conditions H1, H2 and H3.
Consider the diagram \eqref{multi normal deformation}.
There is a $(\RP)^\ell$ action
$$
\mu:\widetilde{X} \times (\RP)^\ell \to \widetilde{X}
$$
which is described in local coordinate system by
$$
((x_1,\ldots,x_n,t_1,\ldots,t_\ell),(c_1,\ldots,c_\ell)) \mapsto \left(c_{J_1}x_1,\ldots,c_{J_n}x_n,{t_1\over c_1},\ldots,{t_\ell\over c_\ell}\right).
$$
More precisely, the $j$-th component of the action is given by
$$
\mu_j:((x_1,\ldots,x_n,t_1,\ldots,t_\ell),c_j) \mapsto \left(c_{1j}x_1,\ldots,c_{nj}x_n,t_1,\ldots,{t_j\over c_j},\ldots,t_\ell\right),
$$
where $c_{ij}=c_j$ if $i\in I_j$ and $c_{ij}=1$ otherwise.



\begin{oss}{\label{rem:stablity-actions}}
These multi-actions are stable on each local mode $\widetilde{X}^{(\alpha)}$.
For the definition of $\widetilde{X}^{(\alpha)}$, see Remark
\ref{rem:coordinate-systems}.
\end{oss}
The zero section $S$
has the natural actions induced from a direct sum of vector bundles.
Hence we consider the actions on $\widetilde{X}$ which are compatible
with those on $S$.

Set
$$
\begin{aligned}
J_{\supsetneq M_j} &:= \{\alpha \in \{1,2,\dots,\ell\};\,  M_\alpha \supsetneq M_j,\,
\text{ there is no $\beta$ with $M_\alpha \supsetneq M_\beta \supsetneq M_j$ } \}\\
J_{\subsetneq M_j} &:= \{\alpha \in \{1,2,\dots,\ell\};\,  M_\alpha \subsetneq M_j,\,
\text{ there is no $\beta$ with $M_\alpha \subsetneq M_\beta \subsetneq M_j$ } \}
\end{aligned}
$$
Note that, by the conditions H1 and H2, the set $J_{\subsetneq M_j}$ either is empty or consists of only one index.


Using the actions $\mu_j$
for $j = 1,2, \dots, \ell$, we define the action
$\tau_j(\lambda): \widetilde{X} \to \widetilde{X}$
by
$$
\tau_j(\lambda) :=
\mu_j(\lambda) \underset{\beta \in J_{\supsetneq M_j}}{\prod} \mu_\beta(\lambda^{-1})
\qquad (j = 1,2,\dots,\ell).
$$
One can easily check that, on the zero section
$S$ of $\widetilde{X}$,
the action $\tau_j(\lambda)$ coincides with
$$
\left.\tau_j(\lambda)\right|_{T_{M_{\alpha}}\iota (M_{\alpha})} =
\left\{
\begin{matrix}
&T_{M_{\alpha}}\iota (M_{\alpha}) \overset{\lambda \cdot}{\to} T_{M_{\alpha}}\iota (M_{\alpha}) \qquad
&(\alpha = j) \\
&\operatorname{id}_{T_{M_{\alpha}}\iota (M_{\alpha})} \qquad
&(\alpha \ne j).
\end{matrix}
\right.
$$

Conversely, we can recover the original actions $\{\mu_j\}$ by $\{\tau_j\}$
in the following way.
$$
\mu_j(\lambda) = \underset{M_j \subseteq M_\alpha}{\prod} \tau_\alpha(\lambda)
\qquad (j=1,2,\dots,l).
$$
Hence both multi-actions on $\widetilde{X}$ associated with $\{\mu_j\}$ and $\{\tau_j\}$ are equivalent.

\begin{es} Let us see some example of multi-actions.
\begin{enumerate}
\item (Majima) Let $X={\mathbb C}^2$
with coordinates $(z_1, z_2)$ and let $M_1 = \{z_1 = 0\}$ and $M_2 = \{z_2 = 0\}$. Then
\begin{eqnarray*}
\mu:((z_1,z_2,t_1,t_2),(c_1,c_2)) & \mapsto & \left(c_1z_1,c_2z_2,{t_1 \over c_1},{t_2 \over c_2}\right) \\
\tau_1(\lambda):(z_1,z_2,t_1,t_2) & \mapsto & (\lambda z_1,z_2,\imin \lambda t_1,t_2) \\
\tau_2(\lambda):(z_1,z_2,t_1,t_2) & \mapsto & (z_1,\lambda z_2,t_1,\imin \lambda t_2).
\end{eqnarray*}

\item (Takeuchi)
Let $X={\mathbb R}^3$ with coordinates $(x_1, x_2, x_3)$ and let
$M_1 = \{0\}$, $M_2 = \{x_2 = x_3 = 0\}$ and $M_3 = \{x_3 = 0\}$. Then
\begin{eqnarray*}
\mu:((x_1,x_2,x_3,t_1,t_2,t_3),(c_1,c_2,c_3)) & \mapsto & \left(c_1x_1,c_1c_2x_2,c_1c_2c_3x_3,{t_1 \over c_1},{t_2 \over c_2},{t_3 \over c_3}\right) \\
\tau_1(\lambda):(x_1,x_2,x_3,t_1,t_2,t_3) & \mapsto & (\lambda x_1,x_2,x_3,\imin \lambda t_1,\lambda t_2,t_3) \\
\tau_2(\lambda):(x_1,x_2,x_3,t_1,t_2,t_3) & \mapsto & (x_1,\lambda x_2,x_3,t_1,\imin \lambda t_2,
 \lambda t_3) \\
\tau_3(\lambda):(x_1,x_2,x_3,t_1,t_2,t_3) & \mapsto & (x_1,x_2, \lambda x_3,t_1,t_2,
\imin \lambda t_3).
\end{eqnarray*}
\item (Mixed)
Let $X={\mathbb R}^3$ with coordinates $(x_1, x_2, x_3)$ and let
$M_1 = \{0\}$, $M_2 = \{x_2 = 0\}$ and $M_3 = \{x_3 = 0\}$. Then
\begin{eqnarray*}
\mu:((x_1,x_2,x_3,t_1,t_2,t_3),(c_1,c_2,c_3)) & \mapsto & \left(c_1x_1,c_1c_2x_2,c_1c_3x_3,{t_1 \over c_1},{t_2 \over c_2},{t_3 \over c_3}\right) \\
\tau_1(\lambda):(x_1,x_2,x_3,t_1,t_2,t_3) & \mapsto & (\lambda x_1,x_2,x_3,\imin \lambda t_1,
 \lambda t_2,\lambda t_3) \\
\tau_2(\lambda):(x_1,x_2,x_3,t_1,t_2,t_3) & \mapsto & (x_1, \lambda x_2,x_3,t_1,
\imin \lambda t_2,t_3) \\
\tau_3(\lambda):(x_1,x_2,x_3,t_1,t_2,t_3) & \mapsto & (x_1,x_2, \lambda x_3,t_1,t_2,
\imin \lambda t_3).
\end{eqnarray*}
\end{enumerate}
\end{es}

\end{section}


\begin{section}{Multi-cones}\label{3}

Let $q \in \underset{1 \le j \le \ell}{\bigcap} M_j$ and
$p_j = (q;\, \xi_j)$ be a point in $T_{M_j}\iota(M_j)$ ($j=1,2,\dots,\ell$).
We set
$$
p = p_1 \underset{X}{\times} \dots \underset{X}{\times} p_\ell
\in S = \underset{X,1\leq j \leq \ell}{\times}T_{M_j}\iota(M_j),
$$
and $\tilde{p}_j = (q;\, \tilde{\xi}_j) \in T_{M_j} X$ designates the
image of the point $p_j$ by the canonical embedding $T_{M_j} \iota (M_j) \hookrightarrow T_{M_j} X$.
We denote by $\operatorname{Cone}_{\chi, j}(p)$ ($j=1,2,\dots,\ell$)
the set of open cones in
$(T_{M_j}X)_q \simeq {\R}^{n - {\rm dim} M_j}$ that contain the point
$\tilde{\xi}_j \in (T_{M_j} X)_q \simeq {\R}^{n - {\rm dim} M_j}$.

%
%
\begin{df}
We say that an open set $G\subset (TX)_q$ is a multi-cone along $\chi$
with direction to
$p \in S_q$ $($the fiber of $S$ over $q$$)$ if
$G$ is written in the form
$$
G = \underset{1\leq j \leq\ell}{\bigcap} \pi_{j,\,q}^{-1}(G_j), \qquad
G_j \in \operatorname{Cone}_{\chi,j}(p)
$$
where $\pi_{j,\,q}: (TX)_q \to (T_{M_j}X)_q$ is the canonical projection.
We denote by $\operatorname{Cone}_\chi(p)$ the set of multi-cones along $\chi$
with direction to $p$.
\end{df}

\begin{es} We now give three typical examples of multi-cones.
\begin{enumerate}
\item (Majima) Let $X={\mathbb C}^2$
with coordinates $(z_1, z_2)$ and let $M_1 = \{z_1 = 0\}$ and $M_2 = \{z_2 = 0\}$.
Then $\operatorname{Cone}_{\chi}(p)$ for $p = (0,0;\,1,1)$ is nothing but
the set of multi sectors along $Z_1 \cup Z_2$ with their direction to $(1,1)$.
\item (Takeuchi)
Let $X={\mathbb R}^3$ with coordinates $(x_1, x_2, x_3)$ and let
$M_1 = \{0\}$, $M_2 = \{x_2 = x_3 = 0\}$ and $M_3 = \{x_3 = 0\}$.
For $p = (0,0,0;\, 1,1,1) \in T_{M_1}M_2 \underset{X}{\times} T_{M_2}M_3 \underset{X}{\times} T_{M_3}X$,
it is easy to see that
a cofinal set of $\operatorname{Cone}_{\chi}(p)$ is, for example, given by
the family of the sets
$$
\{(\xi_1, \xi_2, \xi_3);\,
\vert \xi_2 \vert + \vert \xi_3 \vert < \epsilon \xi_1,\,
\vert \xi_3 \vert < \epsilon \xi_2,\,
\xi_3 > 0
\}_{\epsilon > 0}.
$$
\item (Mixed)
Let $X={\mathbb R}^3$ with coordinates $(x_1, x_2, x_3)$ and let
$M_1 = \{0\}$, $M_2 = \{x_2 = 0\}$ and $M_3 = \{x_3 = 0\}$.
For $p = (0,0,0;\, 1,1,1) \in T_{M_1}(M_2 \cap M_3) \underset{X}{\times} T_{M_2}X \underset{X}{\times} T_{M_3}X$,
a cofinal set of $\operatorname{Cone}_{\chi}(p)$ is, for example, given by
the family of the sets
$$
\{(\xi_1, \xi_2, \xi_3);\,
\vert \xi_2 \vert + \vert \xi_3 \vert < \epsilon \xi_1,\,
\xi_2 > 0,\,
\xi_3 > 0
\}_{\epsilon > 0}.
$$
\end{enumerate}
\end{es}

\end{section}

\begin{section}{Multi-normal cones}\label{4}

\begin{df} Let $Z$ be a subset of $X$. The multi-normal cone to $Z$ along $\chi$ is the set
$$
C_\chi(Z)=\overline{\imin{\widetilde{p}}(Z)} \cap S.
$$
\end{df}

For any $q \in X$,  there exists an isomorphism
$$
\psi: X \iso (TX)_q
$$
near $q$ which satisfies $\psi(q) = (q;\,0)$ and $\psi(M_j) = (TM_j)_q$
for any $j=1,\dots,\ell$.
The existence of such a $\psi$ is guaranteed by Proposition \ref{prop:canonical-coordinate}.

\begin{lem}\label{lem:multi cone}
Let $Z$ be a subset of $X$.
We have the following equivalence:
$p \notin C_\chi(Z)$
if and only if
there exist an open subset $\psi(q) \in U \subset (TX)_q$
and a multi-cone $G \in \operatorname{Cone}_{\chi}(p)$ such that
$$
\psi(Z) \cap G \cap U = \emptyset
$$
holds.
\end{lem}
\dim\ \ The problem is local. Hence we identify $X \simeq (TX)_q \simeq {\mathbb R}^n$ by $\psi$, and
we consider, in what follows, $\operatorname{Cone}_{\chi}(p)$ as the set of cones defined in $X$.
We may assume $q=0$.
Recall that $\hat{I}_j \subset \{1,2,\dots,n\}$ is defined
by (\ref{eq:def-hat-I}).
Set
$B := \underset{1\leq j \leq \ell}{\bigcup} I_j =
\underset{1\leq j \leq \ell}{\sqcup} \hat{I}_j
\subset \{1,2,\dots,n\}.
$
Then the variables $x_i$'s with $i \in B$ give the fiber coordinates of
$S$. It suffices to show the following claim
$$
p \in C_\chi(Z) \iff
\text{For any $p \in U$ and $G \in \operatorname{Cone}_{\chi}(p)$,  }\,
\psi(Z) \cap G \cap U \neq \emptyset.
$$

\

\par\noindent First we will show ($\Rightarrow$):
Assume that
$$
p = (q;\, \{\xi_i\}_{i \in B})
= \underset{X, 1\leq j \leq \ell}{\times} p_j \in C_\chi(Z)
$$
where $p_j$ is a point in $T_{M_j} \iota (M_j)$.
By the definition, we have a sequence
$$
p^{(m)} =
(x_1^{(m)}, \dots, x_n^{(m)}, t_1^{(m)}, \dots, t_\ell^{(m)}) \in
\widetilde{p}^{-1}(Z) \subset \widetilde{X}
\qquad m=1,2, \dots
$$
satisfying that $x_i^{(m)} \to \xi_i$ for $i \in B$,  $x_i^{(m)} \to 0$ for $i \notin B$,
and $t_j^{(m)} \to 0$ for any $j$.
For $j \in \{1,2,\dots,\ell\}$, let us consider the commutative diagram
$$
\begin{matrix}
&& S
& \to  & T_{M_j}X \\
&&\downarrow & & \downarrow \\
&&\varphi:\widetilde{X} & \to &    \widetilde{X}_{M_{j}} \\
&&                         & \searrow & \downarrow \\
&&                         &           & X
\end{matrix}
$$
where the top horizontal arrow is given by the composition of morphisms
$$
S
\to T_{M_j}\iota(M_j) \underset{X}{\times} M
\hookrightarrow   T_{M_j}X
$$
and $\varphi: (x_1, \dots, x_n, t_1, \dots, t_\ell) \to
(\hat{x}_1, \dots, \hat{x}_n, \hat{t})$ is defined by
$$
\begin{array}{llll}
&\hat{t} &= \underset{M_\beta \subset M_j}{\prod} t_\beta, &\\
&\hat{x}_i &= t_{J_i} x_i &\qquad (i \notin I_j),\\
&\hat{x}_i &=
\left(\underset{\{\beta \in J_i;\,\, M_j \subsetneq M_\beta\}}{\prod} t_\beta\right) x_i
&\qquad(i \in I_j).
\end{array}
$$
Here the definitions of $J_i$, etc. were given in the proof
of Proposition \ref{prop:canonical-coordinate}.
Then $\varphi(p^{(m)})$ converges to $\overline{p}_j \in T_{M_j} X$ where
$\overline{p}_j$ is the image of $p_j$ by the canonical embedding $T_{M_j} \iota (M_j) \hookrightarrow T_{M_j} X$.
This implies $\overline{p}_j \in C_{M_j} (Z)$. Therefore it follows from the definition of a usual normal cone
that, for any cone $G_j \in \operatorname{Cone}_{\chi, j}(p)$,  we have
$\widetilde{p}_{M_j}(\varphi(p^{(m)})) \in \pi_{j,\,q}^{-1}(G_j)$ for any sufficiently large $m$.
This entails that, for any multi-cone $G \in \operatorname{Cone}_{\chi}(p)$,
we have $\widetilde{p}(p^{(m)}) \in G$ for any sufficiently large $m$, in particular,
we have $Z \cap G \cap U \neq \emptyset$ for any $U$.

\

\par\noindent Let us  show converse ($\Leftarrow$):  Set, for a subset $I \subset \{1,2,\dots, n\}$,
$$
\vert x \vert_I := \sum_{i \in I} \vert x_i \vert.
$$
Note that, if $I$ is empty,  we set $\vert x \vert_I = 1$.
Let
$$
p = (q;\, \xi) =
\underset{X,1\leq j \leq \ell}{\times} p_j
\in S
$$
where $p_j = (q;\, \xi_j) \in T_{M_j} \iota(M_j)$.
By the conic actions $\tau_j(\cdot)$, we may assume either
$\vert \xi_j \vert = \vert \xi \vert_{\hat{I}_j} = 1$
or
$\vert \xi_j \vert = \vert \xi \vert_{\hat{I}_j} = 0$ for $1\leq j \leq \ell$.

Let $\{G_j^{(m)}\}_{m=1,2,\dots}$ be a cofinal set of
$\operatorname{Cone}_{\chi,j}(p)$
and $\{U^{(m)}\}_{m=1,2,\dots}$ a set of fundamental neighborhoods of $q$. We set
$$
G^{(m)} := \underset{1\leq j \leq \ell}{\bigcap} \pi_{j,\,q}^{-1} (G_j^{(m)}) \qquad m = 1,2, \dots.
$$
Choose points  in $X$ as
$$
q^{(m)} \in Z \cap G^{(m)} \cap U^{(m)} \subset X \qquad m = 1,2, \dots,
$$
and define a sequence in $\widetilde{X}$ by
$
p^{(m)} := (q^{(m)};\, 1,\dots, 1).
$
Clearly we have $\widetilde{p}(p^{(m)}) \in Z \cap G^{(m)} \cap U^{(m)}$.
For each $1\leq j \leq \ell$, by taking a subsequence of $\{q^{(m)}\}$, we may assume
either $\vert q^{(m)} \vert_{\hat{I}_j} \ne 0$  for every $m$ or
$\vert q^{(m)} \vert_{\hat{I}_j} = 0$ for all $m$.
We divide the set $\{1,2,\dots,\ell\}$ into two sets $J'$ and $J''$ as follows:
$$
\begin{aligned}
J' &= \left\{j \in \{1,2,\dots,l\};\,
\vert \xi_j \vert = \vert \xi \vert_{\hat{I}_j} \ne 0\right\}, \\
J'' &= \{1,2,\dots, \ell\} \setminus J'.
\end{aligned}
$$
Note that $\vert q^{(m)} \vert_{\hat{I}_j} \ne 0$ ($m=1,2,\dots$) holds for $j \in J'$.
Let us determine a sequence $\kappa^{(m)} = (\kappa_1^{(m)}, \dots, \kappa_\ell^{(m)})$
of positive real numbers that satisfies the following conditions.
\begin{enumerate}
\item $\kappa^{(m)} \to 0$,
	$\kappa^{(m)}_j = \vert q^{(m)} \vert_{\hat{I}_j}$ for $j \in J'$ and
$
\displaystyle\frac{\vert q^{(m)} \vert_{\hat{I}_j}}{\kappa^{(m)}_j} \to 0
$
for $j \in J''$.
\item For any pair $\alpha,\beta \in \{1,2,\dots,\ell\}$ with
$M_\alpha \subsetneq M_\beta$,
we have
$
\displaystyle\frac{\kappa^{(m)}_\beta} {\kappa^{(m)}_\alpha} \to 0.
$
\end{enumerate}
Set, for $j \in \{1,2,\dots,\ell\}$,
$$
d_j := \max \{k; M_j = M_{j_0} \subsetneq M_{j_1} \subsetneq \dots
\subsetneq M_{j_k},\, M_{j_1},\dots, M_{j_k} \in \chi\}
$$
and
$$
J'_k := J' \cup \{j \in J'';\, d_j \le k\}.
$$
Now we construct a sequence $\kappa^{(m)}$ by induction with respect to
$k$. For this purpose, we introduce
the following condition $\operatorname{Cond}(j,\, k)$ of indices $k$
and $j \in J'_k$.
\begin{enumerate}
\item $\kappa_j^{(m)} \to 0$.
\item $\displaystyle\frac{\kappa^{(m)}_j}{\kappa^{(m)}_{\beta'}} \to 0$
for $\beta' \in J'_k$ with $M_{\beta'} \subsetneq M_j$.
\item
$\displaystyle\frac{\vert q^{(m)} \vert_{\hat{I}_\beta}}{\kappa^{(m)}_j} \to 0$
for $\beta \in \{1,2,\dots,\ell\}$ with $M_j \subsetneq M_\beta$,
and
$\displaystyle\frac{\vert q^{(m)} \vert_{\hat{I}_j}}{\kappa^{(m)}_j} \to 0$
if $j$ belongs to $J''$.
\end{enumerate}

Assume $k = -1$, i.e., $J'_{-1} = J'$.
We set $\kappa^{(m)}_j = \vert q^{(m)} \vert_{\hat{I}_j}$
for $j \in J'$.
It is easy to see that $\operatorname{Cond}(j,\, k)$ is satisfied
for any $j \in J'_{-1}$.
Indeed, as $\xi_j \ne 0$ and $q^{(m)} \in \pi_{j,\,q}^{-1}(G^{(m)}_j)$ hold,
there exists a sequence $\epsilon^{(m)} \to 0$
such that $\vert q^{(m)} \vert_{\hat{I}_\beta}
\le \epsilon^{(m)} \vert q^{(m)} \vert_{\hat{I}_j}$
for $\beta$ with $M_j \subsetneq M_\beta$. This implies 3.~of
$\operatorname{Cond}(j,\, k)$. By the same reason, as $\beta' \in J'$,
the condition 2.~of $\operatorname{Cond}(j,\, k)$ also holds.

Assume that we have constructed $\kappa^{(m)}_j$ for $j \in J'_k$ and
$\operatorname{Cond}(j,\, k)$ holds for any $j \in J'_k$.
Let $j \in J''$ with $d_j = k+1$. Then we will determine
$\kappa^{(m)}_j$ so that $\operatorname{Cond}(j,\, k)$ holds.
First note that
$\chi^*=\{M_{\beta'} \in \chi;\, \beta'  \in J'_k,\, M_{\beta'} \subsetneq M_j\}$
is a totally ordered set with respect to ``$\subset$'', in particular,
we have the maximal submanifold $M_{j^*} \in \chi^*$. If we can determine
$\kappa^{(m)}_j$ satisfying
$\displaystyle\frac{\kappa^{(m)}_j}{\kappa^{(m)}_{j^*}} \to 0$, then
$\displaystyle\frac{\kappa^{(m)}_j}{\kappa^{(m)}_{\beta'}} \to 0$ also holds
for any $M_{\beta'} \in \chi^*$ by induction hypothesis.
It follows from induction hypothesis again that
for any $\beta \in \{1,2,\dots,\ell\}$ with $M_{j^*} \subsetneq M_{\beta}$,
we have
$\displaystyle\frac{\vert q^{(m)} \vert_{\hat{I}_\beta}}{\kappa^{(m)}_{j^*}} \to 0$.
Therefore we can find $\kappa^{(m)}_j$ such that
$\displaystyle\frac{\kappa^{(m)}_j}{\kappa^{(m)}_{j^*}} \to 0$ and
$\displaystyle\frac{\vert q^{(m)} \vert_{\hat{I}_\beta}}{\kappa_j^{(m)}} \to 0$
for $\beta \in \{1,2,\dots,\ell\}$ with $M_j \subset M_{\beta}$ (note that
this contains the case $\beta = j)$.

By repeating the same procedure,
we obtain $\kappa^{(m)}_j$ for any $j \in J''$ with $d_j = k+1$.
As $I_j \cap I_{j'} = \emptyset$ for $j$ and $j'$ with $d_j = d_{j'}$,
we can conclude that
$\operatorname{Cond}(j,\, k+1)$ holds for any $j \in J'_{k+1}$.
Hence we have obtained $\kappa_j^{(m)}$ for all $j \in \{1,2,\dots,\ell\}$.

Let us define points in $\widetilde{X}$ by
$$
\overline{p}^{(m)} := \left(\underset{j \in \{1,2,\dots,\ell\}}{\prod}
\tau_j \left( \frac{1}{\kappa_j^{(m)}}\right)\right) p^{(m)}.
$$
Note that $\widetilde{p}(\overline{p}^{(m)}) \in Z$ still holds.
Then the value of $j$-th coordinate $t_j$ of $\overline{p}^{(m)}$ is given by
that of
$$
\underset{ j \in \{1,2\dots,\ell\}}{\prod}
\mu_j\left( \frac{1}{\kappa^{(m)}_j}\right)
\underset{\beta \in J_{\supsetneq M_j}}{\prod} \mu_\beta
\left(\kappa^{(m)}_j\right) p^{(m)},
$$
which is equal to that of
$$
\mu_j\left( \frac{1}{\kappa^{(m)}_j}\right)
\underset{j \in J_{\supsetneq M_\beta}}{\prod} \mu_j
\left( \kappa^{(m)}_\beta\right)p^{(m)} =
\mu_j\left( \frac{1}{\kappa^{(m)}_j}\right)
\underset{\beta \in J_{\subsetneq M_j}}{\prod} \mu_j
\left(\kappa^{(m)}_\beta\right) p^{(m)}.
$$
Note that $J_{\subsetneq M_j}$ consists of at most 1 element.
If  $J_{\subsetneq M_j}$ is empty, then the
$j$-th component is
$
\kappa_j^{(m)}
$
which clearly tends to $0$ when $m \to \infty$.
If  $J_{\subsetneq M_j} = \{\beta\}$ for some $\beta$, then
it is given by
$
\displaystyle\frac{ \kappa^{(m)}_j }{\kappa^{(m)}_\beta}
$
which also tends to $0$ by the construction of $\kappa^{(m)}$.
As a result, we have $\overline{p}^{(m)} \to p$.
This completes the proof.
\qed

Given the family $\chi = \{M_1, \dots ,M_{\ell}\}$ and a sub-family $\chi_k := \{M_{j_1} ,\dots ,M_{j_k}\}$ of $\chi$. We denote by $S_\chi$ the zero section
$\underset{X,1\leq j\leq\ell}{\times}T_{M_j}\iota_\chi(M_j)$ for the family $\chi$.
We also introduce the set
\begin{equation}{\label{def:S-chi-chi}}
S_{\chi/\chi_k} :=
\left(\underset{X,1\leq \alpha \leq k}{\times}T_{M_{j_\alpha}}\iota_\chi(M_{j_\alpha})
\right)
\underset{X}{\times} M,
\end{equation}
where $M=\bigcap_{j=1}^\ell M_j$.
We emphasize that, in the above definition, we use $\iota_\chi$ (not $\iota_{\chi_k}$).
Then we have the natural embedding
$$
S_\chi \hookleftarrow S_{\chi/\chi_k} \hookrightarrow S_{\chi_k}.
$$




%

\begin{cor} \label{cor: restriction multi-cones}
Let $k \le \ell$ and $\{j_1,\dots,j_k\}$ be a subset of $\{1,2,\dots,\ell\}$. Set $\chi_k=\{M_{j_1},\ldots,M_{j_k}\}$.
Let $Z$ be a subset of $X$. Then we have
$$
C_\chi (Z) \cap S_{\chi/\chi_k}
= C_{\chi_k} (Z) \cap S_{\chi/\chi_k}.
$$
\end{cor}
\dim\ \ Let us prove $\subseteq$. Suppose that $p \in
S_{\chi/\chi_k}$
does not belong to $C_{\chi_k}(Z)$.
By Lemma \ref{lem:multi cone} there exists an open subset $\psi(q) \in U \subset (TM)_q$
and a multi-cone $G' \in \operatorname{Cone}_{\chi_k}(p)$ such that
$
\psi(Z) \cap G' \cap U = \emptyset
$
holds. For $\alpha \notin \{j_1,\ldots,j_k\}$, set
$G_\alpha =(T_{M_\alpha} X)_q$. Hence $G:=\left(\underset{\alpha \notin \{j_1,\dots,j_k\}}{\bigcap} \imin {\pi_\alpha} (G_\alpha)\right) \cap G' = G' \in \operatorname{Cone}_{\chi}(p)$ and $\psi(Z) \cap G \cap U = \emptyset$. By Lemma \ref{lem:multi cone} we obtain $p \notin C_\chi(Z)$.


\

Let us prove $\supseteq$. Suppose that
$p \in S_{\chi/\chi_k}$
does not belong to $C_\chi(Z)$. By Lemma \ref{lem:multi cone} there exists an open subset $\psi(q) \in U \subset (TM)_q$
and a multi-cone $G \in \operatorname{Cone}_{\chi}(p)$ such that
$
\psi(Z) \cap G \cap U = \emptyset
$
holds. For $\alpha \notin \{1,\ldots,k\}$, $p_{j_\alpha}=(q;0)$ and we have
$\operatorname{Cone}_{\chi, j_\alpha}(p) = \{(T_{M_{j_\alpha}} X)_q\}$. Hence 
$G':=G$ can be regarded as an element of $\operatorname{Cone}_{\chi_k}(p)$ and $\psi(Z) \cap G' \cap U = \emptyset$. By Lemma \ref{lem:multi cone} we obtain $p \notin C_{\chi_k}(Z)$.
\qed









\begin{df} Denote by $\op(\widetilde{X})$ the category of open subsets of $\widetilde{X}$, and let $Z$ be a subset of $\widetilde{X}$.

(i) 
We set $\RP_j Z=\mu_j(Z,\RP).$
If $U \in \op(\widetilde{X})$, then $\RP_j U \in \op(\widetilde{X})$ since $\mu_j$ is open for each $j=1,\ldots,\ell$.

(ii)  
Let $J = \{j_1,\dots,j_k\} \subset \{1,\ldots,\ell\}$. We set $$\RP_JZ=\RP_{j_1}\cdots\RP_{j_k}Z=\mu_{j_1}(\cdots\mu_{j_k}(Z,\RP),\dots,\RP).$$
We set $(\RP)^\ell Z=\RP_{\{1,\ldots,\ell\}}Z=\mu(Z,(\RP)^\ell).$
If $U \in \op(\widetilde{X})$, then $\RP_J U \in \op(\widetilde{X})$ since $\mu_j$ is open for each $j \in \{1,\ldots,\ell\}$.

(iii) 
We say that $Z$ is $(\RP)^\ell$-conic ($\ell$-conic for short) if
$Z=(\RP)^\ell Z$. In other words, $Z$ is invariant by the action of
$\mu_j$, $j=1,\dots,\ell$.

\end{df}

\begin{df}
(i)  We say that a subset $Z$ of $\widetilde{X}$ is $\RP_j$-connected if
$Z \cap \RP_j x$ is connected for each $x \in Z$.

(ii) We say that a subset $Z$ of $\widetilde{X}$ is $(\RP)^\ell$-connected if
there exists a permutation $\sigma:\{1,\ldots,\ell\}\to\{1,\ldots,\ell\}$ such that
\begin{equation*}
 \begin{cases}
 \text{$Z$ is $\RP_{\sigma(1)}$-connected,}\\
 \text{$\RP_{\sigma(1)}Z$ is $\RP_{\sigma(2)}$-connected,}\\
 \text{\ \ $\vdots$}\\
 \text{$\RP_{\sigma(1)}\cdots\RP_{\sigma(\ell-1)}Z$ is $\RP_{\sigma(\ell)}$-connected}.
  \end{cases}
\end{equation*}
\end{df}


The proof of the following is almost the same as that of Proposition 4.1.3 of \cite{KS90}, and we shall not repeat it.

\begin{prop}\label{m4.1.3} Let $V$ be a $(\RP)^\ell$-conic open subset of the zero section $S$.
\begin{itemize}
\item[(i)] Let $W$ be an open neighborhood of $V$ in $\widetilde{X}$, and let $U=\widetilde{p}(W \cap \Omega)$. Then $V \cap C_\chi(X \setminus U)=\emptyset$.
\item[(ii)] Let $U$ be an open subset of $X$ such that $V \cap C_\chi(X \setminus U)=\emptyset$. Then $\imin {\widetilde{p}}(U) \cup V$ is an open neighborhood of $V$ in $\overline{\Omega}$.
\end{itemize}
\end{prop}


Set, for $j=1,\dots,\ell$,
$$
S_j :=
\{(q;\,\xi_1,\dots,\xi_\ell) \in S = \underset{X,1\leq \alpha \leq \ell}{\times}
T_{M_{\alpha}}\iota_\chi(M_{\alpha});\, \xi_j = 0\} \subset S
$$
and
$
\pi_j: S \to S_j
$
be the canonical projection induced from one $T_{M_j}\iota_\chi(M_j) \to T_{M_j}M_j = M_j$.
Let $V$ be an $(\RP)^\ell$-conic
subanalytic subset in $S$.
We introduce the following conditions Va.~and Vb.~of $V$ for each $j$.
\begin{enumerate}
\item[{\bf{Va.}}] $V$ does not intersect $S_j$.
\item[{\bf{Vb.}}] $\pi_j(V) \subset \pi_j(V \cap S_j)$.
\end{enumerate}

Then we have the following proposition whose proof will be given in Appendix D.
\begin{prop}\label{m4.1.4}
Let $V$ be an $(\RP)^\ell$-conic subanalytic subset of $S$ which
satisfies the condition either Va.~or Vb.~for each $j$.
Assume that $X=\mathbb{R}^n$ and each $M_j$ is defined
by linear equations as described in 2.~of Proposition \ref{prop:canonical-coordinate}.
Then any subanalytic neighborhood $W$
of $V$ in $\widetilde{X}$ contains $\widetilde{W}$ open and subanalytic in
$\widetilde{X}$ such that:
\begin{equation} \label{mpitilde}
\begin{cases}
\text{(i) $\widetilde{W} \cap \Omega$ is $(\RP)^\ell$-connected,} \\
\text{(ii) $\RP_1 \cdots \RP_{\ell}(\widetilde{W} \cap \Omega)=\imin {\widetilde{p}}(\widetilde{p}(W \cap \Omega))$ is subanalytic in $\widetilde{X}$.}
\end{cases}
\end{equation}
\end{prop}

Let $Z$ be a closed sub-bundle of the vector bundle $S$ over
a closed submanifold of $M$ which is $(\RP)^\ell$-conic.
We have the following corollary.
\begin{cor}\label{cor: Rpl neighborhood}
Let $V$ be a $(\RP)^\ell$-conic open subanalytic subset in $Z$.
Then there exists a locally finite family $\{V^{(\alpha)}\}_\alpha$
of $(\RP)^\ell$-conic open subanalytic subsets in $Z$
which satisfies the following conditions.
\begin{enumerate}
\item $V = \cup_{\alpha} V^{(\alpha)}$ and  each $V^{(\alpha)}$ satisfies the condition either Va.~or Vb.~for each $j$.
\item For any subanalytic open neighborhood $W$ of $V$ in $\widetilde{X}$,
	there exists a subanalytic open neighborhood $\widetilde{W}^{(\alpha)} \subset W$ of $V^{(\alpha)}$ for which
	(\ref{mpitilde}) of Proposition \ref{m4.1.4} holds.
	Furthermore $\{\widetilde{p}(\widetilde{W}^{(\alpha)} \cap \Omega)\}_\alpha$
	is a locally finite family of subanalytic open subsets of $X$.
\end{enumerate}
\end{cor}
\begin{oss}
The following claim also holds.
For any subanalytic open neighborhood $W$ of
a finite intersection $V'$ of $V^{(\alpha)}$'s in $\widetilde{X}$,
there exists a subanalytic open neighborhood $\widetilde{W} \subset W$ of $V'$
for which (\ref{mpitilde}) holds.
\end{oss}
\dim \ It is easy to find a finite family $\{V^{(k)}\}$ of subanalytic open subsets in $Z$ such that $V = \cup_k V^{(k)}$ and
each $V^{(k)}$ satisfies the condition either Va.~or Vb.~for each $j$.
Therefore we may assume that $V$ satisfies these from the beginning.

Let $\{X'^{(\beta)}\}_{\beta \in \Lambda'}$ be a locally finite
family of open subanalytic subsets in $X$ which is a finer covering
$\psi: \Lambda' \to \Lambda$ of $\{X^{(\alpha)}\}_{\alpha \in \Lambda}$
(see Remark {\ref{rem:coordinate-systems}} for $X^{(\alpha)}$)
and such that each $X'^{(\beta)}$ is relatively compact in $X^{ (\psi(\beta))}$.
Set
$$
V^{(\beta)} := V \cap p^{-1}(X'^{(\beta)}),\,
W^{(\beta)} := W \cap p^{-1}(X'^{(\beta)}).
$$
Then, in particular, $W^{(\beta)}$ is a subanalytic open neighborhood of $V^{(\beta)}$ in the local model $\widetilde{X}^{(\psi(\beta))}$.
It follows from Proposition \ref{m4.1.4} that there exists an open neighborhood $\widetilde{W}^{(\beta)} \subset W^{(\beta)}$
satisfying (\ref{mpitilde}) in $\widetilde{X}^{ (\psi(\beta))}$.
Since the multi-actions are stable on the local model, and since $p(\widetilde{W}^{(\beta)}) \subset X'^{(\beta)}$ and
$$
p(\RP_1 \cdots \RP_{\ell}(\widetilde{W}^{(\beta)} \cap \Omega)) =
p(\widetilde{W}^{(\beta)} \cap \Omega) \subset X'^{(\beta)}
$$
hold, both subsets $\widetilde{W}^{(\beta)}$ and $\RP_1 \cdots \RP_{\ell}(\widetilde{W}^{(\beta)} \cap \Omega)$ become subanalytic in $\widetilde{X}$
by Lemma \ref{lem:subanalytic-local-model}.
By stability of the multi-actions on the local model, $\widetilde{W}^{(\beta)}$ is $(\RP)^{\ell}$-connected in $\widetilde{X}$.
The family $\{\widetilde{p}(\widetilde{W}^{(\beta)} \cap \Omega)\}_\beta$ is clearly locally finite as $\{X'^{(\beta)}\}$ is locally finite.
This completes the proof.
\qed

\begin{cor}{\label{cor:intersection-V-U}}
For a $(\RP)^\ell$-conic subanalytic open subset $V$, let
$\mathcal{T}(V)$ denote the family of subanalytic open subsets $U$ in $X$ with
$$
C_\chi(X\setminus U) \cap V = \emptyset.
$$
Let $V_1$ and $V_2$ be $(\RP)^\ell$-conic subanalytic open subsets
in $S$.
Then the family $\{U_1 \cup U_2\}$ $($resp. $\{U_1 \cap U_2\}$$)$
$($$U_k \in \mathcal{T}_k,\,k=1,2$$)$
is equivalent to $\mathcal{T}(V_1 \cup V_2)$
$($resp. $\mathcal{T}(V_1 \cap V_2)$$)$.
Here the equivalence is given by the inclusion of sets.
\end{cor}
\dim \
The assertion for $V_1 \cup V_2$ is clear.
Let us show the assertion for $V_1 \cap V_2$.
It suffices to show that,
for $U \in \mathcal{T}(V_1 \cap V_2)$,
we find $U_1$ and $U_2$ $(U_k \in \mathcal{T}(V_k))$ with $U_1 \cap U_2 \subset U$.

%

Set
$$
\widetilde{U} := \widetilde{p}^{-1}(U) \cup (V_1 \cap V_2) \cup
(\widetilde{X} \setminus \overline{\Omega}).
$$
Then $\widetilde{U}$ is a subanalytic open neighborhood of $V_1 \cap V_2$ by
(ii) of Proposition \ref {m4.1.3}.
It is easy to construct, in $\widetilde{X}$,
a subanalytic open subset $W_k$ ($k=1,2$) such that
$W_k$ is an open neighborhood of $V_k$ and
$W_1 \cap W_2 \subset \widetilde{U}$.
Let $\{V_k^{(\alpha)}\}_\alpha$ and $\{W_k^{(\alpha)}\}_\alpha$ be the families given by
Corollary \ref{cor: Rpl neighborhood} for $V_k \subset W_k$.
Set
$$
\begin{aligned}
U^{(\alpha)}_k &:= \widetilde{p}\left(\widetilde{W}^{(\alpha)}_k \cap \Omega\right)
= \widetilde{p}\left( (\RP)^\ell(\widetilde{W}^{(\alpha)}_k \cap \Omega)\right) \\
&= \widetilde{p}\left( ((\RP)^\ell(\widetilde{W}^{(\alpha)}_k \cap \Omega)) \cap \{t_1= \dots= t_\ell=1\}\right)
\end{aligned}
$$
which is subanalytic in $X$ thanks to (\ref{mpitilde}) of Proposition \ref{m4.1.4}.
As $\{U^{(\alpha)}_k\}$ is locally finite in $X$,
$$
U_k := \underset{\alpha}{\cup} U^{(\alpha)}_k, \qquad (k=1,2)
$$
is also subanalytic in $X$.
It follows from (i) of Proposition \ref {m4.1.3},
each $U_k^{(\alpha)}$ belongs to $\mathcal{T}(V_k^{(\alpha)})$, and hence,
we have $U_k \in \mathcal{T}(\cup_\alpha V_k^{(\alpha)}) = \mathcal{T}(V_k)$
for $k=1,2$.
Finally, as $\widetilde{p}|_{\{t=t^*\}}$
gives an isomorphism for $t^*_j > 0$ ($j=1,\dots,\ell$),
we have
$$
\begin{aligned}
 U_1 \cap U_2
 &= \left(\underset{\alpha}{\cup}\,
 \widetilde{p}\left(\widetilde{W}^{(\alpha)}_1 \cap \Omega\right)\right) \cap
 \left(\underset{\alpha}{\cup}\, \widetilde{p}
 \left(\widetilde{W}^{(\alpha)}_2 \cap \Omega\right)\right) \\
 &\subset \widetilde{p}(W_1 \cap \Omega) \cap
        \widetilde{p}(W_2 \cap \Omega)
 = \widetilde{p}( (W_1 \cap W_2) \cap \Omega)
\subset \widetilde{p}(\widetilde{U} \cap \Omega) = U.
\end{aligned}
$$
This completes the proof.
\qed

\end{section}

\begin{section}{Morphisms between multi-normal deformations}\label{5}

Let $X$ and $Y$ be real analytic manifolds of dimension
$n$ and $m$ respectively
and let
$\chi^M=\{M_j\}_{j=1}^\ell$, $\chi^N=\{N_j\}_{j=1}^\ell$ be
families of smooth closed submanifolds of $X$ and $Y$ 
respectively which satisfy H1, H2 and H3.
Let $f:X \to Y$ be a morphism of real analytic manifolds such that
$f(M_j) \subseteq N_j$, $j=1,\ldots,\ell$.

We want to extend $f$ to a morphism $\widetilde{f}:\widetilde{X} \to \widetilde{Y}$. This is done by repeatedly employing the usual construction
of a morphism between normal deformations,
i.e. we extend $f$ to  $\widetilde{f}_1:\widetilde{X}_{M_1} \to \widetilde{Y}_{N_1}$,
then we extend $\widetilde{f}_1$ to $\widetilde{f}_{1,2}:\widetilde{X}_{M_1,M_2} \to \widetilde{Y}_{N_1,N_2}$.
Then we can define recursively $\widetilde{f}:\widetilde{X} \to \widetilde{Y}$
by extending the morphism $\widetilde{f}_{1,\ldots,\ell-1}:\widetilde{X}_{M_1,\ldots,M_{\ell-1}} \to \widetilde{Y}_{N_1,\ldots,N_{\ell-1}}$
to the normal deformations with respect to $M_\ell$ and $N_\ell$ respectively. We also denote by $S^M$ (resp. $S^N$) the zero section of
$\widetilde{X}$ (resp. $\widetilde{Y}$).

In a local coordinate system set
$$
f(x_1,\ldots,x_n)=(f_1(x_1,\ldots,x_n),\ldots,f_m(x_1,\ldots,x_n)).
$$
Let $I^M_{j}$ (resp. $I^N_{j}$) ($j=1,\dots,\ell$) be a subset of $\{1,\dots,n\}$ (resp. $\{1,\dots,m\}$) satisfying
the conditions (\ref{eq:conditions-indices-set}).
The $J^M_{i}$ (resp. $J^N_i$) was defined in \eqref{eq:def-J_i} for $i \in \{1,\dots,n\}$
(resp. $\{1,\dots,m\}$).
Then, outside $\{t_1t_2\ldots t_\ell = 0\}$ in $\widetilde{X}$,
the morphism $\widetilde{f}:\widetilde{X}\to\widetilde{Y}$ $$
(x_1,\dots,x_n,t_1,\dots,t_\ell) \to
(y_1,\dots,y_m,\hat{t}_1,\dots,\hat{t}_\ell)
$$
is explicitly written in the form
$$
\left\{
\begin{aligned}
&\hat{t}_{J^N_k}y_k=f_k(t_{J^M_1}x_1,\ldots,t_{J^M_n}x_n)
\qquad (k = 1,2,\dots,m),\\
&\hat{t}_j = t_j \qquad (j = 1,2,\dots,\ell),
\end{aligned}
\right.
$$
which naturally extends to the whole $\widetilde{X}$.

In order to understand the restriction of $\widetilde{f}$ to $S^M \subset \widetilde{X}$,
we follow the notations of the proof of Proposition \ref{prop:zero section}. For any $k=1,2,\dots,m$, we set
$$
I^N(k) := \underset{k \in I^N_j}{\bigcap} I^N_j = \underset{j \in J^N_k}
{\bigcap} I^N_j.
$$
Let $k \in \underset{1 \le j \le \ell}{\bigcup} I^N_j$, and let
$j(k) \in \{1,2,\dots,\ell\}$ such that $I^N(k) = I^N_{j(k)}$.
Then, for any $j \in J^N_k$, we get $f(M_j) \subseteq N_{j(k)}$
because of $f(M_j) \subset N_j \subset N_{j(k)}$. We set
$$
I=\bigcup_{j \in J^N_k} I^M_j,\,\, M_I=\bigcap_{j \in J^N_k} M_j.
$$
By expanding $f_k(x_1,\dots,x_n)$ along the submanifold $M_I$,
we obtain
$$
f_k(x_1,\dots,x_n) =
\sum_{i \in I}\left.\frac{\partial f_k}{\partial x_i}\right|_{M_I} x_i
+ \frac{1}{2}\sum_{i_1,i_2 \in I}
\left.
\frac{\partial^2 f_k}{\partial x_{i_1}\partial x_{i_2}}\right|_{M_I}
x_{i_1}x_{i_2}
+ \dots,
$$
as $f_k\vert_{M_I} = 0$ holds due to $f(M_I) \subset N_{j(k)}$.
Then we get, on $t_1\dots t_\ell \ne 0$,
\begin{equation}{\label{eq:morphism-expansions}}
y_k =
\sum_{i \in I}\left.\frac{\partial f_k}{\partial x_i}\right|_{M_I}
\frac{t_{J^M_i}}{t_{J^N_k}}x_i
+ \frac{1}{2}\sum_{i_1,i_2 \in I}
\left.
\frac{\partial^2 f_k}{\partial x_{i_1}\partial x_{i_2}}\right|_{M_I}
\frac{t_{J^M_{i_1}}t_{J^M_{i_2}}}{t_{J^N_k}}x_{i_1}x_{i_2}
+ \dots.
\end{equation}
Let $i_1, \dots, i_p$ be a sequence $(i_1, \dots, i_p \in I)$.
\begin{itemize}
\item[(i)] If there exists $j \in J^N_k$ such that
$\{i_1,\ldots,i_p\} \cap I^M_j=\emptyset$, as
$f_k \vert_{M_j} = 0$ and derivatives
$\displaystyle\frac{\partial}{\partial x_{i_1}}$,
$\ldots$,
$\displaystyle\frac{\partial}{\partial x_{i_p}}$
are tangent to $M_j$, we have
$$
\left.{\partial^p f_k \over \partial x_{i_1} \ldots \partial x_{i_p}}\right|_{M_I}=0.
$$
\item[(ii)] Suppose that $\{i_1, \ldots, i_p\} \bigcap I^M_j \ne \emptyset$
holds for each $j \in J^N_k$. This implies that, for any $j \in J^N_k$, the
$I^M_j$ contains some $i_q \in \{i_1,\ldots, i_p\}$
$\iff$
$j \in J^M_{i_q}$ with some $i_q \in \{i_1,\ldots,i_p\}$.
Hence we obtain
$$
J^N_k \subset J^M_{i_1} \cup \cdots \cup J^M_{i_p}.
$$
Now we have two cases.
    \begin{itemize}
    \item[(a)] If some pair of $J^M_{i_1}$, $\dots$, $J^M_{i_p}$
is not disjoint or if $J^N_k \subsetneqq J^M_{i_1} \sqcup \cdots \sqcup J^M_{i_p}$,
then $t_{J^M_{i_1}} \cdots t_{J^M_{i_p}} \imin {t_{J^N_k}} \to 0$
when $t \to 0$.
    \item[(b)] If $J^N_k = J^M_{i_1} \sqcup \cdots \sqcup J^M_{i_p}$ holds,
	    then the term with its indices $i_1,\dots,i_p$
	    in (\ref{eq:morphism-expansions})
	    becomes, by letting $t$ to $0$,
         $$
 \left.\frac{1}{p!}{\partial^p f_k \over \partial x_{i_1} \cdots \partial x_{i_p}}\right|_{M_I}x_{i_1} \cdots x_{i_p}.
         $$
     \end{itemize}
\end{itemize}
From these observations, the morphism $\widetilde{f}$ is
described by, on $S \subset \widetilde{X}$,
         $$
	 y_k = \sum_{J^M_{i_1} \sqcup \cdots \sqcup J^M_{i_p} = J^N_k} \left.\frac{1}{p!}{\partial^p f_k \over \partial x_{i_1} \cdots \partial x_{i_p}}\right|_{M}x_{i_1} \cdots x_{i_p}
 \qquad (k\in \underset{1 \le j \le \ell}{\bigcup} I^N_j).
         $$
Here $M:=M_1 \cap \dots \cap M_\ell$.  Note that if there is no $\{i_1, \dots, i_p\}$ with
$J^N_k = J^M_{i_1} \sqcup \cdots \sqcup J^M_{i_p}$, then we set $y_k: = 0$.

Let us study the condition $J^N_k = J^M_{i_1} \sqcup \cdots \sqcup J^M_{i_p}$.
For this purpose, we introduce two definitions.
Let $k \in \underset{1 \le j \le \ell}{\bigcup} I^N_j$ and set
$$
\sup{}_\subset J^N_k := \left\{j \in J^N_k;\,
\text{$M_j$ is a maximal submanifold in
$\{M_\beta\}_{\beta \in J^N_k}$}\right\}
$$
and
$$
\inf{}_\subset J^N_k := \left\{j \in J^N_k;\,
\text{$M_j$ is a minimal submanifold
in $\{M_\beta\}_{\beta \in J^N_k}$}\right\}.
$$

Note that $\#\inf_\subset J^N_k \le \#\sup_\subset J^N_k$ holds.

\begin{es}
Let us consider closed submanifolds $\{M_1, M_2, M_3\}$
in $X=\R^n$ and $\{N_1, N_2, N_3\}$ in $Y=\R^m$. Let $f: X \to Y$ be
a morphism satisfying $f(M_j) \subset N_j$ ($j=1,2,3$). We assume that
$N_3 \subset N_2 \subset N_1$, $M_3 \subset M_1$, $M_3 \subset M_2$,
$M_1$ and $M_2$ intersect transversely. For $k \in I^N_1$,
$J^N_{k} = \{1,2,3\}$,
$\sup{}_\subset J^N_k = \{1,2\}$ and
$\inf{}_\subset J^N_k = \{3\}$. For $k \in \hat{I}^N_2=I^N_2 \setminus I^N_1$,
$J^N_k=\{2,3\}$, $\sup{}_\subset J^N_k = \{2\}$ and
$\inf{}_\subset J^N_k = \{3\}$. For $k \in \hat{I}^N_3=I^N_3 \setminus (I^N_1 \cup I^N_2)$,
$J^N_k=\{3\}$, $\sup{}_\subset J^N_k = \inf{}_\subset J^N_k = \{3\}$.
\end{es}

\begin{lem}{\label{lemma:condition-disjoint}}
Let $k \in \underset{1 \le j \le \ell}{\bigcup} I^N_j$ and $\{i_1, \dots,i_p\}$ a subset of
$\underset{1 \le j \le \ell}{\bigcup}I^M_j$.
Then
$J^N_k = J^M_{i_1} \sqcup \cdots \sqcup J^M_{i_p}$ holds if and only if
the conditions (a) and (b) below are satisfied.
\begin{enumerate}
\item[(a)] $k$ satisfies the following condition $(\dagger)_k$
\begin{eqnarray}
&&\#\inf{}_\subset J^N_k = \#\sup{}_\subset J^N_k, \label{eq:condition-disjoint-1} \\
 \notag \\
&&M_\beta \subset M_j\,\,
(j \in J^N_k,\,\beta \in \{1,2,\dots,\ell\})
\Longrightarrow \beta \in J^N_k. \label{eq:condition-disjoint-2}
\end{eqnarray}
\item[(b)] $p = \#\sup{}_\subset J^N_k$ and the indices $i_1,\dots, i_p$ satisfy
\begin{equation}\label{eq:condition-disjoint-3}
i_\alpha \in \hat{I}^M_{\sigma(\alpha)} = I^M_{\sigma(\alpha)} \setminus \left(\bigcup_{I^M_j \subsetneqq I^M_{\sigma(\alpha)}}I^M_j\right) \qquad (\alpha \in \{1,2,\dots,p\})
\end{equation}
for some bijection $\sigma: \{1,2\dots,p\} \to \sup_\subset J^N_k$.
\end{enumerate}
\end{lem}
\dim \ \
We first show the claim that
if $\#\inf_\subset J^N_k < \#\sup_\subset J^N_k$,
then there exists no $\{i_1, \dots, i_p\}$ with
$J^N_k = J^M_{i_1} \sqcup \cdots \sqcup J^M_{i_p}$.
Assume that there exists $\{i_1, \dots, i_p\}$ with
$J^N_k = J^M_{i_1} \sqcup \cdots \sqcup J^M_{i_p}$.
Then it follows from  $\#\inf_\subset J^N_k < \#\sup_\subset J^N_k$ that we can
find indices $j$, $j'$, $j''$ in $J^N_k$ satisfying $I^M_{j'} \subset I^M_{j}$,
$I^M_{j''} \subset I^M_{j}$ and $I^M_{j'} \cap I^M_{j''} = \emptyset$.
By the assumption, there exist $\alpha'$ and $\alpha''$ in $\{1,2,\dots,p\}$
such that $j' \in J^M_{i_{\alpha'}}$ and $j'' \in J^M_{i_{\alpha''}}$.
As $I^M_{j'} \cap I^M_{j''} = \emptyset$, we have $\alpha' \ne \alpha''$.
On the other hand,  $I^M_{j'} \subset I^M_{j}$ and $I^M_{j''} \subset I^M_{j}$
implies $j \in J^M_{i_{\alpha'}}$ and $j \in J^M_{i_{\alpha''}}$, which contradicts
that $J^M_{i_{\alpha'}}$ and $J^M_{i_{\alpha''}}$ are disjoint. Hence we have obtained
the claim.

\medskip

In what follows, we assume
$$
\#\inf{}_\subset J^N_k = \#\sup{}_\subset J^N_k.
$$
Suppose that $\{i_1, \dots, i_p\}$ satisfies the condition
$J^N_k = J^M_{i_1} \sqcup \cdots \sqcup J^M_{i_p}$.

Let $j \in J^N_k$ and $\beta \in \{1,2,\dots,\ell\}$ with
$M_\beta \subset M_j$. Then we can find $i_\alpha$ such that
$j \in J^M_{i_\alpha}$. $M_\beta \subset M_j$ implies
$\beta \in J^M_{i_\alpha} \subset J^N_k$. Therefore we have (\ref{eq:condition-disjoint-2}).

Let $j \in \sup{}_\subset J^N_k$. Then some $J^M_{i_q}$ ($q \in \{1,2,\dots,p\}$)
contains $j$, which implies $i_q \in I^M_j$.
Further we can show that
$i_q$ belongs to $\hat{I}^M_j$. If $i_q \in I^M_j \setminus \hat{I}^M_j$,
then there exists $j'$ with $i_q \in I^M_{j'}$ and
$M_j \subsetneqq M_{j'}$ by the definition of $\hat{I}^M_j$,
from which we have $j' \in J^M_{i_q} \subset J^N_k$. This contradicts
$j \in \sup{}_{\subset} J^N_k$. Hence $i_q \in \hat{I}^M_j$ and
we have obtained that the set $\{i_1, \dots, i_p\}$
contains at least one index
that belongs to $\hat{I}^M_j$ for any $j \in \sup{}_{\subset} J^N_k$.
Further two or more indices in $\hat{I}^M_j$ cannot belong to
$\{i_1, \dots, i_p\}$ at the same time because any pair of
$J^M_{i_1}$, $\cdots$, $J^M_{i_p}$ is disjoint.

Now we show that $\{i_1, \dots, i_p\}$ consists of only these indices.
Let $i$ be an element in $\{i_1, \dots, i_p\}$.
Choosing $j' \in J^M_i$, we can find $j \in \sup{}_{\subset} J^N_k$
with $M_{j'} \subset M_j$. Then, by the above argument,
there exists $i_q \in \hat{I}^M_j$ which belongs to
$\{i_1, \dots, i_p\}$.
As $i_q \in I^M_j \subset I^M_{j'}$,  we have $j' \in J^M_{i_q}$,
which implies $J^M_{i_q} \cap J^M_{i} \ne \emptyset$.
Since each pair is disjoint, we have $i = i_q$.

Therefore we can find a bijection $\sigma: \{1,2,\dots,p\} \to \sup_{\subset} J^N_k$
such that $i_\alpha \in \hat{I}^M_{\sigma(\alpha)}$ and
$J^M_{i_\alpha} \subset J^N_k$($\alpha \in \{1,2,\dots,p\}$).

Conversely if such a $\sigma$ exists, then
$J^N_k = J^M_{i_1} \sqcup \cdots \sqcup J^M_{i_p}$ easily follows
from $\#\inf{}_\subset J^N_k = \#\sup{}_\subset J^N_k$ and
$J^M_{i_\alpha} \subset J^N_k$. The last inclusion can be obtained by the following argument.
Let $\beta \in J^M_{i_\alpha}$. Then, as $i_{\alpha} \in
\hat{I}_{\sigma(\alpha)}$
and $i_{\alpha} \in I_\beta$ hold, we have $I_{\sigma(\alpha)} \subset
I_\beta$,
from which we have $M_\beta \subset M_{\sigma(\alpha)}$ (the inclusion $\supset$ cannot hold because of the maximality of $M_{\sigma(\alpha)}$).
Hence we obtain $\beta \in J^N_k$ by the condition $(\dagger)_k$.
\qed

\begin{es} Let us consider closed submanifolds $\{M_1, M_2, M_3\}$
in $X=\R^3$ and $\{N_1, N_2, N_3\}$ in $Y=\R^m$. Let $f: X \to Y$ be
a morphism satisfying $f(M_j) \subset N_j$ ($j=1,2,3$). We consider these three cases:
\begin{enumerate}
\item $M_i=\{x_i=0\}$, $i=1,2,3$ (Majima),
\item $M_1=\{0\}$, $M_2=\{x_1=x_2=0\}$, $M_3=\{x_1=0\}$ (Takeuchi),
\item $M_1=\{x_1=0\}$, $M_2=\{x_2=0\}$, $M_3=\{0\}$ (Mixed).
\end{enumerate}
Let us see some examples of  conditions (a), (b) of Lemma \ref{lemma:condition-disjoint}.
\begin{itemize}
\item[(i)] Suppose that $J_k^N=\{1,2,3\}$. \\ In 1. $\sup{}_\subset J^N_k = \inf{}_\subset J^N_k = \{1,2,3\}$ and condition $(\dagger)_k$ is clearly satisfied. Moreover $i_1=\sigma(1),i_2=\sigma(2),i_3=\sigma(3)$ satisfy \eqref{eq:condition-disjoint-3} for any permutation $\sigma$ of $\{1,2,3\}$. \\
    In 2. $\sup{}_\subset J^N_k = \{3\}$, $\inf{}_\subset J^N_k = \{1\}$ and condition $(\dagger)_k$ is clearly satisfied. Moreover $i_1=1$ satisfies \eqref{eq:condition-disjoint-3}. \\
     In 3. $\sup{}_\subset J^N_k = \{1,2\}$ and
$\inf{}_\subset J^N_k = \{3\}$. Then condition \eqref{eq:condition-disjoint-1} is not satisfied.
\item[(ii)] Suppose that $J_k^N=\{2,3\}$. \\
    In 1. $\sup{}_\subset J^N_k = \inf{}_\subset J^N_k = \{2,3\}$ and $(\dagger)_k$ is satisfied. It is easy to check that the indices $i_1=2$, $i_2=3$ (or $i_1=3$, $i_2=2$) satisfy \eqref{eq:condition-disjoint-3}. \\
    In 2. $\sup{}_\subset J^N_k = \{3\}$ and $\inf{}_\subset J^N_k = \{2\}$. Then condition \eqref{eq:condition-disjoint-1} is satisfied but $M_2,M_3 \supset M_1$ and $1 \notin \sup{}_\subset J^N_k$, hence \eqref{eq:condition-disjoint-2} does not hold. \\
    In 3. $\sup{}_\subset J^N_k = \{2\}$, $\inf{}_\subset J^N_k = \{3\}$ and $(\dagger)_k$ is satisfied. Moreover the index $i_1=2$ satisfies \eqref{eq:condition-disjoint-3}.
\item[(iii)] Suppose that $J_k^N=\{1,2\}$. \\
    In 1. $\sup{}_\subset J^N_k = \inf{}_\subset J^N_k = \{1,2\}$ and $(\dagger)_k$ is satisfied. It is easy to check that the indices $i_1=1$, $i_2=2$ (or $i_1=2$, $i_2=1$) satisfy \eqref{eq:condition-disjoint-3}. \\
    In 2. $\sup{}_\subset J^N_k = \{2\}$ and $\inf{}_\subset J^N_k = \{1\}$. In this case condition \eqref{eq:condition-disjoint-1} is satisfied and \eqref{eq:condition-disjoint-2} holds too. Moreover the index $i_1=2$ satisfies \eqref{eq:condition-disjoint-3}.\\
    In 3. $\sup{}_\subset J^N_k = \inf{}_\subset J^N_k = \{1,2\}$. Then condition \eqref{eq:condition-disjoint-1} is satisfied but $M_1,M_2 \supset M_3$ and $3 \notin \sup{}_\subset J^N_k$, hence \eqref{eq:condition-disjoint-2} does not hold.
\end{itemize}
\end{es}

As an immediate consequence of Lemma \ref{lemma:condition-disjoint}, we have the following.
\begin{cor}{\label{cor:multilinear-map}}
For $k \in \underset{1 \le j \le \ell}{\bigcup}I^N_j$,
$y_k$ of the morphism $\widetilde{f}$ on $S^M \subset \widetilde{X}$ is given by
\begin{equation}{\label{eq:multi-morphism-zero}}
y_k=\left\{\begin{array}{ll}
 \displaystyle\sum_{i_1 \in \hat{I}^M_{j_1},\, \dots,\, i_p \in \hat{I}^M_{j_p}
}
\left.{\displaystyle\partial^p f_k \over \partial x_{i_1} \cdots
\partial x_{i_p}}\right|_{M}x_{i_1} \cdots x_{i_p}\qquad
&\text{$($$(\dagger)_k$ holds$)$},\\
\\
 0 \qquad &\text{$($otherwise$)$}.
 \end{array}\right.
\end{equation}
Here $p = \#\sup{}_\subset J^N_k$, $\{j_1, \dots, j_p\} = \sup{}_\subset J^N_k$
and the condition $(\dagger)_k$ for $k$ was given in Lemma
\ref{lemma:condition-disjoint}.
\end{cor}

Let $j \in \{1,2,\dots,l\}$.
Then, for any $k$ and $k'$ in $\hat{I}^N_j$, as
$J^N_k = J^N_{k'}$ holds, we have
$\inf{}_\subset J^N_k = \inf{}_\subset J^N_{k'}$ and
$\sup{}_\subset J^N_k = \sup{}_\subset J^N_{k'}$.
This implies that the sets $\inf{}_\subset J^N_k$ and
$\sup{}_\subset J^N_k$ do not depend on a choice of
$k \in \hat{I}^N_j$.
We denote them
by $\underline{\mathcal J}(j)$ and
$\overline{\mathcal J}(j)$ respectively for short. As the submanifolds in $\chi^M$ and $\chi^N$ are connected, the sets $\underline{\mathcal J}(j)$ and
$\overline{\mathcal J}(j)$ are independent of the choice of a local coordinates system.
We also say that $j$ satisfies the condition $(\dagger)$
if the condition $(\dagger)_k$ holds for some $k \in \hat{I}^N_j$.
Note that  this definition is also independent of
a choice of $k \in \hat{I}^N_j$ by
the above observation.
Summing up,
$\underline{\mathcal J}(j)$, $\overline{\mathcal J}(j)$ and the condition $(\dagger)$ are coordinate-invariant notions and
they depend only on the families of submanifolds and the index $j$.

For $j \in \{1,2,\dots,l\}$ that satisfies $(\dagger)$,
by taking  Corollary \ref{cor:multilinear-map} into account,
we have the map
$$
\hat{\varphi}_j:
\left(\underset{X,\, \beta \in \overline{\mathcal J}(j)}{\times} T_{M_\beta}
\iota(M_\beta)\right)
\underset{X}{\times} M
\to
\left(T_{N_{j}}\iota(N_{j})
\underset{Y}{\times} N\right) \underset{Y}{\times} M
$$
where $N=N_1 \cap \dots \cap N_\ell$.
Although $\hat{\varphi}_j$ is not a morphism of vector bundles over $M$,
we still have
$$
\hat{\varphi}_j\left(\tau^X_{j_1}(\lambda_{j_1})\ldots\tau^X_{j_p}(\lambda_{j_p}) p\right)
= \tau^Y_{j}(\lambda_{j_1}\ldots\lambda_{j_p})\hat{\varphi}_j(p)
$$
where $\{j_1, \dots, j_p\} = \overline{\mathcal J}(j)$ and
$\tau_\beta^X$ and $\tau_\beta^Y$
denote the action $\tau_\beta$ on each spaces $X$ and $Y$ respectively.
Hence $\hat{\varphi}_j$ gives a multi-linear map on the fibers.
This implies, in particular, that the image of a multi-conic set is conic and
the inverse image of a conic set is multi-conic.
Summing up, we have
\begin{cor}
The restriction of $\widetilde{f}$ to $S^M \subset \widetilde{X}$ is the map
$$
\varphi_1 \underset{Y}{\times} \dots \underset{Y}{\times} \varphi_{\ell}:
S^M=\left(\underset{X,\, 1 \le \beta \le \ell}{\times} T_{M_\beta}
\iota(M_\beta)\right)
\to
S^N=\left(\underset{Y,\, 1 \le \beta \le \ell}{\times} T_{N_\beta}
\iota(N_\beta)\right),
$$
where
$
\varphi_j: S^M \to T_{N_{j}}\iota(N_{j}) \underset{Y}{\times} N
$
is given as follows.  If $j$ satisfies the condition $(\dagger)$,
the $\varphi_j$ is defined by the composition
$$
S^M
\to
\left(\underset{X,\, \beta \in \overline{\mathcal J}(j)}{\times} T_{M_\beta}
\iota(M_\beta)\right)
\underset{X}{\times} M
\overset{ \hat{\varphi}_j}{\longrightarrow}
T_{N_{j}}\iota(N_{j})
\underset{Y}{\times} N.
$$
For other $j$, the $\varphi_j$ sends a point to the corresponding one in the zero section of
$T_{N_{j}}\iota(N_{j}) \underset{Y}{\times} N$.
\end{cor}



\begin{es} Set $Y=\CC^2$, let $X=\{(z_1, z_2, \xi_1, \xi_2) \in \CC^2 \times \mathbb{P}_\CC^1,\;
        \xi_2z_1 = \xi_1 z_2\}$ and let
\begin{eqnarray*}
\pi:X & \to & \CC^2 \\
(z_1,z_2,\xi_1,\xi_2) & \mapsto & (z_1,z_2)
\end{eqnarray*}
be the desingularization map.  Let $N_1=\{0\}$, $N_2=\{z_2=0\}$, $M_1=\imin \pi(0)$, $M_2=\{\xi_2=0\}$. Locally on $X$, for example on $U_1:=\{\xi_1 \neq 0\}$, set $\lambda=\displaystyle{\xi_2\over\xi_1}$. Then $z_2=\displaystyle{\xi_2\over\xi_1}z_1$ and we have an homeomorphism
\begin{eqnarray*}
\psi:\CC^2 & \iso & U_1 \\
(\lambda,z_1) & \mapsto & (z_1,\lambda z_1,1,\lambda).
\end{eqnarray*}
We have $\imin \psi(M_1)=\{z_1=0\}$ and $\imin \psi(M_2)=\{\lambda=0\}$. We still denote them by $M_1,M_2$. The map $f:=\pi|_{U_1} \circ \psi$ is given by $(\lambda,z_1) \mapsto (\lambda z_1,z_1)$. Let us consider $\widetilde{f}$. On the zero section, let $(w_1,w_2)$ be the coordinates of $\CC^2 \simeq S'$, the zero section of $\widetilde{Y}$. We have $J^N_1=J^M_1=\{1\}$, $J^N_2=\{1,2\}=J^M_1 \sqcup J^M_2$, hence
\begin{eqnarray*}
w_1 & = & {\partial f \over \partial z_1}(0,0) z_1 = z_1 \\
w_2 & = & {\partial^2f \over \partial \lambda \partial z_1}(0,0)\lambda z_1 = \lambda z_1
\end{eqnarray*}
which is a conic map with respect to the $(\RP)^2$-actions on $S$ and $S'$.
\end{es}

One can ask when the morphism thus obtained becomes that of vector bundles.
\begin{prop}
Suppose the conditions 1.~and 2.~below.
\begin{enumerate}
 \item $f(M_j) \subset N_j$ for any $1 \le j \le \ell$.
 \item $M_{j'} \subset M_j$ if and only if $N_{j'} \subset N_j$ for $1 \le j,j' \le \ell$.
\end{enumerate}
Then the map
$
\varphi_1 \underset{Y}{\times} \dots \underset{Y}{\times} \varphi_{\ell}:
S^M
\to
S^N \underset{Y}{\times} M
$
is a morphism of vector bundles over $M$.
\end{prop}
\dim\ \ For any $k \in \{1,2,\dots,m\}$,
we have $\#\inf{}_{\subset} J^N_k = \#\sup{}_{\subset} J^N_k = 1$.
Therefore, by (\ref{eq:multi-morphism-zero}), the map $\varphi_j$ becomes
a bundle map and the result follows. \qed

\begin{oss}\label{oss: chi'<chi} More generally, the morphism exists for the
case $\#\chi^N \leq \#\chi^M$.

\

Let $f: X \to Y$ and
$\chi^M=\{M_1,\dots,M_{\ell}\}$ in $X$
and $\chi^N=\{N_1,\dots, N_{\ell'}\}$ in $Y$
for $\ell' \leq \ell$.

\

Then the situation decomposes into the following:
$$
(X;\, M_1,\dots,M_\ell) \overset{\operatorname{Id}}{\longrightarrow}
(X;\, M_1,\dots,M_{\ell'}) \overset{f}{\longrightarrow}
(Y;\, N_1,\dots,N_{\ell'})
$$
For the second arrow, we have already constructed the morphism.
We will construct the morphism for the first arrow.
In what follows, we assume that
$Y=X$, $f=\operatorname{Id}$, $\ell>\ell'$ and $N_j$ = $M_j$.

\

Locally the morphism between multi-normal deformations
$$
(x_1,\dots,x_n, t_1,\dots,t_\ell) \to
(x'_1,\dots,x'_n, t'_1,\dots,t'_{\ell'})
$$
is given by:
$$
\begin{array}{lll}
t'_j &= t_{\kappa_{\chi^M,\chi^N}(j)}
&\qquad (1 \le j \le \ell'), \\
x'_i &= t_{J_{\chi^M,\chi^N,\, i}} x_i &\qquad (1 \le i \le n).
\end{array}
$$
Here $\kappa_{\chi^M,\chi^N}$ is the map from $\{1,2,\dots,\ell'\}$ to subsets of
$\{1,2,\dots,\ell\}$ defined by
$$
\kappa_{\chi^M,\chi^N}(j) :=
\left\{\beta \in \{1,2,\dots,\ell\};\,
\left(\underset{N_k \subsetneq N_j}{\bigcup} N_k\right)
\subsetneq f(M_\beta) \subset N_j\right\}
$$
and, for $1 \le i \le n$,
$$
J_{\chi^M,\chi^N,\, i} :=
J^M_i \setminus \left(\underset{j \in J^N_i}{\bigcup} \kappa_{\chi^M,\chi^N}(j)\right)
$$
with
$J^M_i: = \{j \in \{1,2,\dots,\ell\};\, i \in I^M_j\}$
and
$J^N_i: = \{j \in \{1,2,\dots,\ell'\};\, i \in I^N_j\}$.

For example,
if $\ell' = 2$ and $M_\ell \subsetneq M_{\ell-1} \subsetneq \dots \subsetneq M_1$ are
satisfied, as $\kappa_{\chi^M,\chi^N}(1) = \{1\}$,
$\kappa_{\chi^M,\chi^N}(2) = \{2,\dots,\ell\}$,
$J_{\chi^M,\chi^N,i} = J^M_i$  $(i \notin I^M_2)$ and
$J_{\chi^M,\chi^N,i} = \emptyset$  $(i \in I^M_2)$,
the local morphism is given by
$$
\begin{array}{llllll}
&t'_1 = t_1, & t'_2 = t_2\dots t_\ell, &\qquad \\
&x'_i = t_{J^M_i} x_i\, (i \notin I^M_2),\,
&x'_i = x_i \, (i \in I^M_2).
\end{array}
$$

We can certainly glue these locally defined morphisms (by the definition
of a multi-normal deformation)
and obtain the morphism between multi-normal deformed manifolds
(say $\widetilde{X}$ and $\widetilde{Y}$) globally. Moreover,
as $J^M_i = \underset{j \in J^N_i}{\sqcup}\kappa_{\chi^M,\chi^N}(j)$ is satisfied by
definition, the following diagram commutes.
$$
\begin{matrix}
\widetilde{X} & \to & \widetilde{Y} \\
            & \searrow & \downarrow \\
            &          &   X=Y
\end{matrix}
$$
Since $\iota(M_j) \subset \iota(N_j)$ holds ($j=1,\dots,\ell'$),
we have the canonical injection
$$
T_{M_j}\iota(M_j) \hookrightarrow T_{N_j}\iota(N_j) \qquad (j \le \ell').
$$
These injections and the zero map on
$
T_{M_j}\iota(M_j)
$
for $j > \ell'$
induce the bundle map over $M:=\underset{1 \le j \le \ell}{\bigcap} M_j$
$$
\varphi: \underset{X,1\leq j\leq \ell}{\times} T_{M_j}\iota(M_j)
\to
M \underset{Y}{\times}
\left(
\underset{Y,1\leq j\leq\ell'}{\times} T_{N_j}\iota({N_j})
\right)
$$

In this case, on the zero section $S^M$, the morphism from $\widetilde{X}$ to
$\widetilde{Y}$ coincides with $\varphi$.
\end{oss}

\end{section}

\begin{section}{Multi-specialization}\label{6}

Let $X$ be a real analytic manifold with $\operatorname{dim}{X} = n$, and let
$\chi = \{M_1,\dots,M_\ell\}$ be
a family of closed submanifolds satisfying H1, H2 and H3. Let $\widetilde{X}$ be the multi-normal deformation of $X$ with respect to the family $\chi$ and consider the diagram \eqref{multi normal deformation}.

Denote by $\op(X_{sa})$ (resp. $\op^c(X_{sa})$) the category of open (resp. open relatively compact) subanalytic subsets of
$X$. One endows $\op(X_{sa})$ with the following topology: $S
\subset \op(X_{sa})$ is a covering of $U \in \op(X_{sa})$ if for
any compact $K$ of $X$ there exists a finite subset $S_0\subset S$
such that $K \cap \bigcup_{V \in S_0}V=K \cap U$. We will call
$X_{sa}$ the subanalytic site and denote by $\rho:X \to X_{sa}$ the natural morphism of sites associated to the inclusion $\op(X_{sa}) \hookrightarrow \op(X)$. Let $\mod(k_{X_{sa}})$ (resp. $D^b(k_{X_{sa}})$ denote the category of sheaves on $X_{sa}$ (resp. bounded derived category of sheaves on $X_{sa}$). Reference for classical sheaf theory are made to \cite{KS90}, for sheaves on subanalytic sites we refer to \cite{KS01} and \cite{Pr08}.
For an exposition on $(\RP)^\ell$-conic sheaves see the Appendix.

\begin{lem} Let $F \in D^b(k_{X_{sa}})$. There is a natural isomorphism
$$
\imin s R\Gamma_\Omega \imin p F \simeq s^!(\imin p F)_\Omega.
$$
\end{lem}

\dim\ \ We prove the assertion in several steps. For $I,J \subseteq \{1,\dots,\ell\}$, $I \cap J = \emptyset$, set $S_I=\{t_i=0,\,i \in I\}$ and $\Omega_J=\{t_j>0,\,j\in J\}$.


(i) We show that if $\sharp J < \ell$, then $(R\Gamma_{\Omega_J}(\imin pF)_\Omega)|_S=0$. By d\'evissage we may reduce to $F \in \mod(k_{X_{sa}})$, so $F=\lind i \rho_* F_i$ with $F_i \in \mod_{\rc}(k_{X_{sa}})$
for each $i$. Since $R^k\Gamma_{\Omega_J}$ commutes with filtrant $\Lind$, we may reduce to $F \in \mod_{\rc}(k_X)$. Taking a suitable triangulation adapted to $\Omega_J$ and $\Omega$ we may assume that $F=k_U$, where $U$ is an open star associated to the triangulation. In this case one checks that
$$(R\Gamma_{\Omega_J}(\imin pF)_{\Omega_J\setminus\Omega})|_S \simeq (R\Gamma_{\Omega_J}\imin pF)|_S.$$

(ii) We show that if $\sharp (I \cup J) < \ell$, then $(R\Gamma_{S_I \cap \Omega_J}(\imin pF)_\Omega)|_S=0$. We argue by induction on $\sharp I$. If $\sharp I=0$ then it follows by (i) that $(R\Gamma_{\Omega_J}(\imin pF)_\Omega)|_S=0$. Suppose that it is true for $\sharp I \leq \ell-2$. Let $i_0 \in I$. We have the distinguished triangle
$$
\dt{(R\Gamma_{S_I \cap \Omega_J}(\imin pF)_\Omega)|_S}{(R\Gamma_{S_{I \setminus \{i_0\}} \cap \Omega_J}(\imin pF)_\Omega)|_S}{(R\Gamma_{S_{I \setminus \{i_0\}} \cap \Omega_{J \cup \{i_0\}}}(\imin pF)_\Omega)|_S}.
$$
The second and the third terms of the  triangle are zero since $\sharp (I \setminus \{i_0\}) \leq \ell-2$. Then the first one is zero.

(iii) We show that if $I \cup J=\ell$ then $(R\Gamma_{S_I \cap \Omega_J}(\imin p F)_\Omega)|_S \simeq (R\Gamma_S(\imin p F)_\Omega)|_S[\sharp J]$. We argue by induction on $\sharp J$. Let $j_0 \in J$. We have the distinguished triangle
$$
\dt{(R\Gamma_{S_{I \cup \{j_0\}} \cap \Omega_{J \setminus \{j_0\}}}(\imin pF)_\Omega)|_S}{(R\Gamma_{S_I \cap \Omega_{J \setminus \{j_0\}}}(\imin pF)_\Omega)|_S}{(R\Gamma_{S_I \cap \Omega_J}(\imin pF)_\Omega)|_S}
$$
where the second term is zero by (ii). If $\sharp J=1$ the result follows immediately.
Suppose that it is true for $\sharp J \leq \ell-1$. Then $(R\Gamma_{S_{I \cup \{j_0\}} \cap \Omega_{J \setminus \{j_0\}}}(\imin pF)_\Omega)|_S[1] \simeq (R\Gamma_{S_I \cap \Omega_J}(\imin pF)_\Omega)|_S$ by (ii) and $(R\Gamma_{S_{I \cup \{j_0\}} \cap \Omega_{J \setminus \{j_0\}}}(\imin p F)_\Omega)|_S \simeq (R\Gamma_S(\imin p F)_\Omega)|_S[\sharp J-1]$ by the induction hypothesis.

(iv) By (iii) we obtain $s^!(\imin p F)_\Omega[\ell] \simeq \imin sR\Gamma_\Omega\imin pF$. It remains to show that $(\imin p F)_\Omega[\ell] \simeq (p^!F)_\Omega$. We may reduce to the case $F \simeq \lind i \rho_*F_i$ with $F_i\in\mod_{\rc}(k_X)$. In this case we have
\begin{eqnarray*}
H^{k+\ell}(\imin p F)_\Omega & \simeq & \lind i\rho_*H^{k+\ell}(\imin p F_i)_\Omega \\
& \simeq & \lind i \rho_*H^k(p^!F_i)_\Omega \\
& \simeq & H^k(p^!F)_\Omega
\end{eqnarray*}
where the second isomorphism follows since $\imin {\widetilde{p}} [\ell] \simeq \widetilde{p}^!$ in $\mod(k_X)$. \qed

\begin{df} The (multi-)specialization along $\chi$ is the functor
\begin{eqnarray*}
\nu^{sa}_\chi:D^b(k_{X_{sa}}) & \to & D^b(k_{S_{sa}}) \\
F & \mapsto & \imin s\mathrm{R}\Gamma_\Omega\imin p F.
\end{eqnarray*}
\end{df}

\begin{teo} \label{m4.2.3} Let $F \in D^b(k_{X_{sa}})$.
\begin{itemize}
\item[(i)] $\nu^{sa}_\chi F \in D^b_{(\RP)^\ell}(k_{S_{sa}})$ (see Definition \ref{def:conic sheaf}).
\item[(ii)] Let $V$ be an $(\RP)^{\ell}$-conic open subanalytic subset
in $S$.
Then:
$$H^j(V;\nu^{sa}_\chi F) \simeq \lind U H^j(U;F),$$
where $U$ ranges through the family of $\op(X_{sa})$ such that
$C_\chi(X \setminus U) \cap V=\emptyset$.
\end{itemize}
\end{teo}
\dim\ \ (i) We may reduce to the case $F \in \mod(k_{X_{sa}})$.
Hence $F=\lind i \rho_* F_i$ with $F_i \in \mod_{\rc}(k_{X_{sa}})$
for each $i$.
We have $\imin p \lind i \rho_* F_i \simeq \lind i \rho_* \imin p
F_i$ and $\imin p F_i$ is $\R$-constructible and $(\RP)^\ell$-conic for each
$i$. Hence $\imin p F$ is $(\RP)^\ell$-conic. Since the functors
$\mathrm{R}\Gamma_\Omega$ and $\imin s$ send $(\RP)^\ell$-conic sheaves to $(\RP)^\ell$-conic sheaves
we obtain $\imin s \mathrm{R}\Gamma_\Omega \imin pF=\nu^{sa}_\chi F \in
D^b_{(\RP)^\ell}(k_{S_{sa}})$.

(ii) Let $U \in \op(X_{sa})$ such that $V \cap C_\chi(X \setminus
U)=\emptyset$. We have the chain of morphisms
\begin{eqnarray*}
\mathrm{R}\Gamma(U;F) & \to & \mathrm{R}\Gamma(\imin p(U);\imin p F) \\
& \to & \mathrm{R}\Gamma(\imin p(U) \cap \Omega;\imin p F) \\
& \to & \mathrm{R}\Gamma(\imin {\widetilde{p}}(U) \cup V;\mathrm{R}\Gamma_\Omega\imin p F) \\
& \to & \mathrm{R}\Gamma(V;\nu^{sa}_\chi F)
\end{eqnarray*}
where the third arrow exists since $\imin {\widetilde{p}}(U) \cup
V$ is a neighborhood of $V$ in $\overline{\Omega}$ by Proposition
\ref{m4.1.3} (ii).

Let us show that it is an isomorphism.
By Corollaries \ref{cor: Rpl neighborhood} and \ref{cor:intersection-V-U},
we may assume that $V$ is one of the $V^{(\alpha)}$'s given in Corollary \ref{cor: Rpl neighborhood}.
We have
\begin{eqnarray*}
H^k(V;\nu^{sa}_\chi F) & \simeq & \lind W H^k(W;\mathrm{R}\Gamma_\Omega\imin p F) \\
& \simeq & \lind W H^k(W \cap \Omega;\imin p F),
\end{eqnarray*}
where $W$ ranges through the family of subanalytic open
neighborhoods of $V$ in $\widetilde{X}$. By Corollary \ref{cor: Rpl neighborhood},
we may assume that $W$ satisfies \eqref{mpitilde}.
Since $\imin p F$ is $(\RP)^\ell$-conic, we have
\begin{eqnarray*}
H^k(W \cap \Omega;\imin p F) & \simeq & H^k(\imin p(p(W \cap
\Omega));\imin p F) \\
& \simeq & H^k(p(W \cap \Omega) \times \{(1)_\ell\};\imin p F) \\
& \simeq & H^k(p(W \cap \Omega);F),
\end{eqnarray*}
where $(1)_\ell=(1,\dots,1) \in \R^\ell$. The second isomorphism follows since every subanalytic
neighborhood of $p(W \cap \Omega) \times \{(1)_\ell\}$ contains an
$(\RP)^\ell$-connected subanalytic neighborhood (the proof is similar to
that of Proposition \ref{m4.1.4}). By Proposition \ref{m4.1.3} (i) we have that
$p(W \cap \Omega)$ ranges through the family of subanalytic open
subsets $U$ of $X$ such that $V \cap C_\chi(X \setminus
U)=\emptyset$ and we obtain the result.
\qed

Remember that a sheaf $F \in \mod(k_{X_{sa}})$ is said to be quasi-injective if the restriction morphism $\Gamma(U;F) \to \Gamma(V;F)$ is surjective for each $U,V \in \op^c(X_{sa})$ with $U \supseteq V$.

\begin{cor} Let $F \in \mod(k_{X_{sa}})$ be quasi-injective. Then $F$ is $\nu^{sa}_\chi$-acyclic.
\end{cor}
\dim\ \ The result follows from Theorem \ref{m4.2.3} and the fact that quasi-injective sheaves are $\Gamma(U;\cdot)$-acyclic for each $U \in \op(X_{sa})$. \qed

\begin{cor} Let $F \in D^b(k_{X_{sa}})$ and let $p=(q;\xi) \in S$. Then
$$
(\imin \rho H^j\nu^{sa}_\chi F)_p \simeq \lind {W,\epsilon} H^j(W \cap B_\epsilon;F),
$$
where $W$ ranges through the family $\operatorname{Cone}\chi(p)$ and $B_\epsilon$ ranges through the family of open balls of radius $\epsilon>0$ containing $q$. Here we locally identify $X$ with a vector space as in Section \ref{4}.
\end{cor}
\dim\ \ The result follows since for any subanalytic conic neighborhood $V$ of $p$, any $U \in \op(X_{sa})$ such that $C_\chi(X \setminus U) \cap V=\emptyset$ contains $W \cap B_\epsilon$, $q \in B_\epsilon$, $\epsilon>0$, $W \in \operatorname{Cone}_\chi(p)$. \qed


\begin{cor} Let $F \in D^b(k_{X_{sa}})$. Let $k \le \ell$ and $\{j_1,\dots,j_k\}$ be a subset of $\{1,2,\dots,\ell\}$.
Set $\chi_k=\{M_{j_1},\dots,M_{j_k}\}$.
Then we have
$$
(\nu^{sa}_\chi F)|_{S_{\chi/\chi_k}}
\simeq (\nu^{sa}_{\chi_k}F)|_{S_{\chi/\chi_k}}
$$
Here the subset $S_{\chi/\chi_k}$ was defined by (\ref{def:S-chi-chi}).
\end{cor}
\dim\ \
Set $Z=S_{\chi/\chi_k}$.
On $Z$ the multi-actions induced by $\{\mu_1,\dots,\mu_\ell\}$ and $\{\mu_{j_1},\dots,\mu_{j_k}\}$ coincide, and
we can regard $Z$ as an $(\RP)^\ell$-conic
(resp.  $(\RP)^k$-conic) sub-bundle of
$S_\chi$ (resp. $S_{\chi_k} \underset{X}{\times}{M}$).
Let $V$ be a $(\RP)^\ell$-conic subanalytic open subset of $Z$.
We may assume that
$V$ is one of $V^{(\alpha)}$'s given in Corollary \ref{cor: Rpl neighborhood}.
Let $U \in \op(X_{sa})$ such that $V \cap C_\chi(X \setminus U)=\emptyset$.
Let $W$ be the complement of $C_\chi(X \setminus U)$.
Then $W$ is a $(\RP)^\ell$-conic open subanalytic neighborhood of $V$ such that $C_\chi(X \setminus U) \cap W=\emptyset$.
Then for $j \in \Z$
\begin{eqnarray*}
H^j(V;\nu^{sa}_\chi F) & \simeq & \lind W H^j(W;\nu^{sa}_\chi F) \\
& \simeq & \lind U H^k(U;F) \\
& \simeq & \lind {U'} H^k(U';F),
\end{eqnarray*}
where $W$ ranges through the family of $(\RP)^\ell$-conic open subanalytic neighborhoods of $V$, $U \in \op(X_{sa})$ is such that $C_\chi(X\setminus U) \cap W=\emptyset$ and $U' \in \op(X_{sa})$ is such that $C_{\chi}(X\setminus U') \cap V=\emptyset$. Remark that the first isomorphism is a consequence of the fact that $\nu^{sa}_\chi F$ is conic and the fact that, by 2.~of Corollary \ref{cor: Rpl neighborhood},
every open subanalytic neighborhood of $V$ contains a $(\RP)^\ell$-connected one.

Clearly $V$ satisfies the condition either Va.~or Vb.~for each $j_1,\dots, j_k$
with respect to the family $\chi_k$. Hence, in the same way, we see that
\begin{eqnarray*}
H^j(V;\nu^{sa}_{\chi_k} F)
& \simeq & \lind {U'} H^k(U';F),
\end{eqnarray*}
where $U' \in \op(X_{sa})$ is such that $C_{\chi_k}(X\setminus U') \cap V=\emptyset$. To show the result it is enough to see that
$$
C_\chi(X \setminus U) \cap V = \emptyset \ \ \Longleftrightarrow \ \ C_{\chi_k}(X \setminus U) \cap V = \emptyset.
$$
This follows from Corollary \ref{cor: restriction multi-cones}, which says that
$$
C_\chi(X \setminus U) \cap Z = C_{\chi_k}(X \setminus U) \cap Z.
$$
Hence we have the isomorphism
$$
H^j(V;\nu^{sa}_\chi F) \simeq
H^j(V;\nu^{sa}_{\chi_k} F)
$$
for each $(\RP)^\ell$-conic globally subanalytic open subset of $Z$. By Lemma \ref{lem: RPl decomposition} and the fact that both $\nu^{sa}_\chi F$ and $\nu^{sa}_{\chi_k} F$ are multi-conic we obtain $(\nu^{sa}_\chi F)|_Z \simeq (\nu^{sa}_{\chi_k} F)|_Z$.
\qed

Let $f:X \to Y$ be a morphism of real analytic manifolds, $\chi^M=\{M_1,\ldots,M_\ell\}$, $\chi^N=\{N_1,\ldots,N_\ell\}$ two families of closed analytic submanifolds of $X$ and $Y$ respectively satisfying the hypothesis H1, H2 and H3. Suppose that $f(M_i) \subseteq N_i$, $i=1,\ldots,\ell$.
 We call $T_\chi f$ the induced map on the zero section. Let $J \subseteq \{1,\ldots,\ell\}$ and set $\chi_J=\{M_{j_1},\ldots,M_{j_k}\}$, $j_k \in J$. The map $T_{\chi_J}f$ denotes the restriction of $\widetilde{f}$ to $\{t_{j_k}=0,\,j_k\in J\}$. In the following we will denote with the same symbol $C_{\chi_J}(Z)$ the normal cone with respect to $\chi_J$ and its inverse image via the map $\widetilde{X} \to \widetilde{X}_{M_{j_1},\dots,M_{j_k}}$.

\begin{prop} Let $F \in D^b(k_{X_{sa}})$.

\begin{itemize}
\item[(i)] There exists a commutative diagram of canonical
morphisms
$$
\xymatrix{R(T_\chi f)_{!!}\nu^{sa}_{\chi^M} F  \ar[d] \ar[r] & \nu^{sa}_{\chi^N}Rf_{!!}F \ar[d] \\
R(T_\chi f)_*\nu^{sa}_{\chi^M} F & \nu^{sa}_{\chi^N}Rf_*F. \ar[l]}
$$
\item[(ii)] Moreover if $f:\supp F \to Y$ and $T_{\chi_J}f:C_{\chi_J}(\supp F) \to
\{t_{j_1}=\cdots=t_{j_k}=0\}$
for each $J = \{j_1,\dots,j_k\}$ are proper, and if $\supp F \cap \imin f(N_j) \subseteq M_j$, $j \in \{1,\ldots,\ell\}$, then
the above morphisms are isomorphisms.
\end{itemize}
\end{prop}
\dim\ \ (i) The existence of the arrows is done as in \cite{KS90} Proposition 4.2.4.

(ii) If $\widetilde{p}^{-1}(\supp F)$ is proper over $\widetilde{X}$, then all the morphisms are isomorphisms. We have to prove that for a closed subset $Z$ of $X$, the restriction of $\widetilde{f}$ to $\overline{\imin {\widetilde{p}}(Z)}$ is proper, if $Z \to Y$ and $C_{\chi_J}(Z) \to
\{t_{j_1}=\cdots=t_{j_k}=0\}$
for $k\leq \ell$ are proper, and if $Z \cap \imin f(N_j) \subseteq M_j$, $j \in \{1,\ldots,\ell\}$. We argue by induction on $\sharp \chi$. If $\sharp \chi=1$ this is Proposition 4.2.4 of \cite{KS90}. Suppose it is true for $\sharp\chi\leq\ell-1$. It follows from the hypothesis that the fibers of $\widetilde{f}$ restricted to $\overline{\widetilde{p}^{-1}(Z)}$ are compact (if $t_{j_1}=\cdots=t_{j_k}=0$ this is a consequence of the fact that $C_{M_{j_1}\dots M_{j_k}}(Z) \to
\{t_{j_1}=\cdots=t_{j_k}=0\}$
 is proper). Then it remains to prove that it is a closed map. Let $\{u_n\}_{n\in\N}$ be a sequence in $\overline{\imin{\widetilde{p}}(Z)}$ and suppose that $\{\widetilde{f}(u_n)\}_{n\in\N}$ converges. We shall find a convergent subsequence of $\{u_n\}_{n\in\N}$. We may also assume that $\{\widetilde{p}(u_n)\}_{n\in\N}$ converges. The map $\imin{\widetilde{p}}(Z) \setminus \{t_1=\cdots=t_\ell=0\} \to \widetilde{Y} \setminus \{t_1=\cdots=t_\ell=0\}$ is proper. Indeed, let $K$ be a compact subset of $\widetilde{Y} \setminus \{t_1=\cdots=t_\ell=0\}$, and
reduce to the case that $K$ is contained in $\{c \leq t_j \leq d\} \subset \widetilde{Y}$, $c,d>0$, $j \in \{1,\dots,\ell\}$. Suppose without loss of generality that $j=1$. Then $K_1:=\widetilde{p}_{N_1}(K)$ is a compact subset of $\widetilde{Y}_{N_2\ldots N_\ell}$. Let us identify $K_1$ with $\imin{\widetilde{p}_{N_1}}(K_1) \cap \{t_1=1\}$. Then $\imin{\widetilde{f}}(K_1)$ is compact by the induction hypothesis. Hence $\imin{\widetilde{f}}(K) \subseteq \mu_1(\imin{\widetilde{f}}(K_1),[c,d])$ is compact since it is closed and contained in a compact subset. 
We may assume that $\{\widetilde{f}(u_n)\}_{n\in\N}$ converges to a point of $\{t_1=\cdots=t_\ell=0\}$. Then $\{\widetilde{p}(u_n)\}_{n\in\N}$ converges to a point of $Z \cap \imin f(N_1 \cap \cdots \cap N_\ell) \subseteq M_1 \cap \cdots \cap M_\ell$.

\

Taking local coordinates systems of $X$ and $Y$, let $u_n=(x_{1n},\ldots,x_{mn},t_{1n},\ldots,t_{\ell n})$, $t_{jn}>0$, $j=1,\dots,\ell$. Then $t_{jn} \to 0$, $j=1,\dots,\ell$ and $t_{J^M_in}x_{in}\to 0$, $i=1,\ldots,m$. It is enough to show that $\{|x_{in}|\}_{n\in\N}$ is bounded for each $i=1,\dots,m$. We argue by contradiction. Suppose without loss of generality that $|x_{1n}|\to+\infty$ and that $\{x_{1n}/|x_n|\}_{n\in\N}$ is not infinitesimal (i.e. its limit is not 0).
Set $u_n^{(\beta)}=(x_n^{(\beta)},t_n^{(\beta)})$, where $x_{in}^{(\beta)}=t_{in\beta}x_{in}$ with $t_{in\beta}=t_{\beta n}$ if $i \in I_\beta$ and $t_{in\beta}=1$ if $i \notin I_\beta$, and where $t^{(\beta)}_{jn}=t_{jn}$ if $j \neq \beta$ and $t^{(\beta)}_{\beta n}=1$. Then $u_n$ belongs to $\imin{\widetilde{p}}(Z) \setminus \{t_1=\cdots=t_\ell=0\}$ and $\{\widetilde{f}(u_n^{(\beta)})\}_{n\in\N}$ converges. Since $\widetilde{f}$ is proper on $\imin{\widetilde{p}}(Z) \setminus \{t_1=\cdots=t_\ell=0\}$, then $\{|x_n^{(\beta)}|\}_{n\in\N}$ is bounded.



(a) Suppose that there exists $\beta \in \{1,\dots,\ell\}$ such that $1 \notin I^M_\beta$. The sequence $\{|x_n^{(\beta)}|\}_{n\in\N}$ is bounded. In particular $\{|x_{1n}|\}_{n\in\N}$ is bounded which is a contradiction.

(b) Suppose that $1 \in I^M_j$ for $j=1,\dots,\ell$. Suppose without loss of generality that $I^M_1$ is the biggest $I^M_j$ containing $1$. Set $\widetilde{u}_n=(\widetilde{x}_n,\widetilde{t}_n)$ where $\widetilde{x}_{in}=x_{in}/|x_n|$ if $i \in I^M_1$ and $\widetilde{x}_{in}=x_{in}$ otherwise and where $\widetilde{t}_{1n}=t_{1n}|x_n|$, $\widetilde{t}_{jn}=t_{jn}$ if $j \neq 1$. Remark that $|t_{1n}x_n|=|t_{1n}x_{1n}||x_n|/|x_{1n}|$ is bounded since the sequences $\{|x^{(1)}_{1n}|\}_{n\in\N}$ and $\{|x_n|/x_{1n}\}_{n\in\N}$ are bounded. By extracting a subsequence $\{\widetilde{u}_n\}_{n\in\N}$ converges to a non-zero vector $v$ and $\widetilde{u}_n$ belongs to $\imin{\widetilde{p}}(Z)$. On the other hand $\widetilde{f}(u_n)$ converges and hence $f_k(t_{J^M_1n}x_{1n},\dots,t_{J^M_mn}x_{mn})/t_{J^N_kn}|x_n|$ converges to $0$ for each $k$ such that 
$k \in I_1^N$. Moreover for each $k \in I^N_2 \setminus I^N_1$ by \eqref{eq:morphism-expansions} the sequence $f_k(t_{J_1^mn}x_{1n},\dots,t_{J_m^Mn}x_{mn})/t_{J_k^Nn}$ converges to 0. Hence $\widetilde{f}(\mu_2(v,\lambda))=\widetilde{f}(v)$ for each $\lambda \in \RP$ which contradicts the fact that the fibers of $\overline{\imin{\widetilde{p}}(Z)} \to \widetilde{Y}$ are compact.
\qed

\begin{prop} Let $F \in D^b(k_{Y_{sa}})$.

\begin{itemize}
\item[(i)] There exists a commutative diagram of canonical
morphisms
$$
\xymatrix{\omega_{S_X/S_Y}\otimes\imin {(T_\chi f)}\nu^{sa}_{\chi^N} F  \ar[d] \ar[r] & \nu^{sa}_{\chi^M}(\omega_{S_X/S_Y}\otimes \imin f F) \ar[d] \\
T_{\chi} f^!\nu^{sa}_{\chi^N} F & \nu^{sa}_{\chi^M}f^!F. \ar[l]}
$$
\item[(ii)]
The above morphisms are isomorphisms on the open sets where $T_{\chi_J} f$ is smooth for each $J \subseteq \{1,\ldots,\ell\}$.
\end{itemize}
\end{prop}
\dim\ \ (i) The existence of the arrows is done as in \cite{KS90} Proposition 4.2.5. When (ii) is satisfied the function $\widetilde{f}$ is smooth at any point of the boundary of $\Omega$ and all the above morphisms become isomorphisms.
\qed

\end{section}

\begin{section}{Multi-asymptotic expansions}\label{7}

Let $I_j$ ($j=1,\dots,\ell)$ be a subset of $\{1,\dots,n\}$ which
satisfies the conditions (\ref{eq:conditions-indices-set}), and
let $X$ be the $n$-dimensional complex vector space $\CC^n$ with coordinates
$(z_1, z_2, \dots, z_n)$. Let $\chi = \{Z_1, \dots, Z_{\ell}\}$ be
a family of $\ell$ complex submanifolds
defined by
$$
Z_j := \{z \in X;\, z_{i} = 0 \text{ for $i \in I_j$}\}.
$$
Remember the definition of $\hat{I}_j$ given by (\ref{eq:def-hat-I})
and that of $J_i$ by (\ref{eq:def-J_i}).
Then, by the condition (\ref{eq:conditions-indices-set}),
we have $T_{Z_j}\iota(Z_j)\simeq Z_j \times \CC^{\#\hat{I}_j}$ and
\begin{equation}\label{eq:disjoint-union-indices}
\underset{1 \le j \le \ell}{\bigcup} I_j = \hat{I}_1 \sqcup \hat{I}_2 \sqcup \dots
\sqcup \hat{I}_\ell.
\end{equation}
For a subset $I = \{i_1, i_2, \dots, i_m\}$ of $\{1,2,\dots,n\}$,
we denote by $z_I$ the coordinates $(z_{i_1}, z_{i_2}, \dots, z_{i_m})$ and
by $\CC^I$ the complex space $\CC^{\#I}$ with coordinates $z_I$.
The map $\pi_{I}$ designates the projection from $X$ to $\CC^I$.

Let $\mu_j(z,\lambda): X \times \R \to X$ ($j=1,2,\dots,\ell$) be the action defined by
$$
\mu_j((z_1, \dots, z_n), \lambda) = (\lambda_{j1}z_1, \dots, \lambda_{jn}z_n)
$$ where $\lambda_{ji} = \lambda$
if $i \in I_j$ and $\lambda_{ji} = 1$ otherwise.
Let $U$ be an open subset of product type defined by
\begin{equation}\label{eq:def-open-set-product-type}
U=B_1 \times \dots \times B_\ell \times W
\subset \mathbb{C}^{\hat{I}_1} \times \dots \times
\mathbb{C}^{\hat{I}_\ell} \times \mathbb{C}^d = X,
\end{equation}
where $d = n - \displaystyle\sum_{j=1}^\ell \# \hat{I}_j$ and
$B_j$ is an open ball with center at the origin in
$\mathbb{C}^{\hat{I}_j}$ and $W$ is a convex open subset in $\mathbb{C}^d$.
Let $G_j$ ($j=1,2, \dots,\ell$) be an open proper convex cone in $\CC^{\hat{I}_j}$.
Note that, in our paper, we only consider a proper convex cone for simplicity.

We define an open proper multi-cone
$S(U,\{G_j\}, \epsilon)$ ($\epsilon > 0$) in $X$
as follows. Set, for $j=1,2,\dots,\ell$,
$$
V_j :=
\left\{
z \in U;\,
z_{\hat{I}_j} \in G_j, \,
\vert z_{\hat{I}_\beta} \vert < \epsilon \vert z_{\hat{I}_j}\vert
\text { for } \beta \text
{ with } I_\beta \subsetneq I_j
\right\}
\subset X
$$
Then we define
$$
S(U,\{G_j\}, \epsilon) := V_1 \cap V_2 \cap \dots \cap V_\ell \subset X.
$$
If $\ell = 0$, i.e., no submanifolds in $\chi$,
then we set $S(U, \{G_j\}, \epsilon) := U$ by convention.
Note that $S:=S(U, \{G_j\}, \epsilon)$ is non-empty and
$$
S \cap \left(\underset{1\le j \le \ell}{\bigcup} Z_j\right) = \emptyset
$$
follows from the condition (\ref{eq:disjoint-union-indices}).
For example, $S$ gives just
a poly-sector when ${\underset{1 \le j \le \ell}{\bigcup} Z_j}$ forms a
normal-crossing divisor. The multi-cone $S$ plays the same
geometrical role as that of a sector in an asymptotic expansion.

Let $\Z_\ell$ denote the set of integers $\{1,2,\dots, \ell\}$ and
$\Powl$ be the set of all the subsets of
$\Z_\ell$ except for the empty set.
For a $J = \{j_1, j_2, \dots, j_m\} \in \Powl$ and
a proper multi-cone $S:=S(U, \{G_j\}, \epsilon)$,
we define
$$
\begin{aligned}
Z_{J} &:= Z_{j_1} \cap Z_{j_2} \cap \dots \cap Z_{j_m}, \\
I_{J} &:= I_{j_1} \cup I_{j_2} \cup \dots \cup I_{j_m}, \\
S_{J} &:= \operatorname{Int}_{Z_J}\left(\overline{S} \cap Z_{J}\right)
\end{aligned}
$$
where $\overline{A}$ is the closure of a set $A$ in $X$ and
$\operatorname{Int}_{Z_J}(B)$ denotes the interior of
a set $B$ in $Z_J$.  We often write $S_J$ by $S_j$ if $J = \{j\}$ from now on. 
We give some fundamental properties of a multi-cone which are needed later.
\begin{lem}{\label{lemma:fundamental-geo}}
For $S := S(U,\{G_j\}, \epsilon)$ and $J \in \Powl$,
we have the followings.
\begin{enumerate}
\item $S_J$ is also a non-empty open proper multi-cone in $Z_J$
for submanifolds $\{Z_j \cap Z_J\}_{j \in J^*}$. To be more precise,
we have, in $Z_J$,
\begin{equation}{\label{eq:subsector-set-def}}
S_J = S(U\cap Z_J, \{G_j\}_{j \in J^*}, \epsilon)
\end{equation}
where $J^*$ is given by
\begin{equation}{\label{eq:def-J-star}}
J^* := \{j \in \Z_\ell;\, Z_J \nsubseteq Z_j\} = \{j \in \Z_\ell;\,
I_j \nsubseteq I_J\}.
\end{equation}
Note that
$
I^*_j := I_j \setminus I_J$ $(j \in J^*)
$
also satisfies the conditions (\ref{eq:conditions-indices-set}), and
it defines the submanifold $Z_j \cap Z_J$.
\item $\overline{S_J} = \overline{S} \cap Z_J$ and $Z_{J'} \cap S_J = \emptyset$
for $J' \in \Powl$ with $Z_{J} \nsubseteq Z_{J'}$.
\item $\{S_J\}_{J \in \Powl}$
covers the edge $\overline{S} \cap \left(\underset{1\le j \le \ell}{\bigcup} Z_j\right) \cap U$ of $S$, i.e.,
\begin{equation}
\overline{S} \cap \left(\underset{1\le j \le \ell}{\bigcup} Z_j\right) \cap U=
\underset{J \in \Powl}{\bigcup} S_J.
\end{equation}
\end{enumerate}
\end{lem}
\dim\ \
Note that $Z_J \nsubseteq Z_j$ is equivalent to
$\hat{I}_j \cap I_{J} = \emptyset$ because of
the conditions (\ref{eq:conditions-indices-set}).
We first show 1.~of the lemma.

\begin{lem}\label{lem:lemma-multicone-closure}
we have $\overline{S} = \overline{V}_1 \cap \dots \cap \overline{V}_\ell$.
\end{lem}
\dim\ \
Clearly we have
$$
\begin{aligned}
&\overline{V}_1 \cap \dots \cap \overline{V}_\ell\\
&\quad = \bigcap_{1 \le j \le \ell} \left\{
z \in \overline{U};\,
z_{\hat{I}_j} \in \overline{G}_j, \,
\vert z_{\hat{I}_\beta} \vert \le \epsilon \vert z_{\hat{I}_j}\vert
\text { for } \beta \text
{ with } I_\beta \subsetneq I_j
\right\} \\
&\quad =
\left\{z \in \overline{\left(G_1 \times \dots \times G_\ell \times \mathbb{C}^d\right) \cap U};\,
\vert z_{\hat{I}_\beta} \vert \le \epsilon \vert z_{\hat{I}_\alpha}\vert
\text { for } I_\beta \subsetneq I_\alpha
\right\},
\end{aligned}
$$
where $d = n - \displaystyle\sum_{j=1}^\ell \#\hat{I}_j$ and
$X = \mathbb{C}^{\hat{I}_1} \times \dots \times \mathbb{C}^{{I}_\ell}
\times \mathbb{C}^d$.
In the same way, $S$ is defined by
$$
\left\{z \in \left(G_1 \times \dots \times G_\ell \times \mathbb{C}^d\right) \cap U;\,
\vert z_{\hat{I}_\beta} \vert < \epsilon \vert z_{\hat{I}_\alpha}\vert
\text { for } I_\beta \subsetneq I_\alpha
\right\}.
$$
The claim of the lemma follows from these expressions.
\qed

By the lemma, we get
$$
S_J = \operatorname{Int}_{Z_J}(\overline{V}_1 \cap \dots \cap \overline{V}_\ell \cap Z_{J})
=
\operatorname{Int}_{Z_J}(\overline{V}_1 \cap Z_{J}) \cap
\dots \cap
\operatorname{Int}_{Z_J}(\overline{V}_\ell \cap Z_{J}).
$$
As
$$
\operatorname{Int}_{Z_J}(\overline{V}_j \cap Z_{J})
=\left\{
\begin{aligned}
V_j \cap Z_{J}   \qquad & (\hat{I}_j \cap I_{J} = \emptyset),\\
U \cap Z_{J}     \qquad & (\text{otherwise}),
\end{aligned}
\right.
$$
holds, we have
$$
S_J = \underset{\hat{I}_j \cap I_J = \emptyset}{\bigcap} (V_j \cap Z_J).
$$
Hence $S_{J}$ is
the multi-cone $S(U \cap Z_{J}, \{G_j\}_{j \in J^*}, \epsilon)$
of $Z_J$.

Next we show 2.~of the lemma.
We have, by the same argument used in the proof for
Lemma \ref{lem:lemma-multicone-closure},
$$
\overline{\underset{\hat{I}_j \cap I_{J} = \emptyset}{\bigcap} V_j \cap Z_{J}}
= \underset{\hat{I}_j \cap I_{J} = \emptyset}{\bigcap} \overline{V_j} \cap Z_{J}.
$$
Hence, by noticing
$
\overline{V}_j \cap Z_{J} = \overline{U} \cap Z_{J}
$
if $\hat{I}_j \cap I_{J} \ne \emptyset$, we obtain
$$
\overline{S_J}
= \overline{\underset{\hat{I}_j \cap I_{J} = \emptyset}{\bigcap} V_j \cap Z_{J}}
= \underset{\hat{I}_j \cap I_{J} = \emptyset}{\bigcap} \overline{V_j} \cap Z_{J}
= \underset{1 \le j \le \ell}{\bigcap} \overline{V_j} \cap Z_{J}
= \overline{S} \cap Z_{J}.
$$
For $Z_J \nsubseteq Z_{J'}$,
we can find $j \in J'$ satisfying $Z_J \nsubseteq Z_j$. Then $j$ belongs to
$J^*$ by definition,
from which we obtain $(Z_j \cap Z_J) \cap S_J = \emptyset$
since $S_J$ is an open proper
multi-cone in $Z_J$ by 1.~of the lemma.

Noticing that $\overline{S_J} = \overline{S} \cap Z_J$ and $S_J$ is still
a multi-cone, 3.~of the lemma can be shown by induction.
\qed

The following properties of $S := S(U, \{G_j\}, \epsilon)$ are easily verified.
\begin{enumerate}
\item For any $j$ and $z \in S$, we have $\mu_j(z,(0,1]) \subset S$.
\item For any $\tilde{z} \in S_j$ ($1 \le j \le \ell$),
	there exists a point $z \in S$ satisfying
	$\mu_j(z,0) = \tilde{z}$.
\end{enumerate}
For $J = \{j_1, j_2, \dots, j_m\} \in \Powl$,
we set
$$
\mu_J(z,(\lambda_{j_1}, \lambda_{j_2}, \dots, \lambda_{j_m}))
=\mu_{j_1}(\mu_{j_2}(\dots\mu_{j_m}(z,\lambda_{j_m}),\dots,\lambda_{j_2}),\lambda_{j_1}).
$$
Then the above properties 1.~and 2.~also hold for these multi-actions.
\begin{lem}{\label{lemma:mutli-cone-geometric}}
Let $J \in \Powl$ and $S := S(U, \{G_j\}, \epsilon)$.
\begin{enumerate}
	\item For any $z \in S$, we have
		$\mu_J(z, (0,1] \times \cdots \times (0,1]) \subset S$.
	\item For any $\tilde{z} \in S_J$, there exists $z \in S$ with
		$\mu_J(z,(0,\dots,0)) = \tilde{z}$.
\end{enumerate}
\end{lem}
\dim\ \
The claim 1.~is easy. We prove 2.~by induction with respect to the number of elements
of $J$. Assume $\#J > 1$ and fix an element $k \in J$.
The set $S_k$ is
the multi-cone $S(U \cap Z_k, \{G_j\}_{j \in J^*}, \epsilon)$
in $Z_k$.
Here $J^*$ was given in (\ref{eq:subsector-set-def}).
Noticing that $Z_k \subset Z_j$ if $\hat{I}_j \cap I_{k} \ne \emptyset$,
we have
$$
S_{J} = (S_k)_{J^* \cap J}.
$$
Then, by induction hypothesis, we can find $z_0 \in S_k$ satisfying
$$
\tilde{z} = \mu'_{J^* \cap J}(z_0,(0,\dots,0))
$$
where $\mu'_{J^* \cap J}$
is the multi action in $Z_k$ defined by $\{I^*_j\}_{j \in J^*}$.
We can also find $z_1 \in S$ with $\mu_k(z_1,0) = z_0$.
As the restriction of
$\mu_{J^* \cap J}$ to $Z_k$ coincides with $\mu'_{J^* \cap J}$, we have
$$
\tilde{z} = \mu_{(J^* \cap J) \cup \{k\}}(z_1,(0,\dots,0)).
$$
Since $\mu_j(z,0)$ ($j \in J$) is the identity map on $Z_J$, we finally obtain
$\tilde{z} = \mu_{J}(z_1,(0,\dots,0))$, which completes the proof.
\qed

Thanks to the lemma, we have
\begin{cor}{\label{cor:multi-cone-domain}}
Let $J \in \Powl$ and $S := S(U, \{G_j\}, \epsilon)$.
We have $\pi_J(S) = S_J$ where $\pi_J: X \to Z_J$ is the
canonical projection defined by $(z) \mapsto (z_{I_J^{\operatorname{C}}})$
$($$I_J^{\operatorname{C}}$ is the complement set of $I_J$$)$.
\end{cor}
\dim\ \
As $\pi_J(z) = \mu_J(z,(0,\dots,0))$ holds, the inclusion
$S_J \subset \pi_J(S)$ comes from 2.~of Lemma \ref{lemma:mutli-cone-geometric}.
Let us show a converse inclusion. Let $\tilde{z} \in \pi_J(S)$ with
$\tilde{z} = \pi_J(z')$ ($z' \in S$). Then
$z'_\lambda := \mu_J(z', (\lambda, \lambda, \dots, \lambda))$
$(0 < \lambda \le 1)$ is contained in $S$, and
$$
\lim_{\lambda \to 0^+} z'_\lambda = \mu_J(z',(0,\dots,0)) = \pi_J(z') = \tilde{z}.
$$
Moreover, as $\mu_J(z,(0,\dots,0)): X \to Z_J$ is an open map,
an open neighborhood of $z'$ contained in $S$ is mapped to
an open neighborhood of $\tilde{z}$ in $Z_J$ by $\mu_J(z,(0,\dots,0))$.
This implies $\tilde{z} \in S_{J}$.
\qed

\begin{lem}\label{lemma:1-regular-S}
The multi-cone $S := S(U,\{G_j\}, \epsilon)$ is 1-regular, that is,
there exists a positive constant $C$ such that,
for any points $z$ and $w$ in $S$, there exists a continuous subanalytic curve $s$ in $S$ which joins
$z$ and $w$ and which satisfies
$$
(\text{the length of $s$})  \le C|z - w|.
$$
\end{lem}
\dim\ \
We may assume
$
X = \mathbb{C}^{\hat{I}_1} \times \dots \times {\mathbb{C}}^{\hat{I}_\ell}
\times \mathbb{C}^d
$
with $d = n - \displaystyle\sum_{j=1}^{\ell} \# \hat{I}_j$ and
$$
z = (z_{\hat{I}_1}, \dots, z_{\hat{I}_\ell}, z'),\
w = (w_{\hat{I}_1}, \dots, w_{\hat{I}_\ell}, w') \in S
$$
where $z_{\hat{I}_j}$ (resp. $w_{\hat{I}_j}$) is a point in $\mathbb{C}^{\hat{I}_j}$
and $z'$ (resp. $w'$) is one in $\mathbb{C}^d$. Set
$$
r_j := \min\{\vert z_{\hat{I}_j}\vert,\, \vert w_{\hat{I}_j} \vert \} > 0.
$$
We define the points $\tilde{z}$, $\tilde{w}$ and $\tilde{\tilde{w}}$ by
$$
\begin{aligned}
\tilde{z} &= \left(
\frac{r_1}{\vert z_{\hat{I}_1}\vert} z_{\hat{I}_1}, \dots,
\frac{r_\ell}{\vert z_{\hat{I}_\ell}\vert} z_{\hat{I}_\ell}, z'\right), \qquad
\tilde{w} = \left(
\frac{r_1}{\vert w_{\hat{I}_1}\vert} w_{\hat{I}_1}, \dots,
\frac{r_\ell}{\vert w_{\hat{I}_\ell}\vert} w_{\hat{I}_\ell}, z'\right) \\
\tilde{\tilde{w}} &= \left(
w_{\hat{I}_1}, \dots, w_{\hat{I}_\ell}, z'\right)
\end{aligned}
$$
respectively. Note that these points belong to $S$ by definition of $S$.

We will construct a path with the desired property by joining these points.
The points $z$ and $\tilde{z}$,
$\tilde{w}$ and $\tilde{\tilde{w}}$,
$\tilde{\tilde{w}}$ and $w$ are connected by
the straight paths.  We can easily see that they are also contained in $S$.

Let $H_j$ be the sphere in $\mathbb{C}^{\hat{I}_j}$
with radius $r_j$ and center at the origin. Then
the points $\tilde{z}$ and $\tilde{w}$ are in
$$
H := H_1 \times \dots \times H_{\ell} \times \{z'\} \subset
 \mathbb{C}^{\hat{I}_1} \times \dots \times {\mathbb{C}}^{\hat{I}_\ell}
 \times \mathbb{C}^d = X.
$$
Then we join the points $\tilde{z}$ and $\tilde{w}$ by the shortest path
in $H$. This path also belongs to $S$.
Now the points $w$ and $z$ are joined by the paths constructed above,
for which the estimate of the lemma holds
with the constant $C = \sqrt{\ell}(1 + \pi) + 1$.
This completes the proof.
\qed

Now we introduce a formal power series appearing in our multi-asymptotic expansions.
Let $J = \{j_1, \dots, j_m\} \in \Powl$, and let
$\lambda = (\lambda_1, \dots, \lambda_\ell)$ be
variables of parameters for the multi-action
$\mu:=\mu_{\{1,2,\dots,\ell\}}$ in $X$.
We denote by $\Z_{\ge 0}^{J}$ the subset of $\Z_{\ge 0}^\ell$
consisting of $\beta = (\beta_1, \beta_2, \dots, \beta_\ell) \in \Z_{\ge 0}^\ell$ with
$\beta_j = 0$ if $j \notin J$. Note that $\Z_{\ge 0}^J$ is isomorphic
to $\Z_{\ge 0}^{\#J}$.
We set
$\partial_{\lambda}^\beta :=
\displaystyle\frac{\partial^\beta}{\partial \lambda^\beta}$
for $\beta \in \Z_{\ge 0}^\ell$. When $\beta$ belongs to $\Z_{\ge 0}^J$,
we sometimes denote it by $\partial_{\lambda, J}^\beta$ to emphasize the fact
that $\beta$ is an element of $Z_{\ge 0}^J$.

For $\beta \in \Z^J_{\ge 0}$, we introduce the polynomial of
the variables $z_{I_J}$ with constant coefficients by
\begin{equation}
T^\beta_{J}(z_{I_J}) :=
\left.\frac{1}{\beta!}\exp(-\mu(z,\lambda))
\partial^\beta_{\lambda, J}
\exp(\mu(z,\lambda))\right|_{\lambda = e_{\lambda,J}}.
\end{equation}
Here $e_{\lambda, J} \in \CC^\ell$ is the point $(\lambda_1, \dots, \lambda_\ell)$
with $\lambda_j = 0$ ($j \in J$) and $\lambda_j = 1$ ($j \notin J$) and
the exponential function $\exp$ on $\CC^n$ is defined by
$\exp(z_1,\dots,z_n):= e^{z_1} \dots e^{z_n}$ as usual.

Then, for $N = (n_1,n_2, \dots, n_\ell) \in \Z_{\ge 0}^\ell$, we define
a polynomial of $z_{I_J}$ with constant coefficients by
\begin{equation}
T^{<N}_{J}(z_{I_J}) := \sum_{\beta <_J N,\, \beta \in \Z_{\ge 0}^J} T^\beta_{J}(z_{I_J})
\end{equation}
where $\beta <_J N$ if and only if $\beta_k < n_k$ for any $k \in J$.
Note that if there exists no index $\beta <_J N$, then we set
$T^{<N}_{J}(z_{I_J}) := 0$ as usual convention.

We now give a concrete form of the polynomial $T^{<N}_{J}(z_{I_J})$.
Let $I=\{i_1, \dots,i_m\}$ be a subset of $\{1,2,\dots,n\}$.
We denote by $\Z^I_{\ge 0}$ the subset of $\Z^n_{\ge 0}$ whose
element is $(\alpha_1, \dots, \alpha_n)$ with $\alpha_i = 0$ if $i \notin I$.
For $\alpha \in \Z_{\ge 0}^l$,
a monomial $z^\alpha$ designates $z_1^{\alpha_1}\dots z_n^{\alpha_n}$. If
$\alpha$ belongs to $\Z^I_{\ge 0}$, we also write it by $z_I^\alpha$.
We set
$$
\vert \alpha \vert_{I} = \alpha_{i_1} + \dots + \alpha_{i_m}.
$$
\begin{lem}{\label{lem:T_N-polynomial}}
We have
\begin{equation}{\label{eq:lemma-polynomial-TN}}
T^{<N}_{J}(z_{I_J}) = \sum_{\alpha \in A(N)} \frac{1}{\alpha!}z_{I_J}^\alpha
\end{equation}
where the indices set $A(N)$ is given by
$$
A(N) := \{\alpha \in \Z^{I_J}_{\ge 0};\, \vert \alpha \vert_{I_j} < n_j
\text { for any $j \in J$}\}.
$$
\end{lem}
\dim \ \
Let $i \in I_J$ and
$\beta = (\beta_1, \dots, \beta_n) \in Z^J_{\ge 0}$
($\beta_j = 0$ if $j \notin J$). We first assume
$\beta_j \ge 1$ for any $j \in J$. Then we have
$$
\begin{aligned}
\frac{\partial}{\partial {z_i}}T^{\beta}_J
& =
\frac{\partial}{\partial {z_i}}
\left( \left.\frac{1}{\beta!}\exp(-\mu(z,\lambda))
\partial^\beta_{\lambda, J}
\exp(\mu(z,\lambda))\right|_{\lambda = e_{\lambda,J}}\right) \\
&=
\left.\frac{1}{\beta!}\exp(-\mu(z,\lambda))
\partial^\beta_{\lambda, J}
\left(
\underset{j \in J_i}{\prod}\lambda_j
\exp(\mu(z,\lambda))\right)\right|_{\lambda = e_{\lambda,J}} \\
&=
\left.\frac{1}{(\beta - e_{J_i})!}\exp(-\mu(z,\lambda))
\partial^{\beta - e_{J_i}}_{\lambda, J}
\exp(\mu(z,\lambda))\right|_{\lambda = e_{\lambda,J}}
= T^{\beta-e_{J_i}}_J
\end{aligned}
$$
where $e_{J_i} = (e_{1,\, J_i}, \dots, e_{\ell,\, J_i}) \in \Z^l_{\ge 0}$
is determined
by $e_{j,\, J_i} = 1$ for $j \in J_i$ and $e_{j,\, J_i} = 0$ for $j \notin J_i$
(the definition of $J_i$ is given by (\ref{eq:def-J_i})).
This formula also holds for any $\beta \in \Z^l_{\ge 0}$
if we set $T^\beta_J(z_{I_J}) := 0$
for $\beta = (\beta_1, \dots, \beta_l) \in \Z^J$
with some $\beta_j < 0$ ($j \in J$).
Hence we obtained, for $N=(n_1, \dots, n_\ell) \in \Z_{\ge 0}^l$,
$$
\frac{\partial}{\partial {z_i}} T^{<N}_J = T^{<N - e_{J_i}}_J.
$$

We prove the lemma by induction with $n = n_1 + \dots + n_\ell$.
If $n = 0$, as both sides of the equation
(\ref{eq:lemma-polynomial-TN}) are zero by definition,
the lemma is true.

Now we prove the lemma for a general $n > 0$.
Let us consider the system of partial differential equations
of an unknown function
$u(z_{I_J})$ defined by
\begin{equation}{\label{eq:partial-diff}}
\frac{\partial}{\partial {z_i}} u = T^{<N - e_{J_i}}_J \quad (i \in I_J)
\end{equation}
Then, by induction hypothesis, the right hand side of the above equation
is given by
$
\displaystyle\sum_{\alpha \in A(N-e_{J_i})} \frac{1}{\alpha!}z_{I_J}^\alpha.
$
Clearly both $u=T^{<N}_{J}(z_{I_J})$ and
$
u=\displaystyle\sum_{\alpha \in A(N)} \frac{1}{\alpha!}z_{I_J}^\alpha
$
satisfy the same equation (\ref{eq:partial-diff}).
The solution of (\ref{eq:partial-diff}) is uniquely determined
if the initial values at $z_{I_J} = 0$ is given.
It is easy to see that $T^{<N}_{J}(0) = 1$ if $A(N) \ne \emptyset$
and $T^{<N}_{J}(0) = 0$ if $A(N) = \emptyset$.
Hence we have obtained
(\ref{eq:lemma-polynomial-TN}) for $N$.
This completes the proof.
\qed

\begin{df}
Let $S:=(U, \{G_i\},\epsilon)$ be a multi-cone in $X$.
We say that $F:=\{F_{J}\}_{J \in \Powl}$
is a total family of coefficients of multi-asymptotic expansion
along ${\underset{1 \le j \le \ell}{\bigcup} Z_j}$ on $S$ if each $F_J$ consists of
a family of holomorphic functions
$\{f_{J, \alpha}(z_{I^{\operatorname{C}}_J})\}_{\alpha \in \Z_{\ge 0}^{I_J}}$
defined on $S_J$. Here $I^{\operatorname{C}}_J$ is the complement set
of $I_{J}$ $($note that $z_{I^{\operatorname{C}}_J}$ are the coordinates of
the submanifold $Z_J$$)$.
\end{df}

Let $F$ be a total family of coefficients defined in the above definition.
We introduce a map $\tau_{F, J}$ from polynomials of the variables $z_{I_J}$
to those with coefficients in holomorphic functions on $S_J$ in the following
way. Let $p(z_{I_J}) = \displaystyle\sum_\alpha c_\alpha z_{I_J}^\alpha$
be a polynomial of the variable $z_{I_J}$ with constant coefficients. Then
we define $\tau_{F, J}(p)(z)$ by replacing
a monomial $z_{I_J}^\alpha$ in $p(z_{I_J})$ with
$f_{J, \alpha}(z_{I_J^{\operatorname{C}}})z^\alpha_{I_J}$ where $f_{J, \alpha}$
is given in $F_J = \{f_{J, \alpha}\}$, that is,
\begin{equation}
\tau_{F, J}(p)(z) := \sum_\alpha c_\alpha f_{J,\alpha}(z_{I_J^{\operatorname{C}}})
z_{I_J}^\alpha.
\end{equation}
Note that $\tau_{F, J}(p)(z)$ is a holomorphic function on
$\pi_J^{-1}(S_J)$, in particular,
it is defined on $S \subset \pi_J^{-1}(S_J)$
by Corollary \ref{cor:multi-cone-domain}.
We set
\begin{equation}
\begin{aligned}
T^{\beta}_{J}(F;\, z) &= \tau_{F, J}(T^\beta_J),\\
T^{<N}_{J}(F;\, z) &:= \sum_{\beta <_J N} T^\beta_{J}(F;\, z),
\end{aligned}
\end{equation}
and
\begin{equation}
\begin{aligned}
\operatorname{App}^{<N}(F; z)
:= \sum_{J \in \Powl}
(-1)^{(\#J + 1)} T^{<N}_{J}(F;\, z).
\end{aligned}
\end{equation}
It follows from Lemma \ref{lem:T_N-polynomial} that,
if ${\underset{1 \le j \le \ell}{\bigcup} Z_j}$
forms a normal crossing divisor,
$\operatorname{App}^{<N}(F; z)$ coincides with one defined by
Majima in \cite{Ma84}.

\medskip

Let us recall the definitions of the sets $J_{\supsetneq Z_j}$ and
$J_{\subsetneq Z_j}$ for $j \in \Z_\ell$.
$$
\begin{aligned}
J_{\supsetneq Z_j} &:= \{k \in \Z_\ell;\,  Z_k \supsetneq
Z_j,\,
\text{ there is no $m$ with $Z_k \supsetneq Z_m \supsetneq Z_j$ } \}\\
J_{\subsetneq Z_j} &:= \{k \in \Z_\ell;\,  Z_k \subsetneq Z_j,\,
\text{ there is no $m$ with $Z_k \subsetneq Z_m \subsetneq Z_j$ } \}.
\end{aligned}
$$
Then the function $w_j: \Z^\ell_{\ge 0} \to \Z$ ($j = 1,2, \dots, \ell$)
is defined by
\begin{equation}
w_j(N) = n_j - \sum_{k \in J_{\supsetneq Z_j}} n_k
\end{equation}
for $N=(n_1,n_2, \dots,n_\ell) \in \Z^\ell_{\ge 0}$. Note that $w_j$ takes
not only positive values but also negative ones.

Let $S'$ be another open proper multi-cone $S(U',\{G'_j\}, \epsilon')$.
We say that $S'$ is properly contained in $S$ if
$\epsilon' < \epsilon$, $U'$ is relatively compact in $U$ and
$G'_j$ ($j=1,2, \dots,\ell$) is properly contained in $G_j$ as a conic cone.

\begin{df}
Let $f$ be a holomorphic function on $S = S(U, \{G_j\}, \epsilon)$.
We say that
$f$ is multi-asymptotically developable along
${\underset{1 \le j \le \ell}{\bigcup} Z_j}$  on $S$ if
there exists a total family $F$ of coefficients of
multi-asymptotic expansion along ${\underset{1 \le j \le \ell}{\bigcup} Z_j}$ on $S$
that satisfies the following condition.

For any open proper multi-cone $S'=S(U', \{G'_j\}, \epsilon')$
properly contained in $S$ and for any
$N = (n_1, \dots, n_\ell) \in {\mathbb Z}_{\ge 0}^\ell$,
there exists a constant $C_{S',N} > 0$ for which we have an estimate
\begin{equation}{\label{eq:multi-asymptotic-formula}}
\left\vert f(z) - \operatorname{App}^{<N}(F; z) \right\vert \le C_{S', N}
\underset{1 \le j \le \ell}{\prod}
\operatorname{dist}(z,\, Z_j)^{w_j(N)}
\qquad (z \in S').
\end{equation}
\end{df}

Let us see some typical examples.

\begin{es} \label{Example:Majima}
Let $X=\mathbb C^n$ with the coordinates
$(z_1, z_2, \dots, z_n) = (z_1, z_2, z')$ and
$Z_j$ $(j=1,2)$ submanifolds defined by $\{z_j = 0\}$ $($i.e. $I_j = \{j\}$$)$.
Let $G_j$ $(j=1,2)$ be a proper open sector in $\CC$ and $U_R:=B^1_R \times B^1_R
\times B^{n-2}_R$ where $B^k_R$ designates an open ball in $\CC^k$
with its radius $R > 0$ and center $0$.
In this case, the multi-cone $S=S(U_R,\{G_1, G_2\}, \epsilon)$ is nothing but
a poly-sector $(G_1 \times G_2 \times \CC^{n-2}) \cap U_R$ and
we have
$$
S_1 = (G_2 \cap B^1_R) \times B^{n-2}_R,\,
S_2 = (G_1 \cap B^1_R) \times B^{n-2}_R,\,
S_{\{1,2\}} = B^{n-2}_R.
$$
Let $F$ be a total family of coefficients of multi-asymptotic expansion, that is,
$$
\begin{aligned}
F &= (F_{\{1\}}, F_{\{2\}}, F_{\{1,2\}}) \\
&=\left(
\left\{f_{\{1\},k}(z_2, z')\right\}_{k \ge 0},\,
\left\{f_{\{2\},k}(z_1, z')\right\}_{k \ge 0},\,
\left\{f_{\{1,2\},\alpha}(z')\right\}_{\alpha \in \Z_{\ge 0}^2}\right)
\end{aligned}
$$
where $f_{\{1\},k}$ $($resp. $f_{\{2\},k}$ and $f_{\{1,2,\}\alpha}$$)$ is
a holomorphic function on $S_1$ $($resp. $S_2$ and $S_{\{1,2\}}$$)$.
For $N=(n_1, n_2) \in \Z^2_{\ge 0}$, an asymptotic expansion
$\operatorname{App}^{<N}(F; z)$ is given by
$$
\begin{aligned}
	T_{\{1\}}^{<N}(F;\, z) &= \sum_{k < n_1} f_{\{1\},k}(z_2, z')\frac{z_1^k}{k!}, \\
	T_{\{2\}}^{<N}(F;\, z) &= \sum_{k < n_2} f_{\{2\},k}(z_1, z')\frac{z_2^k}{k!}, \\
	T_{\{1,2\}}^{<N}(F;\, z) &= \sum_{\alpha_1 < n_1, \alpha_2 < n_2}
	f_{\{1,2\},\alpha}(z') \frac{z_1^{\alpha_1}z_2^{\alpha_2}}{\alpha_1! \alpha_2!} \\
	\operatorname{App}^{<N}(F;\, z) &=
T_{\{1\}}^{<N}(F;\, z) + T_{\{2\}}^{<N}(F;\, z) - T_{\{1,2\}}^{<N}(F;\, z).
\end{aligned}
$$
As $w_j(N) = n_j$ and $\operatorname{dist}(z, Z_j) = \vert z_j \vert$ $(j=1,2)$,
a holomorphic function $f$ is multi-asymptotically developable to $F$ if,
for any poly-sector $S'$ properly contained in $S$ and
for any $N=(n_1, n_2) \in \Z_{\ge 0}^2$,
there exists a positive constant $C_{S',N}$ such that
$$
\left|
f(z) - \operatorname{App}^{<N}(F;\, z)
\right|
\le C_{S', N} \vert z_1 \vert^{n_1} \vert z_2 \vert^{n_2}
\qquad (z \in S').
$$
Hence our definition coincides with that of strongly asymptotic
developability established by Majima in \cite{Ma84}. We give some examples of asymptotics:
\begin{eqnarray*}
N=(0,0) &&
|f(z)| \leq C_{S',N},
\\
N=(1,0) &&
|f(z)-f_{\{1\},0}(z_2,z')| \leq C_{S',N}|z_1|,
\\
N=(0,1) &&
|f(z)-f_{\{2\},0}(z_1,z')| \leq C_{S',N}|z_2|,
\\
N=(1,1) &&
|f(z)-f_{\{1\},0}(z_2,z')-f_{\{2\},0}(z_1,z')+f_{\{1,2\},(0,0)}(z')| \leq C_{S',N}|z_1||z_2|.
\end{eqnarray*}

\end{es}

\begin{es} \label{Example:Takeuchi}
Let $X=\mathbb C^n$ with the coordinates
$(z_1, z_2, \dots, z_n) = (z_1, z_2, z')$, and let
$Z_1 = \{z_1 = 0\}$ and $Z_2 = \{z_1 = z_2 = 0\}$.
In this case,  $I_1$ and $I_2$ are given by
$\{1\}$ and $\{1,2\}$ respectively, and we have
$\hat{I}_1 = \{1\}$, $\hat{I}_2 = \{2\}$.
Let $G_1$ and $G_2$ be proper open sectors in $\CC$.
We set $U_R:=B^1_R \times B^1_R \times B^{n-2}_R$.
The multi-cone $S=S(U_R,\{G_1, G_2\}, \epsilon)$ is defined by
$$
\{(z_1, z_2, z') \in (G_1 \times G_2 \times \CC^{n-2}_{z'}) \cap U;\,
\vert z_1 \vert < \epsilon \vert z_2 \vert \}.
$$
Then we have
$$
S_1 = (G_2 \cap B^1_R) \times B^{n-2}_R,\,
S_2 = S_{\{1,2\}} = B^{n-2}_R.
$$
Let $F$ be a total family of coefficients of multi-asymptotic expansion, which consists of
$$
\begin{aligned}
F &= (F_{\{1\}}, F_{\{2\}}, F_{\{1,2\}}) \\
&=\left(
\left\{f_{\{1\},k}(z_2, z')\right\}_{k \ge 0},\,
\left\{f_{\{2\},\alpha}(z')\right\}_{\alpha \in \Z^2_{\ge 0}},\,
\left\{f_{\{1,2\},\alpha}(z')\right\}_{\alpha \in \Z^2_{\ge 0}}\right)
\end{aligned}
$$
where $f_{\{1\},k}$ $($resp. $f_{\{2\},\alpha}$ and $f_{\{1,2\},\alpha}$$)$ is
a holomorphic function on $S_1$ $($resp. $S_2$ and $S_{\{1,2\}}$$)$.
Then we have, for $N=(n_1, n_2) \in \Z^2_{\ge 0}$,
$$
\begin{aligned}
	T_{\{1\}}^{<N}(F;\, z) &= \sum_{k < n_1} f_{\{1\},k}(z_2, z')\frac{z_1^k}{k!}, \\
	T_{\{2\}}^{<N}(F;\, z) &= \sum_{\alpha_1 + \alpha_2 < n_2}
	 f_{\{2\},\alpha}(z')\frac{z_1^{\alpha_1}z_2^{\alpha_2}}{\alpha_1!\alpha_2!}, \\
	T_{\{1,2\}}^{<N}(F;\, z) &= \sum_{\alpha_1 < n_1,\, \alpha_1 + \alpha_2 < n_2}
	f_{\{1,2\},\alpha}(z') \frac{z_1^{\alpha_1}z_2^{\alpha_2}}{\alpha_1! \alpha_2!}, \\
	\operatorname{App}^{<N}(F;\, z) &=
T_{\{1\}}^{<N}(F;\, z) + T_{\{2\}}^{<N}(F;\, z) - T_{\{1,2\}}^{<N}(F;\, z).
\end{aligned}
$$
It follows from the definition that we have
$$
w_1(N) = n_1 \text{ and } w_2(N) = n_2 - n_1
$$
and
$$
\operatorname{dist}(z, Z_1) = \vert z_1 \vert
\text{ and }
\operatorname{dist}(z, Z_2) \simeq \vert z_1 \vert + \vert z_2 \vert.
$$
Hence a holomorphic function $f$ is multi-asymptotically developable to $F$ if,
for any proper multi-cone $S'$ properly contained in $S$ and
for any $N=(n_1, n_2) \in \Z_{\ge 0}^2$,
there exists a positive constant $C_{S',N}$ such that
$$
\left|
f(z) - \operatorname{App}^{<N}(F;\, z)
\right|
\le C_{S', N} \vert z_1 \vert^{n_1}
(\vert z_1 \vert + \vert z_2 \vert)^{n_2 - n_1}
\qquad (z \in S').
$$
We give some examples of asymptotics:
\begin{eqnarray*}
N=(0,0) &&
|f(z)| \leq C_{S',N},
\\
N=(1,0) &&
|f(z)-f_{\{1\},0}(z_2,z')| \leq C_{S',N}{|z_1| \over |z_1|+|z_2|},
\\
N=(0,1) &&
|f(z)-f_{\{2\},0}(z')| \leq C_{S',N}(|z_1|+|z_2|),
\\
N=(1,1) &&
|f(z)-f_{\{1\},0}(z_2,z')-f_{\{2\},0}(z')+f_{\{1,2\},(0,0)}(z')| \leq C_{S',N}|z_1|.
\end{eqnarray*}

\end{es}

One of important features of multi-asymptoticity is stability for differentiations.
\begin{prop}{\label{prop:multi-asymptotic-derivative}}
Let $S:=(U, \{G_i\},\epsilon)$ be a proper multi-cone in $X$ and $f$ a holomorphic function
on $S$. If $f(z)$ is multi-asymptotically developable along
${\underset{1 \le j \le \ell}{\bigcup} Z_j}$ on
$S$, then any derivative of $f$ is also multi-asymptotically developable
along ${\underset{1 \le j \le \ell}{\bigcup} Z_j}$ on $S$.
\end{prop}
\dim\ \
It suffices to show that, for an $i \in \{1,2,\dots, n\}$,
$\displaystyle\frac{\partial f}{\partial z_i}$
is also multi-asymptotically developable. We need the lemma below.
Let $F:=\{F_{J}\}_{J \in \Powl}$
be a total family of coefficients of multi-asymptotic expansion with
$F_J = \{f_{J,\alpha}\}_{\alpha \in \Z_{\ge 0}^{I_J}}$.
\begin{lem}{\label{lemma:multi-asymptotic-derivative}}
Let $N=(n_1, n_2, \dots, n_\ell) \in \Z_{\ge 0}^\ell$.  Then we have
$$
\frac{\partial} {\partial z_i}
\operatorname{App}^{<N_+}(F;\, z)
= \operatorname{App}^{<N}(F';\, z).
$$
Here $N_+=(n'_1,\dots,n'_\ell)$ is determined by $n'_j = n_j + 1$ for $j$ with $i \in I_j$
and $n'_j = n_j$ otherwise. The total family $F'=\{F'_J\}$
is given by
$$
F'_J := \{f'_{J,\alpha}\} =
\left\{\displaystyle\frac{\partial f_{J,\alpha}}{\partial z_i}\right\}
\qquad (\text{for $J \in \Powl$ with $i \notin I_J$})
$$
and
$$
F'_J := \{f'_{J,\alpha}\} := \{f_{J, \alpha + e_{i}}\}
\qquad(\text{for $J \in \Powl$ with $i \in I_J$})
$$
where $e_i \in \Z_{\ge 0}^n$ is the unit vector whose $i$-th element is equal to $1$.
\end{lem}
\dim\ \
Let $J = (j_1, j_2, \dots, j_m) \in \Powl$.
For $J$ with $i \notin I_J$, as the polynomial $T_J^{<N}$ does not contain
the variable $z_i$ and $n'_j = n_j$ for $j \in J$, we have
$$
\frac{\partial}{\partial {z_i}} T^{<N_+}_J(F;\, z)
=
\frac{\partial}{\partial {z_i}} \tau_{F, J}\left(T^{<N_+}_J\right)
=
\tau_{F', J}\left(T^{<N_+}_J\right)
=
T^{<N}_J(F';\, z).
$$

Assume that $J$ satisfies $i \in I_J$.
Then, by the proof of Lemma \ref{lem:T_N-polynomial}, we have
$$
\frac{\partial}{\partial {z_i}} T^{<N_+}_J = T^{<N}_J.
$$
As
$$
\displaystyle\frac{\partial}{\partial z_i}\tau_{F,J}(p) =
\tau_{F',J}\left(\displaystyle\frac{\partial p}{\partial z_i}\right)
$$ holds for a polynomial $p(z_{I_J})$ of the variables $z_{I_J}$,
we have obtained
$$
\begin{aligned}
\frac{\partial}{\partial {z_i}} T^{<N_+}_J(F;\, z) &=
\frac{\partial}{\partial {z_i}} \tau_{F, J}(T^{<N_+}_J) =
\tau_{F', J}\left(\frac{\partial}{\partial {z_i}} T^{<N_+}_J\right) \\
&= \tau_{F', J}\left(T^{<N}_J\right)
= T^{<N}_J(F';\, z).
\end{aligned}
$$
\qed

We continue the proof of the proposition.
It suffices to consider the case for $i \in I_{\{1,2,\dots,\ell\}}$.
By the lemma, we have
$$
\left|
\frac{\partial f}{\partial z_i} - \operatorname{App}^{<N}(F';\, z) \right|
=
\left| \frac{\partial}{\partial z_i}
\left( f(z) - \operatorname{App}^{<N_+}(F;\, z) \right) \right|.
$$
Let $j_0$ be a unique integer with $i \in \hat{I}_{j_0}$ and
$S':=S(U', \{G'_j\}, \epsilon')$
(resp. $S'':=S(U'', \{G''_j\}, \epsilon'')$)
a proper multi-cone properly
contained in $S$ (resp. $S'$). Then, as
$\operatorname{dist}(z, Z_{j_0}) \le \epsilon''
\operatorname{dist}(z, Z_j)$ ($z \in S''$) holds
for $Z_j \subset Z_{j_0}$,
there exists a positive constant $\kappa > 0$
such that for any point $z^* \in S''$
\begin{equation}{\label{diff-caucy-domain}}
z^* + \left\{z \in \CC^n;\, z_k = 0\, (k \ne i), \,
\vert z_i \vert \le \kappa \min\{1,\, \operatorname{dist}(z^*, Z_{j_0})\}
\right\}
\subset S'.
\end{equation}
Assume that $z \in S''$. Set
$$
D = \left\{(\zeta_1,\dots,\zeta_n);\, \zeta_\alpha = z_\alpha\, (\alpha \neq i),\, |\zeta_i - z_i| = {\kappa \over 2}\operatorname{dist}(z,Z_{j_0})\right\}.
$$
(Here we may assume $\kappa < 1$). $D$ is a circle in the $z_i$-plane and the other coordinates are fixed. By \eqref{diff-caucy-domain}, $D$ is contained in $S'$.

Now remark the following fact. Let $j \in \{1, 2, \dots , \ell\}$. Then
\begin{equation}
\operatorname{dist}(z,Z_j) \simeq \sum_{i \in \hat{I}_j}|z_i|
\end{equation}
(note that sum is taken over indices in $\hat{I}_j$ and not in $I_j$) as long as $z \in S'$, which
comes from the fact that each cone $G_j$ is proper. Hence, for $j \neq j_0$, we may assume
$$
\operatorname{dist}(z,Z_j) = \operatorname{dist}(\zeta,Z_j) \qquad (\zeta \in D)
$$
because $\hat{I}_j$ does not contain the index $i$ and $D$ is a circle in the $z_i$-plane (the other coordinates are fixed). Further, for $j = j_0$, we have
$$
\left(1 - {\kappa \over 2}\right) \operatorname{dist}(z,Z_j) \leq \operatorname{dist}(\zeta,Z_j) \leq \left(1 + {\kappa \over 2}\right) \operatorname{dist}(z,Z_j) \qquad (\zeta \in D)
$$
(note that both inequalities are needed in the estimation \eqref{eq:estimation h} below, which
depends on the sign of $w_{j_0}(N)$).

Since $f(z)$ is multi-asymptotically developable, there exists a constant
$C_{S', N_+}$ satisfying
$$
\left| h(z):=f(z) - \operatorname{App}^{<N_+}(F;\, z) \right|
\le C_{S', N_+}
\underset{1 \le j \le \ell}{\prod}
\operatorname{dist}(z,\, Z_j)^{w_j(N_+)} \qquad (z \in S').
$$
Therefore we obtain
\begin{equation}\label{eq:estimation h}
|h(\zeta_1,\dots,\zeta_i,\dots,\zeta_n)| \leq C_\kappa\prod_{1\leq j\leq \ell}\operatorname{dist}(z,Z_j)^{w_j (N+)} \qquad (\zeta \in D)
\end{equation}
for some positive constant $C_\kappa$. By the Cauchy integral formula, we have
$$
{\partial h \over \partial z_i}(z) =
\frac{1}{2\pi\sqrt{-1}}
\int_{\partial D} {h(z_1,\dots,\zeta_i,\dots,z_n) \over
(\zeta_i - z_i)^2} d\zeta_i,
$$
where the path of the integration is the projection of $D$ to the $z_i$-plane. Putting
the estimation of $h(z)$ and the above estimations into the formula,  we get
$$
\left| \frac{\partial h(z)}{\partial z_i} \right| \le
C\frac{
\underset{1 \le j \le \ell}{\prod}
\operatorname{dist}(z,\, Z_j)^{w_j(N_+)}}
{\operatorname{dist}(z, Z_{j_0})} \qquad (z \in S'')
$$
for a constant $C > 0$.
As $w_{j_0}(N) = w_{j_0}(N_+) - 1$ and $w_j(N) = w_j(N_+)$ for $j \ne j_0$,
we finally obtain
$$
\left|
\frac{\partial f}{\partial z_i} - \operatorname{App}^{<N}(F';\, z) \right|
\le C
\underset{1 \le j \le \ell}{\prod} \operatorname{dist}(z,\, Z_j)^{w_j(N)}
\qquad (z \in S'').
$$
This completes the proof.
\qed

By this proposition and integration by parts, we can obtain the following theorem.
\begin{teo}{\label{theorem:multi-asymptotic}}
Let $S:=S(U, \{G_i\}, \epsilon)$ be a proper multi-cone and $f$ a holomorphic
function on $S$. Then the following conditions are equivalent.
\begin{enumerate}
\item $f$ is multi-asymptotically developable along ${\underset{1 \le j \le \ell}{\bigcup} Z_j}$
on $S$.
\item  For any open proper multi-cone $S':=S(U', \{G'_i\}, \epsilon')$ properly
contained in $S$ and for any $\alpha \in \Z_{\ge 0}^n$,
$\left|\displaystyle\frac{\partial^\alpha f}{\partial z_\alpha}\right|$ is
bounded on $S'$.
\item For any open proper multi-cone $S':=S(U', \{G'_i\}, \epsilon')$ properly
contained in $S$,
the holomorphic function $f\vert_{S'}$ on $S'$ can be extended to a $\C^\infty$-function
on $X_{\mathbb R}$ ($X_{\mathbb R}$ denotes the underlying real analytic manifold of $X$).
\end{enumerate}
\end{teo}
\dim\ \
We first show 1.~implies 2.
By (\ref{eq:multi-asymptotic-formula}) with $N=(0,\dots,0)$, we
obtain that $f$ is bounded on $S'$. Since each higher derivative of $f$
is still multi-asymptotically developable thanks to
Proposition \ref{prop:multi-asymptotic-derivative},
$\displaystyle\frac{\partial^\alpha f}{\partial z^\alpha}$
is also bounded for any
$\alpha \in \Z^n_{\ge 0}$ on $S'$.

\medskip
As $S'$ is 1-regular by Lemma \ref{lemma:1-regular-S}
and as $f$ is holomorphic (i.e.
$\displaystyle\frac{\partial^\alpha f}{\partial \bar{z}^\alpha} = 0$),
the claim 3.~follows from 2.~by the result of Whitney in \cite{Wh34}.
Clearly 3.~implies 2. Hence the claim 2.~and 3.~are equivalent.

\medskip

Now we will show 3.~implies 1.
Assume that $f$ satisfies 3. In particular, any derivative of $f$
extends to $\overline{S'}$ and is bounded on $\overline{S'}$.
It follows from Lemma {\ref{lemma:mutli-cone-geometric}} that
for $z \in S'$ and $0 \le \lambda_j \le 1$ ($j=1,2,\dots,\ell$),
we get $\mu(z,\lambda) \in \overline{S'}$. Therefore,
for $N = (n_1, \dots, n_\ell) \in \Z^\ell_{\ge 0}$,
$$
\varphi_N(f;z) :=
\int_{0}^{\infty}\dots\int_{0}^{\infty}
\frac{\partial^N f(\mu(z,\lambda))}{\partial \lambda^N}
\underset{1 \le j \le \ell}{\prod}
K_{n_j - 1}(1 - \lambda_j) d\lambda_1 \dots d\lambda_\ell
$$
where
$$
K_n(t) :=
\left\{
\begin{aligned}
\delta(t)\qquad &(n = -1), \\
\frac{t_+^n}{n!}      \qquad &(n \ge 0)
\end{aligned}
\right.
$$
is well-defined on $S'$. Here $\delta$ denotes the Dirac delta function and $t_+= t$ if $t \ge 0$, $t_+=0$ if $t < 0$.
We define
$e_J = (e_{1,J}, \dots, e_{\ell,J}) \in \Z^\ell$
by $e_{j,J} = 0$ if $j \in J$
and $e_{j,J} = 1$ otherwise. Then, by integration by parts, we have
$$
\varphi_N(f;z) = f(z) - \sum_{J \in \Powl}
(-1)^{\#J + 1}\sum_{\beta \in \Z^{J}, \beta<_J N}
\left.
\frac{1}{\beta!}\partial^\beta_{\lambda, J}f(\mu(z,\lambda))
\right\vert_{\lambda = e_J}.
$$
Hence it suffices to show  that $\varphi_N(f; z)$ has an estimation
which appears in a multi-asymptotic expansion.
Let us consider coordinates transformation of $\lambda$ as
$\lambda_j = \nu_j(z) \tilde{\lambda}_j$ ($j=1,2,\dots,\ell$) where
$\nu_j(z)$ is given by
$$
\nu_j (z):= \left\{
\begin{aligned}
	\frac{1}{\vert z_{I_j} \vert} \qquad&
	(J_{\subsetneq Z_j} = \emptyset), \\
	\frac{\vert z_{I_k} \vert }{\vert z_{I_j} \vert}
	\qquad&  (k \in J_{\subsetneq Z_j}).
\end{aligned}
\right.
$$
Note that $J_{\subsetneq Z_j}$ consists of at most one element.
Then we have
$$
\frac{\partial^N f(\mu(z,\lambda))}{\partial \lambda^N} =
\frac{1}{\underset{1\le j \le \ell}{\prod} \nu_j^{n_j}(z)}
\frac{\partial^N }{\partial \tilde{\lambda}^N}
f\left(\frac{z_1}{\vert z_{I_{j(1)}} \vert} \tilde{\lambda}_{J_1},
\dots,
\frac{z_n}{\vert z_{I_{j(n)}} \vert}  \tilde{\lambda}_{J_n}\right).
$$
Here, for $k \in \{1,2,\dots,n\}$, we denote by $j(k)$ the integer $j$
that satisfies $k \in \hat{I}_{j}$ and we set
$J_k = \{j \in \Z_\ell;\, k \in I_j\}$.
Since $\displaystyle{\frac{1}{\nu_j(z)}}$ is bounded on $S'$,
$\tilde{\lambda}_j$ is also bounded when
$0 \le \lambda_j \le 1$. Hence we get
$$
\sup_{\lambda \in [0,1]^\ell }
\left| \frac{\partial^N f(\mu(z,\lambda))}{\partial \lambda^N} \right| \le
\frac{C}{\underset{1\le j \le \ell}{\prod} \nu^{n_j}_j(z)} \qquad (z \in S')
$$
for some constant $C > 0$. By noticing
$$
\frac{1}{\underset{1\le j \le \ell}{\prod} \nu_j^{n_j}(z)} =
\underset{1\le j \le \ell}{\prod} \operatorname{dist}(z, Z_j)^{w_j(N)},
$$
we have obtained a multi-asymptotic expansion of $f$ with desired estimation.
The proof has been completed.
\qed

Let $S:=S(U, \{G_i\}, \epsilon)$ be a proper multi-cone
and $f$ a holomorphic function that is multi-asymptotically developable
along ${\underset{1 \le j \le \ell}{\bigcup} Z_j}$ on $S$.
Let $S^{(k)} := S(U^{(k)}, \{G^{(k)}_{i}\}, \epsilon^{(k)})$ ($k=1,2,3,\dots$)
be a family of
proper multi-cones properly contained in $S$, which satisfies
$$
S^{(1)} \subset S^{(2)} \subset S^{(3)} \subset \dots \text{ and }
\underset{k \ge 1}{\bigcup} S^{(k)} = S.
$$

Let $J \in \Powl$.
Then
$\left.
\displaystyle\frac{\partial^\alpha f}{\partial z^\alpha}
\right|_{S^{(k)}}$
has a unique Whitney extension $\varphi^{(k)}$ on $\overline{S^{(k)}}$
by the above theorem. As $\underset{k \ge 1}{\bigcup} S^{(k)}_J = S_J$,
a family $\left\{\varphi^{(k)}\vert_{S^{(k)}_J}\right\}$ determines
a function on $S_J$,  which is holomorphic.
We denote it by
$\left.
\displaystyle\frac{\partial^\alpha f}{\partial z^\alpha}
\right|_{S_J}.
$
Note that, as $S_J$ is a multi-cone in $Z_J$ for a family of submanifolds
$\chi_J := \{Z_j \cap Z_J\}_{j \in J^*}$
where $J^*$ is given by
$\{j \in \Z_\ell; Z_J \nsubseteq Z_j\}$,
it follows again from the theorem that the holomorphic function
$
\left.
\displaystyle\frac{\partial^\alpha f}{\partial z^\alpha}
\right|_{S_J}
$
is multi-asymptotically developable along
$\underset{j\in J^*}{\bigcup} (Z_j \cap Z_J)$ on $S_J$.

\begin{prop}{\label{prop:consistent-coeffs}}
Let $S:=S(U, \{G_i\}, \epsilon)$ be a proper multi-cone and
$f$ a holomorphic function on $S$. Assume that $f$ is multi-asymptotically
developable to a family
$F:=\{F_J\}_{J \in \Powl}$
of coefficients of multi-asymptotic expansion with
$F_J =\{f_{J,\alpha}\}_{\alpha \in \Z^{I_J}_{\ge 0}}$. Then we have
\begin{equation}{\label{eq:prop-boundary}}
\left.\frac{\partial^\alpha f}{\partial z^\alpha} \right|_{S_J} = f_{J, \alpha}
\qquad (J \in \Powl,\, \alpha \in \Z^{I_J}_{\ge 0}).
\end{equation}
In particular, coefficients of multi-asymptotic expansion are unique and they are
multi-asymptotically developable functions themselves.
\end{prop}
\dim\ \
It suffices to show the proposition for a proper multi-cone
$S':=S(U', \{G'_i\}, \epsilon')$ that is properly contained in $S$.
Moreover, by Proposition \ref{prop:multi-asymptotic-derivative} and
Lemma \ref{lemma:multi-asymptotic-derivative},
it enough to show (\ref{eq:prop-boundary}) with $\alpha = 0$.
We prove the lemma by induction with respect to the number of elements in $J$.
We use the same symbol $f$ for a unique continuous extension of $f$
on $\overline{S'}$ in what follows.

Let $J=\{j\}$.
By taking $N=(n_1,n_2,\dots,n_\ell)$ with $n_j = 1$ and
$n_k = 0$ ($k \ne j$), we have
\begin{equation}{\label{eq:prop-xxx-xxx-xxx}}
\vert f(z) - f_{J,0}(z_{I_J}) \vert \le C
\frac{\operatorname{dist}(z, Z_j)}{\delta_j(z)}
\end{equation}
where $\delta_j(z) = 1$ if
$J_{\subsetneq Z_j}$ is an empty set and
$\delta_j(z) = \operatorname{dist}(z, Z_k)$ if
$J_{\subsetneq Z_j} = \{k\}$. Note that
the set $J_{\subsetneq Z_j}$ consists of at most one element.
For $\tilde{z}^* \in S'_J$, by Lemma \ref{lemma:mutli-cone-geometric},
we can find a point $z^* \in S'$ with $\mu_j(z^*,0) = \tilde{z}^*$.
Then there exists $\delta > 0$
such that $\operatorname{dist}(\mu_j(z^*,[0,1]), Z_k) > \delta$
for $k$ with $Z_k \subsetneq Z_j$ because of $\tilde{z}^* \notin Z_k$
(2.~of Lemma \ref{lemma:fundamental-geo}).
Hence, by putting $z = \mu(z^*,\lambda)$ into (\ref{eq:prop-xxx-xxx-xxx})
and letting $\lambda \to 0^+$, we get $f(\tilde{z}^*) = f_{J,0}(\tilde{z}^*)$.
This shows (\ref{eq:prop-boundary}) with $\alpha = 0$.

Let $J \in \Powl$ with $\#J > 1$,
and let
$J^*$ be a subset of $J$ consisting of $j \in J$
such that $Z_j$ is a minimal submanifold
in $\{Z_k\}_{k \in J}$ with respect to the order $\subset$.
Set $N=(n_1,n_2,\dots,n_\ell)$ with $n_j = 1$ if $j \in J$ and $n_j = 0$
otherwise. Then, by noticing
$\operatorname{dist}(z, Z_k) \le \epsilon' \operatorname{dist}(z, Z_j)$
if $Z_j \subsetneq Z_k$,
we have
\begin{equation}{\label{eq:tmp-tmp-xxx-yyy}}
\begin{aligned}
\Biggl\vert f(z) + \sum_{J' \subsetneq J, J' \ne \emptyset}
&(-1)^{\#J'}f_{J',0}(z_{I_{J'}})
+ (-1)^{\#J}f_{J,0}(z_{I_J}) \Biggr\vert  \\
&\le C
\underset{j \in J^*}{\prod}
\frac{\operatorname{dist}(z, Z_j)}{\delta_j(z)}
\end{aligned}
\end{equation}
By induction hypothesis, for $J' \subsetneq J$, the coefficient
$f_{J',0}$ also has a unique continuous extension on $\overline{S'_{J'}}$ satisfying
$f_{J',0} = f\vert_{\overline{S'_{J'}}}$.
Hence, noticing $S'_J \subset \overline{S'_{J'}} \cap Z_J$
due to 2.~of Lemma \ref{lemma:fundamental-geo},
we have
$$
f(z) + \sum_{J' \subsetneq J, J' \ne \emptyset}(-1)^{\#J'}f_{J',0}(z_{I_{J'}})
= (-1)^{\#J+1}f(z)\qquad z \in S'_J.
$$
Let $\tilde{z} \in S'_J$ and $z^*$ a point in $S$ with
$\mu_J(z^*,(0,\dots,0)) = \tilde{z}$. Then,
by putting $z = \mu_J(z^*,(\lambda,\dots,\lambda))$ into
(\ref{eq:tmp-tmp-xxx-yyy}) and letting $\lambda \to 0^+$, as
the right hand side of (\ref{eq:tmp-tmp-xxx-yyy}) tends to 0
by the same reason as that for $\#J = 1$,
we have $f_{J,0}(\tilde{z}) = f(\tilde{z})$.
This completes the proof.
\qed

Let us extend the notion of {\it{a consistent family of
coefficients of strongly asymptotic expansion}} defined by Majima in \cite{Ma84}
to our case.
Let $S:=S(U, \{G_i\}, \epsilon)$ be a proper multi-cone and
$F:=\{F_J\}_{J \in \Powl}$
a family of  coefficients of multi-asymptotic expansion with
$F_J =\{f_{J,\alpha}\}_{\alpha \in \Z^{I_J}_{\ge 0}}$.
For a proper multi-cone $S':=S(U', \{G'_i\}, \epsilon')$ properly
contained in $S$, we can consider the natural restriction $F\vert_{S'}$ of $F$
to $S'$. Moreover we can also define the restriction of $F$ to a submanifold.
Let $J \in \Powl$ and $\alpha \in \Z_{\ge 0}^{I_J}$.
Then the restriction $F\vert_{J,\alpha}$ of $F$ to $Z_J$ is defined as follows.

Let $J^* := \{j \in \Z_\ell; Z_J \nsubseteq Z_j\}$ and
$I^*_j = I_j \setminus I_J$ ($j \in J^*$).
We assume $J^* \ne \emptyset$ and
denote by $S^*$ the proper multi-cone $S_J$ in $Z_J$.
Then for $K \in \P(J^*) \setminus \emptyset$,
a family
$F^*_{K, \alpha} := \{f^*_{K, \alpha, \gamma} \}_{\gamma \in \Z_{\ge 0}^{I^*_K}}$
of holomorphic functions on $S^*_K \subset Z_J$
is defined by
$$
f^*_{K, \alpha, \gamma} = f_{J \cup K, (\alpha, \gamma)}\vert_{Z_J}
$$
where $(\alpha, \gamma) \in \Z^{I_J}_{\ge 0} \times \Z^{I^*_K}_{\ge 0}$
can be considered as an element in $\Z^{I_{J \cup K}}_{\ge 0}$ by
the identification $\Z^{I_J} \times \Z^{I^*_K} = \Z^{I_{J \cup K}}$.
Then the restriction $F\vert_{J,\alpha}$ of $F$ is defined by a family
$\{F^*_{K,\alpha}\}_{K \in \P(J^*) \setminus \emptyset}$.

The restriction $F\vert_{J,\alpha}$ gives a family of coefficients of multi-asymptotic expansion along
submanifolds $\{Z_j \cap Z_J\}_{j \in J^*}$ in $Z_J$ on the proper
multi-cone $S^* = S(U \cap Z_J, \{G_j\}_{j \in J^*}, \epsilon) \subset Z_J$.
\begin{df}
Let $S:=S(U, \{G_i\}, \epsilon)$ be a proper multi-cone and
$F := \{F_J\}_{J \in \Powl}$
a family of  coefficients of multi-asymptotic expansion with
$F_J =\{f_{J,\alpha}\}_{\alpha \in \Z^{I_J}_{\ge 0}}$.
We say that $F$ is a consistent family of coefficients of multi-asymptotic expansion
along ${\underset{1 \le j \le \ell}{\bigcup} Z_j}$ on $S$ if the following conditions are satisfied.
\begin{enumerate}
\item For any $J \in \Powl$ with
$J^* := \{j \in \Z_\ell; Z_J \nsubseteq Z_j\} \ne \emptyset$
and for any $\alpha \in \Z^{I_J}_{\ge 0}$,
a holomorphic function $f_{J,\alpha}$ on $S_J$
is multi-asymptotically developable to the family $F\vert_{J,\alpha}$ along
$\{Z_j \cap Z_J\}_{j \in J^*}$ on $S_J$.
\item For any $J$ and $J' \in \Powl$ with $S_J = S_{J'}$, we have $F_J = F_{J'}$.
\end{enumerate}
\end{df}

\begin{es}
Let us see some typical examples of consistent families.
\begin{enumerate}
\item
(Majima) Let us consider the Example \ref{Example:Majima}. The consistent families $F=(F_{\{1\}},F_{\{2\}},F_{\{1,2\}})$ are those satisfying

    - $F_{\{1\}}$ is asymptotic to $F_{\{1,2\}}$ when $z_2 \to 0$,

    - $F_{\{2\}}$ is asymptotic to $F_{\{1,2\}}$ when $z_1 \to 0$.
\item
(Takeuchi) Let us consider the Example \ref{Example:Takeuchi}. The consistent families $F=(F_{\{1\}},F_{\{2\}},F_{\{1,2\}})$ are those satisfying

     - $F_{\{1\}}$ is asymptotic to $F_{\{1,2\}}$ when $z_2 \to 0$,

     - $F_{\{2\}}=F_{\{1,2\}}$.

Therefore the case $N=(1,1)$ in Example \ref{Example:Takeuchi} is equivalent to
$$
|f(z)-f_{\{1\},0}(z_2,z')| \leq C_{S',N}|z_1|.
$$
\end{enumerate}
\end{es}

The following corollary immediately follows from
Theorem \ref{theorem:multi-asymptotic} and
Proposition \ref{prop:consistent-coeffs}.
\begin{cor}
Let $S:=S(U, \{G_i\}, \epsilon)$ be a proper multi-cone and
$F$ a family of  coefficients of multi-asymptotic expansion along ${\underset{1 \le j \le \ell}{\bigcup} Z_j}$
on $S$. If some holomorphic function on $S$ is multi-asymptotically developable to $F$,
then $F$ is a consistent family of multi-asymptotic expansion.
\end{cor}

Let $S':=S(U', \{G'_i\}, \epsilon')$ be a proper multi-cone
properly contained in $S$, and let $F:= \{F_J\}$ be a consistent family
on $S$ with $F_J = \{f_{J,\alpha}\}_{\alpha \in \Z^{I_J}_{\ge 0}}$.
By Theorem \ref{theorem:multi-asymptotic}, each $F_J$ can be regarded
as a $  \C^\infty$-Whitney jet on $\overline{S'_J}$.
(see \cite{M66} for Malgrange's definition of a $\C^\infty$-Whitney jet).
As $\overline{S'_J} \cap \overline{S'_{J'}} = \overline{S'_{J \cup J'}}$ holds
($J$, $J'$ $\in \hat{\P}(\Z_\ell)$),
Whitney jets defined by $F_J$ and $F_{J'}$ coincide
on the set $\overline{S'_{J \cup J'}}$
by Proposition \ref{prop:consistent-coeffs}.

Therefore, as $\underset{J \in \hat{\P}(\Z_\ell)}{\bigcup} \overline{S'_J} =
\overline{S'} \cap \left(\underset{1\le j \le \ell}{\bigcup} Z_j\right)$ holds by Lemma \ref{lemma:fundamental-geo},
it follows from Theorem 5.5 of \cite{M66} that we obtain the Whitney
jet defined on $\overline{S'} \cap \left(\underset{1\le j \le \ell}{\bigcup} Z_j\right)$
whose restriction to $\overline{S'_J}$ is equal to the one
defined by $F_J$.
Hence we have obtained the following proposition.

\begin{prop}{\label{prop:consisten-equiv}}
Let $S:=S(U, \{G_i\}, \epsilon)$ be a proper multi-cone and
$F:= \{F_J\}$ a family of  coefficients of multi-asymptotic expansion along
${\underset{1 \le j \le \ell}{\bigcup} Z_j}$ on $S$ with
$F_J = \{f_{J,\alpha}\}_{\alpha \in \Z^{I_J}_{\ge 0}}$.
Then the following conditions are equivalent.
\begin{enumerate}
\item $F$ is consistent.
\item For any proper multi-cone $S':=S(U',\{G_i\}, \epsilon')$ properly contained in $S$,
the restriction $F\vert_{S'}$ of $F$ to $S'$ can be
extended to a $\C^\infty$-function on $X_{\R}$, that is,
there exists a $\C^\infty$-function $\varphi(z)$ on $X_{\R}$ such that for any
$J \in \Powl$
and $\alpha \in \Z^{I_J}_{\ge 0}$ we have
$
f_{J, \alpha} =
\left. \displaystyle\frac{\partial^\alpha \varphi}
{\partial z^\alpha}\right|_{S'_J}
$
for $z \in S'_J$.
\end{enumerate}
\end{prop}

\end{section}

\begin{section}{Multi-specialization and asymptotic expansions}\label{8}

Let $X$ be a real analytic manifold. We consider a slight generalization of the sheaf of Whitney $\C^\infty$-functions of \cite{KS01}. As usual, given $F \in D^b(\CC_{X})$ we set $D'F=\rh(F,\CC_X)$. Remember that an open subset $U$ of $X$ is locally cohomologically trivial (l.c.t. for short) if $D'\CC_U \simeq \CC_{\overline{U}}$.

\begin{df} Let $F \in \mod_{\rc}(\CC_X)$ and let $U \in \op(X_{sa})$. We define the presheaf $\CW_{X|F}$ as follows:
$$U \mapsto \Gamma(X;H^0D'\CC_U \otimes F \wtens \C^\infty_X).$$
\end{df}
Let $U,V \in \op(X_{sa})$, and consider the exact sequence
$$\exs{\CC_{U\cap V}}{\CC_U \oplus \CC_V}{\CC_{U\cup V}},$$
applying the functor $\ho(\cdot,\CC_X)=H^0D'(\cdot)$ we obtain
$$\lexs{H^0D'\CC_{U\cup V}}{H^0D'\CC_U \oplus H^0D'\CC_V}{H^0D'\CC_{U\cap V}},$$
applying the exact functors $\cdot \otimes F$, $\cdot\wtens \C^\infty_X$ and taking global sections
we obtain
$$\lexs{\CW_{X|F}(U\cup V)}{\CW_{X|F}(U) \oplus \CW_{X|F}(V)}{\CW_{X|F}(U\cap V)}.$$
This implies that $\CW_{X|F}$ is a sheaf on $X_{sa}$. Moreover if $U \in \op(X_{sa})$ is l.c.t., the morphism $\Gamma(X;\CW_{X|F}) \to \Gamma(U;\CW_{X|F})$ is surjective and $R\Gamma(U;\CW_{X|F})$ is concentrated in degree zero. Let $\exs{F}{G}{H}$ be an exact sequence in $\mod_{\rc}(\CC_X)$, we obtain an exact sequence in $\mod(\CC_{X_{sa}})$
\begin{equation}\label{exsFGH}
\exs{\CW_{X|F}}{\CW_{X|G}}{\CW_{X|H}}.
\end{equation}

We can easily extend the sheaf $\CW_{X|F}$ to the case of $F \in D^b_{\rc}(\CC_X)$, taking a finite resolution of $F$ consisting of locally finite sums $\oplus \CC_V$ with $V$ l.c.t. in $\op^c(X_{sa})$. In fact, the sheaves $\CW_{X|\oplus \CC_V}$ form a complex quasi-isomorphic to $\CW_{X|F}$ consisting of acyclic objects with respect to $\Gamma(U;\cdot)$, where $U$ is l.c.t. in $\op^c(X_{sa})$.

As in the case of Whitney $\C^\infty$-functions one can prove that, if $G \in D^b_{\rc}(\CC_X)$ one has
$$\imin \rho \rh(G,\CW_{F|X}) \simeq D'G \otimes F \wtens \C^\infty_X.$$


\begin{es} Setting $F=\CC_X$ we obtain the sheaf of Whitney $\C^\infty$-functions. Let $Z$ be a closed subanalytic subset of $X$. Then $\CW_{X|\CC_{X\setminus Z}}$ is the sheaf of Whitney $\C^\infty$-functions vanishing on $Z$ with all their derivatives.
\end{es}

\begin{nt} Let $S$ be a locally closed subanalytic subset of $X$. We set for short $\CW_{X|S}$ instead of $\CW_{X|\CC_S}$.
\end{nt}

Let $\chi=\{M_1,\ldots,M_\ell\}$ be a family of closed analytic submanifolds of $X$ satisfying H1, H2 and H3, set $M= \bigcup_{i=1}^\ell M_i$ and consider $\widetilde{X}$. Consider the diagram \eqref{multi normal deformation}.

Set $F=\CC_{X\setminus M}$, $G=\CC_X$, $H=\CC_M$ in \eqref{exsFGH}. The exact sequence
$$\exs{\CW_{X|X\setminus M}}{\CW_X}{\CW_{X|M}}$$
induces an exact sequence
\begin{equation}\label{asym}
\exs{\nu^{sa}_\chi\CW_{X|X\setminus M}}{\nu^{sa}_\chi\CW_X}{\nu^{sa}_\chi \CW_{X|M}},
\end{equation}
in fact let $V$ be a l.c.t. conic subanalytic subset of the zero section of $\widetilde{X}$ and $U \in \op(X_{sa})$ such that
$C_\chi(X \setminus U) \cap V=\emptyset$, then we can find a l.c.t. $U'\subset U$ satisfying the same property. 
Applying the functor $\imin \rho$ to the exact sequence \eqref{asym} we obtain the exact sequence
$$
\exs{\imin \rho\nu^{sa}_\chi\CW_{X|X\setminus M}}{\imin \rho\nu^{sa}_\chi\CW_X}{\imin \rho \nu^{sa}_\chi \CW_{X|M}},
$$
where the surjective arrow is the map which associates to a function its asymptotic expansion.\\

Let $X$ be a complex manifold and let $X_\R$ denote the underlying real analytic manifold
of $X$. Let $\chi=\{Z_1,\ldots,Z_\ell\}$ be a family of complex submanifolds of $X$ satisfying H1, H2 and H3 and set $Z= \bigcup_{i=1}^\ell Z_i$. Let $F \in D^b_{\rc}(\CC_X)$. We denote by $\OW_{X|F}$ the sheaf defined as follows:
$$
\OW_{X|F} := \rh_{\rho_!\D_{\overline{X}}}(\rho_!\OO_{\overline{X}},\CW_{X_\R|F}).
 $$
Let $\exs{F}{G}{H}$ be an exact sequence in $\mod_{\rc}(\CC_X)$. Then the exact sequence \eqref{exsFGH} gives rise to the distinguished triangle
\begin{equation}\label{dtFGH}
\dt{\OW_{X|F}}{\OW_{X|G}}{\OW_{X|H}}.
\end{equation}

Setting $F=\CC_{X\setminus Z}$, $G=\CC_X$, $H=\CC_Z$ in \eqref{dtFGH} and applying the functor of specialization, we have the distinguished triangle
\begin{equation}\label{dtColin}
\dt{\imin \rho \nu^{sa}_\chi\OW_{X|X \setminus Z}}{\imin \rho \nu^{sa}_\chi\OWX}{\imin \rho \nu^{sa}_\chi \OW_{X|Z}}.
\end{equation}
The sheaf $\imin \rho \nu^{sa}_\chi\OWX$ is concentrated in degree zero.
This follows from the following result of \cite{DS96,Du79}: if $U \in \op(X_{sa})$ is convex, then $R\Gamma(X;\CC_{\overline{U}} \wtens \OO_X)$ is concentrated in degree zero.

Recall that we set $S:=  \underset{X,1\leq k\leq\ell}{\times}T_{Z_k}\iota(Z_k)$. We also set
\begin{equation}\label{eq:def-S-circ}
S^\circ := \{(q;\,\zeta_1,\dots,\zeta_\ell) \in S;\, \zeta_j \ne 0,\, j=1,\dots,\ell\}.
\end{equation}
We often say that $p$ is outside of the zero section of $S$ if $p \in S^\circ$.
\begin{prop}\label{vanishing} On $S^\circ$, the flat specialization along
$Z := Z_1 \cup \cdots \cup Z_\ell$
$$
\imin\rho\nu_\chi^{sa}\OW_{X|X \setminus Z}
$$
is concentrated in degree zero.
\end{prop}
\dim\ \
Set
$$I_j = \{i_{j,1},\dots,i_{j,m_j}\} \ \ (j = 1, 2, \dots , \ell).$$
Moreover, by a rotation of the coordinates, we may assume that
$$
\hat{I}_j = \{i_{j,1},\dots,i_{j,m'_j}\} \ \ (m'_j \leq m_j).
$$
We can identify $\CC^{{\hat{I}_j}}$ with $(T_{Z_j} \iota(Z_j))_0$. Let $p_j = (1, 0, \dots , 0) \in (T_{Z_j} \iota(Z_j))_0$.
Under the above identification, $p_j$ is the point in $\CC^{\hat{I}_j}$ with $z_{i_{j,1}} = 1$ and
$z_{i_{j,k}} = 0$  ($1 < k \leq m'_j$).
Let $G_j$ be a closed proper convex cone in $\CC^{\hat{I}_j}$ with direction $p_j$, which is defined by
$$
\{(z_{j,1},\dots,z_{j,m'_j}) \in \CC^{\hat{I}_j};\; |z_{j,1}|\leq\epsilon,\; |\arg(z_{j,1})|\leq\epsilon,\; |z_{j,k}|\leq\epsilon|z_{j,1}|,\;
k\geq 2\}
$$
for $0 < \epsilon < \frac{\pi}{2}$. Here we set $\operatorname{arg}(0) = 0$.
Let us define $G$ as the product
$$
G:=G_1 \times G_2 \times \cdots \times G_\ell \times B_\epsilon,
$$
where $B_\epsilon$ is a closed ball of radius $\epsilon$ and center at the origin
in $X' := \CC^{n-\sum_{j=1}^\ell m'_j}$. Set
\begin{eqnarray*}
T &:=& \{z \in \CC^n;\;z_{i_{1,1}}=0\} \cup \cdots \cup \{z \in \CC^n;\; z_{i_{\ell,1}}=0\} \\
&=& \{z \in \CC^n;\; z_{i_{1,1}} \cdots z_{i_{\ell,1}}=0\}.
\end{eqnarray*}
As
$$
G_j \cap \{(z_{i_{j,1}},\dots,z_{i_{j,m'_j}}) \in \CC^{\hat{I}_j};\;z_{i_{j,1}}=0\} = \{0\}
$$
holds, we have
$$
G \setminus T = (G_1 \setminus \{0\}) \times \cdots \times (G_\ell \setminus \{0\}) \times B_\epsilon.
$$
Therefore, by Proposition 5.4 of \cite{KS96}
$R\Gamma(X,\CC_{G\setminus T} \wtens \OO_X)$ is quasi-isomorphic to
$$
R\Gamma(\CC^{\hat{I}_1};\CC_{G_1\setminus\{0\}}\wtens\OO) \widehat{\boxtimes} \cdots \widehat{\boxtimes} R\Gamma(\CC^{\hat{I}_\ell};\CC_{G_\ell\setminus\{0\}}\wtens\OO) \widehat{\boxtimes} R\Gamma(X';\CC_{B_\epsilon} \wtens \OO),
$$
where $\widehat{\boxtimes}$ denotes the topological tensor product of \cite{Gr55}.
%

Each $R\Gamma(\CC^{\hat{I}_j};\CC_{G_j\setminus\{0\}} \wtens\OO)$,
$j=1,\dots,\ell$, is concentrated in degree zero
by the following lemma whose proof
will be given later.
\begin{lem}[Propostion 6.1.1 \cite{KS97}]{\label{vanishing-for-single-cone}}
For the closed set
$$
G :=
\{(z_1,\dots,z_n) \in \CC^{n};\; |z_1|\leq\epsilon,\; |\arg(z_1)|\leq\epsilon,\;
|z_{k}|\leq\epsilon|z_{1}|,\;
k\geq 2\},
$$
we have $H^k(\CC^{n};\CC_{G\setminus\{0\}} \wtens\OO_{\CC^n}) = 0$
for $k \ne 0$.
\end{lem}
Hence we have obtained
$$
H^k(X,\CC_{G\setminus T} \wtens \OO_X) = 0, \qquad (k \ne 0).
$$

Let us consider the the holomorphic map $f : X \to X$
$$
f(z_1,\dots,z_n) = (f_1(z), \dots, f_n(z)) :=
\left(\left( \prod_{j \in \hat{J}_1}z_{i_{j,1}}\right)z_1,\dots,\left(\prod_{j \in \hat{J}_n}z_{i_{j,1}}\right)z_n\right),
$$
where, for $k \in \{1, 2, \dots , n\}$,
$$
\hat{J}_k := \{j \in \{1, 2, \dots , \ell\};\; k \in I_j,\; k \notin \hat{I}_j\}.
$$
Clearly $f|_G$ is a proper map ($G$ is compact). We can also prove that $f$ induces
an isomorphism on $X \setminus T$.
Moreover we have
$$
p \in G \cap T \Longleftrightarrow f(p) \in f(G) \cap (Z_1 \cup \cdots \cup  Z_\ell).
$$
Therefore we get
$$
Rf_!\CC_{G\setminus T} = \CC_{f(G)\setminus (Z_1 \cup \cdots \cup Z_\ell)}.
$$
It follows from Theorem 5.7 of \cite{KS96} that we have
$$
R\Gamma(X;\rh_{\D_X}(\D_{X \stackrel{f}{\to} X},\CC_{G \setminus T} \wtens \OO_X) \simeq R\Gamma(X;\CC_{f(G)\setminus (Z_1 \cup \cdots \cup Z_\ell)} \wtens \OO_X).
$$
Set $\varphi := z_{i_{1,1}}z_{i_{2,1}}\cdots z_{i_{\ell,1}}$.
Let $\OO_{X,\varphi}$ denote the sheaf of meromorphic functions whose poles are contained in $\{\varphi = 0\}$
(i.e. $\OO_{X, \varphi} =
\OO_X[z_{i_{1,1}}^{-1},\, \dots,\, z_{i_{\ell,1}}^{-1}]$) and
we set $\D_{X,\varphi} := \OO_{X, \varphi} \underset{\OO_X}{\otimes} \D_{X}$.
As $\CC_{G\setminus T} \wtens \OO_X$ is a $\D_{X,\varphi}$-module and
$\D_{X \stackrel{f}{\to}X}$
is left
quasi-coherent over $\D_X$, we have
$$
\rh_{\D_X}(\D_{X \stackrel{f}{\to}X}, \CC_{G\setminus T}
\wtens \OO_X) \simeq  \rh_{\D_{X,\varphi}}(\D_{X,\varphi} \underset{\D_X}{\otimes}
 \D_{X \stackrel{f}{\to}X},\CC_{G\setminus T} \wtens \OO_X).
$$
We also have
$$
\D_{X,\varphi} \underset{\D_X}{\otimes}\D_{X \stackrel{f}{\to}X}\simeq \D_{X,\varphi}.
$$
As a matter of fact, $\D_{X \stackrel{f}{\to}X}$ is given by
$$
\D_{X \stackrel{f}{\to}X} = \frac{\D_{X\times X}}
{\D_{X \times X}(w_1 - f_1(z),\,\dots,\, w_n - f_n(z),\,
\theta_1,\, \dots,\, \theta_n)}
$$
where $(w_1, \dots, w_n)$ is a system of coordinates of the second $X$ and
the vector field $\theta_k$ ($k=1,2,\dots,n$) is defined by
$$
\theta_k := \frac{\partial}{\partial z_k} +
\sum_{i = 1}^n \frac{\partial f_i}{\partial z_k}\frac{\partial}{\partial w_i}.
$$
We denote by $J_f$ the Jacobian of $f$. Since $f$ gives an isomorphism
outside $T$,
we have $J_f \ne 0$ outside $T$, which implies
$J_f = z_{i_{1,1}}^{\beta_1}\dots z_{i_{\ell,1}}^{\beta_l} h$ for some
$\beta_j \ge 0$
$j=1,\dots,\ell$ and
for some $h$ with $h(0) \ne 0$.
Hence, as $J_f$ has an inverse in $\OO_{X, \varphi}$, we conclude that
$\D_{X,\varphi} \underset{\D_X}{\otimes}\D_{X \stackrel{f}{\to}X}$
is a free $\D_{X,\varphi}$ module of rank $1$.

Hence we have obtained
$$
R\Gamma(X;\CC_{G\setminus T} \wtens \OO_X) \simeq R\Gamma(X;\CC_{f(G)\setminus(Z_1\cup\cdots\cup Z_\ell)} \wtens \OO_X).
$$
This implies that $R\Gamma(X;\CC_{f(G)\setminus(Z_1\cup\cdots\cup Z_\ell)} \wtens \OO_X)$ is also concentrated in degree zero. Since
$f(G)$ contains $\overline{S}$ for some multi-cone $S$, and since, for a given multi-cone $S$,
there exists a $G$ such that $\overline{S}$ contains $f(G)$, we have the required result. \qed

To complete the proof of Proposition \ref{vanishing},
we give the proof of Lemma \ref{vanishing-for-single-cone}.\\

\noindent {\bf Proof of Lemma \ref{vanishing-for-single-cone}.}\ \
Let us consider the holomorphic map $f: \CC^n \to \CC^n$ defined by
$$
f(z_1, z_2, \dots, z_n) = (z_1, z_1z_2, \dots, z_1z_n).
$$
Then, by Proposition 5.4 and Theorem 5.7 of \cite{KS96} and
applying the same argument as that in the proof of Proposition
\ref{vanishing} to $f$,
we have
$$
R\Gamma(\CC^n;\, \CC_{G\setminus\{0\}}\wtens\OO_{\CC^n})
\simeq
R\Gamma(\CC;\, \CC_{(K \cap B^1_\epsilon)\setminus\{0\}}\wtens\OO_{\CC})
\widehat{\boxtimes} R\Gamma(\CC^{n-1};\, \CC_{B^{n-1}_\epsilon} \wtens \OO_{\CC^{n-1}})
$$
where $K$ is a closed cone defined by $\{z \in \CC; |\arg(z)|\leq\epsilon\}$
and $B^k_\epsilon$ is a closed ball in $\CC^k$
of radius $\epsilon$ and center at the origin.
Therefore it suffices to show
$$
H^1(\CC;\, \CC_{(K \cap B^1_\epsilon)\setminus\{0\}}\wtens\OO_{\CC}) = 0.
$$
We have the commutative diagram
$$
\begin{matrix}
\Gamma({\mathbb P}^1;\, \CC_{\{0\}}\wtens\OO_{{\mathbb P}^1})
&\to
&H^1({\mathbb P}^1;\, \CC_{\overline{K}\setminus{\{0\}}}\wtens\OO_{{\mathbb P}^1}) \\

 \|
&
& \downarrow
 \\

\Gamma(\CC;\, \CC_{\{0\}}\wtens\OO_{\CC})
& \to
& H^1(\CC;\, \CC_{(K\cap B^1_\epsilon)\setminus\{0\}}\wtens\OO_{\CC})
& \to
& 0
\end{matrix}
$$
Here the second row is exact because
$H^1(\CC;\, \CC_{(K\cap B^1_\epsilon)}\wtens\OO_{\CC}) = 0$  since $K \cap B^1_\epsilon$ is convex (see \cite{DS96,Du79}) and
$\overline{K}$ denotes the closure of $K$ in ${\mathbb P}^1$.
By applying Proposition 6.1.1 of \cite{KS97} to the case $V=\CC$ and
employing the coordinate transformation which
exchanges the origin and the infinity, we get
$H^1({\mathbb P}^1;\, \CC_{\overline{K}\setminus{\{0\}}}\wtens\OO_{{\mathbb P}^1}) = 0$.
Then
$
H^1(\CC;\, \CC_{(K\cap B^1_\epsilon)\setminus{\{0\}}}\wtens\OO_{\CC}) = 0
$
follows from the exactness of the second row.
\qed



\begin{prop}\label{prop: exact sequence} The distinguished triangle \eqref{dtColin} induces an exact sequence
	on $S^\circ$ (whose definition was given in $(\ref{eq:def-S-circ})$).
\begin{equation}\label{exspT}
\exs{\imin \rho H^0 \nu^{sa}_\chi\OW_{X|X \setminus Z}}{\imin \rho H^0 \nu^{sa}_\chi\OWX}{\imin \rho H^0 \nu^{sa}_\chi \OW_{X|Z}}.
\end{equation}
All the complexes $\imin \rho \nu^{sa}_\chi\OW_{X|X \setminus Z}$, $\imin \rho \nu^{sa}_\chi\OWX$ and $\imin \rho \nu^{sa}_\chi \OW_{X|Z}$ are concentrated in degree zero on $S^\circ$.
\end{prop}
\dim\ \ The exactness of the sequence \eqref{exspT} follows from Proposition \ref{vanishing}. The fact that outside the zero section $\imin \rho \nu^{sa}_\chi\OW_{X|X \setminus Z}$ and $\imin \rho \nu^{sa}_\chi\OW_X$ are concentrated in degree zero proves the vanishing of the cohomology of degree $\geq 1$ of the three terms.
\qed

We give a functorial construction of multi-asymptotically  developable
functions using Whitney holomorphic functions. Let $U$ be an open l.c.t. subanalytic subset in $X$.
Then $\Gamma(U;\CW_{X_\R}) \simeq \Gamma(X; \CC_{\overline{U}} \wtens C^\infty_{X_\R})$ is nothing
but the set of $\C^\infty$-Whitney jets on $\overline{U}$. Further
$H^0(U;\OW_X) \simeq H^0(X; \CC_{\overline{U}} \wtens \OO_X)$ consists of
$\C^\infty$-Whitney jets on $\overline{U}$ that satisfy the Cauchy-Riemann system.
Therefore, if a proper multi-cone $S:=S(U,\{G_j\}, \epsilon)$
is subanalytic, the set of holomorphic functions on $S$
that are multi-asymptotically developable along
$Z$ is equal to
$$
\varprojlim_{S'}
\Gamma(X;\CC_{\overline{S'}} \wtens \OO_X)
$$
where $S'$ runs through the family of open subanalytic proper multi-cones
$S(U',\, \{G'_j\},\, \epsilon')$ properly contained in $S$.
Moreover, by Proposition \ref{prop:consisten-equiv}, the set of consistent families
of coefficients of multi-asymptotic expansion is given by
$$
\varprojlim_{S'}
\Gamma(X; \CC_{\overline{S'} \cap Z} \wtens \OO_X).
$$

Let $W$ be an open convex subset of $Z_1 \cap \dots \cap Z_\ell$ and
$G_j$ an open proper convex cone of
$\CC^{\hat{I}_j}$
($j=1,2,\dots,\ell$).
Using these $W$ and $\{G_j\}$, we define
an $(\RP)^\ell$-conic open subset $V(W, \{G_j\})$ of
$S^\circ$ (its definition  was given in (\ref{eq:def-S-circ}))
by
$$
W \times G_1 \times \dots \times G_\ell
\subset \left(\bigcap Z_j\right) \times (\CC^{*})^{\hat{I}_1} \times \dots (\CC^{*})^{\hat{I}_\ell}
\simeq S^\circ.
$$
Note that a family of open sets of type $V(W, \{G_j\})$ forms a basis of the
topology on
$S^\circ$ for $(\RP)^\ell$-conic sheaves.
We define the presheaves
$\widetilde{\A}^{<0}_{\chi}$, $\widetilde{\A}_{\chi}$ and $\widetilde{\A}^{CF}_{\chi}$
on $S^\circ$. For an $(\RP)^\ell$-conic open set $V:=V(W, \{G_j\})$, we set
$$
\begin{aligned}
\Gamma(V;\widetilde{\A}^{<0}_{\chi}) &=  \lind {U,\, \epsilon > 0}
\{f \in \OO(S(U,\, \{G_j\}, \epsilon));\, \text{$f$ is multi-asymptotically}\\
&\qquad \text{developable to zero along $Z$ on $S(U,\, \{G_j\}, \epsilon)$}\},
\end{aligned}
$$
$$
\begin{aligned}
\Gamma(V;\widetilde{\A}_{\chi})
&= \varinjlim_{U,\, \epsilon > 0}
\{f \in \OO(S(U,\, \{G_j\}, \epsilon));\, \text{$f$ is multi-asymptotically}\\
&\qquad \text{developable along $Z$ on $S(U,\, \{G_j\}, \epsilon)$}\},
\end{aligned}
$$
$$
\begin{aligned}
\Gamma(V;\widetilde{\A}^{CF}_{\chi})
&= \varinjlim_{U,\, \epsilon > 0}
\{F;\, \text{$F$ is a consistent family of coefficients of}\\
&\qquad \text{multi-asymptotic expansion along $Z$
on $S(U,\, \{G_j\}, \epsilon)$}\},
\end{aligned}
$$
where
$U$ ranges through the family open neighborhoods of $W$
which has a form defined in (\ref{eq:def-open-set-product-type}).
Let us consider the multi-specialization of Whitney holomorphic functions. We have
$$
\begin{aligned}
\Gamma(V;\imin \rho \nu^{sa}_\chi \OW_{X|X\setminus Z}) &= \lpro {V'}\,\lind {U'} \Gamma(U';\OW_{X|X\setminus Z}),
\end{aligned}
$$
$$
\begin{aligned}
\Gamma(V;\imin \rho \nu^{sa}_\chi \OW_X) &= \lpro {V'}\,\lind {U'} \Gamma(U';\OW_X),
\end{aligned}
$$
$$
\begin{aligned}
\Gamma(V;\imin \rho \nu^{sa}_\chi \OW_{X|Z}) &= \lpro {V'}\,\lind {U'} \Gamma(U';\OW_{X|Z}),
\end{aligned}
$$
where $V'$ ranges through the family of subanalytic open cones such that $\overline{V'} \setminus Z \subset V$, $U'$ ranges through the family of $\op(X_{sa})$ such that
$C_\chi(X \setminus U') \cap V'=\emptyset$.

The identity morphism induces morphisms of presheaves
\begin{eqnarray*}
\widetilde{\A}^{<0}_\chi & \to & \imin \rho \nu^{sa}_\chi \OW_{X|X\setminus Z}, \\
\widetilde{\A}_\chi & \to & \imin \rho \nu^{sa}_\chi \OW_X, \\
\widetilde{\A}^{CF}_\chi & \to & \imin \rho \nu^{sa}_\chi \OW_{X|Z}.
\end{eqnarray*}
The above morphisms are isomorphisms in the stalks (i.e. in the limit of multi-cones containing a given direction
$p \in S^\circ$). Let $\A^{<0}_\chi$ (resp. $\A_\chi$, resp. $\A^{CF}_\chi$) be the sheaves associated to $\widetilde{\A}^{<0}_\chi$ (resp. $\widetilde{\A}_\chi$, resp. $\widetilde{\A}^{CF}_\chi$). We get \begin{eqnarray*}
\A^{<0}_\chi & \iso & \imin \rho \nu^{sa}_\chi \OW_{X|X\setminus Z}, \\
\A_\chi & \iso & \imin \rho \nu^{sa}_\chi \OW_X, \\
\A^{CF}_\chi & \iso & \imin \rho \nu^{sa}_\chi \OW_{X|Z}.
\end{eqnarray*}
By Proposition \ref{prop: exact sequence}, on $S^\circ$ we have the exact sequence of sheaves
\begin{equation}\label{eq:Borel-Ritt}
\exs{\A^{<0}_{\chi}}{\A_{\chi}}{\A^{CF}_{\chi}}.
\end{equation}
When $Z$ forms a normal crossing divisor,
these sheaves are nothing but the sheaves of strongly asymptotically developable
functions defined by Majima in \cite{Ma84}.

\begin{oss} The exact sequence \eqref{eq:Borel-Ritt} is nothing but a Borel-Ritt exact sequence for multi-asymptotically developable functions. This was already proven for formal specialization in the single divisor case in \cite{Co01} and for Majima's asymptotic in \cite{GS99,HS00}. Both results are based on the Borel-Ritt theorem in dimension one (for a proof, see \cite{Wa65}).
Thanks to Proposition \ref{vanishing} we obtained a purely cohomological proof of the Borel-Ritt exact sequence.
\end{oss}

\end{section}

\appendix
\begin{section}*{Appendix}
\end{section}
\begin{section}{Conic sheaves}

Let $k$ be a field. Let $X$ be a real analytic manifold endowed
with a subanalytic action $\mu$ of $\RP$. In other words we have a
subanalytic map
$$\mu: X \times \RP \to X,$$
which satisfies, for each $t_1,t_2 \in \RP$:
$$
  \begin{cases}
    \mu(x,t_1t_2)=\mu(\mu(x,t_1),t_2), \\
    \mu(x,1)=x.
  \end{cases}
$$
Note that $\mu$ is open, in fact let $U \in \op(X)$ and $(t_1,t_2)
\in \op(\RP)$. Then $\mu(U,(t_1,t_2))=\bigcup_{t \in
(t_1,t_2)}\mu(U,t)$, and $\mu(\cdot,t):X \to X$ is a homeomorphism
(with inverse $\mu(\cdot,\imin t)$). We have a diagram
$$\xymatrix{X \ar[r]^{\hspace{-5mm}j} & X \times \RP \ar@ <2pt>
[r]^{\hspace{0.3cm}\mu} \ar@ <-2pt> [r]_{\hspace{0.3cm}p}& X,}$$
where $j(x)=(x,1)$ and $p$ denotes the projection. We have $\mu
\circ j=p \circ j=\id$.

\begin{df} (i) Let $S$ be a subset of $X$. We set $\RP S=\mu(S,\RP).$
If $U \in \op(X)$, then $\RP U \in \op(X)$ since $\mu$ is open.

(ii) Let $S$ be a subset of $X$. We say that $S$ is conic if
$S=\RP S$. In other words, $S$ is invariant by the action of
$\mu$.

(iii) An  orbit of $\mu$ is the set $\RP x$ with $x \in X$.
\end{df}

We assume that the orbits of $\mu$ are contractible. For each $x \in X$ there are two possibilities: either $\RP x=x$ or $\RP x \simeq \R$.

\begin{df} We say that a subset $S$ of $X$ is $\RP$-connected if
$S \cap \RP x$ is connected for each $x \in S$.
\end{df}

\begin{lem} \label{lem: RP-connectedness} (i) Let $S_1,S_2 \subset X$ and suppose that $S_2$ is conic. Then $\RP(S_1 \cap S_2)=\RP S_1 \cap S_2$. (ii) If $S_1$ and $S_2$ are $\RP$-connected, then $S_1 \cap S_2$ is $\RP$-connected.
\end{lem}
\dim\ \ (i) The inclusion $\subseteq$ is true since $\RP(S_1 \cap S_2) \subset \RP S_i$, $i=1,2$. Let us prove $\supseteq$. Let $x \in \RP S_1 \cap S_2$. Then there exists $a \in \RP$ such that $\mu(x,a) \in S_1$. Since $S_2$ is conic $\mu(x,a) \in S_2$ and the result follows.

(ii) Let $x_1,x_2 \in S_1 \cap S_2 \cap \RP x$ for some $x \in X$. Then $x_i=\mu(x,a_i)$ and some $a_i \in \RP$, $i=1,2$. Suppose $a_1 \leq a_2$. Since $S_1,S_2$ are $\RP$-connected, $i=1,2$ and the orbits of $\mu$ are either points or isomorphic to $\R$, we have $\mu(x,[a_1,a_2]) \subseteq S_i$, $i=1,2$ and the result follows.
\qed

Let $X,Y$ topological spaces endowed with an action ($\mu_{X}$ and $\mu_{Y}$ respectively) of $\RP$.

\begin{df} A continuous function $f:X \to Y$ is said to be conic if for each $x \in X$, $a \in \RP$ we have
$f(\mu_{X}(x,a))=\mu_{Y}(f(x),a)$.
\end{df}

\begin{lem} \label{lem: conic maps} Let $f:X \to Y$ be a conic map. (i) Suppose that $S \subset Y$ is $\RP$-connected (resp. conic). Then $\imin f(S)$ is $\RP$-connected (resp. conic). (ii) Suppose that $Z \subset X$ is conic. Then $f(Z)$ is conic.
\end{lem}
\dim\ \ (i) Let $x_1,x_2 \in \imin f(S)$ and suppose that there exists $x \in X$ such that $x_1,x_2 \in \mu_X(x,\RP)$, i.e. $x_i=\mu_X(x,a_i)$, $a_i \in \RP$, $i=1,2$. Since $f$ is conic we have $f(x_i)=\mu_Y(f(x),a_i) \in \mu_Y(f(x),\RP) \cap S=f(\mu_X(x,\RP) \cap \imin f(S))$. Suppose $a_1 \leq a_2$. Since $S$ is $\RP$-connected $\mu_Y(f(x),[a_1,a_2]) \subseteq S$, hence $\mu_X(x,[a_1,a_2]) \subseteq \imin f (f(\mu_X(x,[a_1,a_2])))=\imin f(\mu_Y(f(x),[a_1,a_2])) \subseteq \imin f(S)$.

If $S$ is conic then $f(\mu_X(\imin f(S),\RP))=\mu_Y(f(\imin f(S)),\RP)=\mu_Y(S,\RP)=S$, hence $\mu_X(\imin f(S),\RP) =\imin f(S)$.

(ii) We have $\mu_Y(f(Z),\RP)=f(\mu_X(Z,\RP))=f(Z)$.
\qed

Let $X$ be a real analytic manifold endowed with a subanalytic action of $\RP$. Denote by $X_{sa}$ the associated subanalytic site.

\begin{df} A sheaf of $k$-modules $F$ on $X_{sa}$ is conic if the restriction morphism
$\Gamma(\RP U;F) \to \Gamma(U;F)$ is an isomorphism for each
$\RP$-connected $U \in \op^c(X_{sa})$ with $\RP U \in
\op(X_{sa})$.
\begin{itemize}
\item[(i)]We denote by $\mod_{\RP}(k_{X_{sa}})$ the subcategory of
$\mod(k_{X_{sa}})$ consisting of conic sheaves.
\item[(ii)] We denote by $D^b_{\RP}(k_{X_{sa}})$,
the subcategory of $D^b(k_{X_{sa}})$ consisting of objects $F$
such that $H^j(F)$ belongs to $\mod_{\RP}(k_{X_{sa}})$ for all $j
\in \Z$.
\end{itemize}
\end{df}

Assume the hypothesis below:
\begin{equation}\label{hypsa}
  \begin{cases}
 \text{(i) every $U \in \op^c(X_{sa})$ has a finite covering consisting }\\
 \text{\ \ of $\RP$-connected subanalytic open subsets,}\\
 \text{(ii) for any $U \in \op^c(X_{sa})$ we have $\RP U \in \op(X_{sa})$,}\\
 \text{(iii) for any $x \in X$\ \ the set $\RP x$ is contractible,}\\
 \text{(iv) there exists a covering $\{V_n\}_{n \in \N}$ of $X_{sa}$ such that}\\
 \text{\ \ $V_n$ is $\RP$-connected and $V_n \subset\subset V_{n+1}$ for each $n$}.
  \end{cases}
\end{equation}

The following result was proven in \cite{Pr11}.

\begin{prop}\label{RPU} Assume \eqref{hypsa}. Let $U \in \op(X_{sa})$ be $\RP$-connected and such that $\RP U \in \op(X_{sa})$. Let $F \in D^b_{\RP}(k_{X_{sa}})$. Then
$$
R\Gamma(\RP U;F) \iso R\Gamma(U;F).
$$
\end{prop}

\end{section}

\begin{section}{Multi-conic sheaves}

Let $X$ be a topological space
with $\ell$ actions $\{\mu_i\}_{i=1}^\ell$ of $\RP$ such that $\mu_i(\mu_j(x,t_j),t_i)=\mu_j(\mu_i(x,t_i),t_j)$. We have a map
\begin{eqnarray*}
\mu: X \times (\RP)^\ell & \to & X \\
(x,(t_1,\dots,t_\ell)) & \mapsto & \mu_1(\cdots\mu_\ell(x,t_\ell),\dots,t_1).
\end{eqnarray*}

\begin{df} (i) Let $S$ be a subset of $X$. We set $\RP_i S=\mu_i(S,\RP).$
If $U \in \op(X)$, then $\RP_i U \in \op(X)$ since $\mu_i$ is open for each $i=1,\ldots,\ell$.

(ii)  Let $S$ be a subset of $X$. Let $J = \{i_1,\ldots,i_k\} \subseteq \{1,\ldots,\ell\}$. We set $$\RP_JS=\RP_{i_1}\cdots\RP_{i_k}S=\mu_{i_1}(\cdots\mu_{i_k}(S,\RP),\dots,\RP), \ \ i_1,\ldots i_k \in J.$$
We set $(\RP)^\ell S=\RP_{\{1,\ldots,\ell\}}S=\mu(S,(\RP)^\ell).$
If $U \in \op(X)$, then $\RP_J U \in \op(X)$ since $\mu_i$ is open for each $i \in \{1,\ldots,\ell\}$.

(iii) Let $S$ be a subset of $X$. We say that $S$ is $(\RP)^\ell$-conic if
$S=(\RP)^\ell S$. In other words, $S$ is invariant by the action of
$\mu_i$, $i=1,\dots,\ell$.

\end{df}

\begin{df} (i)  We say that a subset $S$ of $X$ is $\RP_i$-connected if
$S \cap \RP_i x$ is connected for each $x \in S$.

(ii) We say that a subset $S$ of $X$ is $(\RP)^\ell$-connected if
there exists a permutation $\sigma:\{1,\ldots,\ell\}\to\{1,\ldots,\ell\}$ such that
\begin{equation}\label{RPl connected}
 \begin{cases}
 \text{$S$ is $\RP_{\sigma(1)}$-connected,}\\
 \text{$\RP_{\sigma(1)}S$ is $\RP_{\sigma(2)}$-connected,}\\
 \text{\ \ $\vdots$}\\
 \text{$\RP_{\sigma(1)}\cdots\RP_{\sigma(\ell-1)}S$ is $\RP_{\sigma(\ell)}$-connected}.
  \end{cases}
\end{equation}
\end{df}

The following results follow from the case $\ell=1$.

\begin{lem} \label{lem: l-connectedness} (i) Let $S_1,S_2 \subset X$ and suppose that $S_2$ is $(\RP)^\ell$-conic. Then $(\RP)^\ell(S_1 \cap S_2)=(\RP)^\ell S_1 \cap S_2$. (ii) If moreover $S_1$ is $(\RP)^\ell$-connected then $S_1 \cap S_2$ is $(\RP)^\ell$-connected.
\end{lem}


\begin{oss} In (ii) of Lemma \ref{lem: l-connectedness} we have to assume that $S_2$ is $(\RP)^\ell$-conic. Indeed it is not true that the intersection of two $(\RP)^\ell$-connected is $(\RP)^\ell$-connected in general.
\end{oss}

Let $X,Y$ topological spaces endowed with $\ell$ actions $\{\mu_{Xi}\}_{i=1}^\ell$, $\{\mu_{Yi}\}_{i=1}^\ell$ of $\RP$.

\begin{df} A continuous function $f:X \to Y$ is said to be $(\RP)^\ell$-conic if for each $x \in X$, $a \in \RP$ we have
$f(\mu_{Xi}(x,a))=\mu_{Yi}(f(x),a)$, $i=1,\ldots,\ell$.
\end{df}

\begin{lem} \label{lem: l-conic maps} Let $f:X \to Y$ be a $(\RP)^\ell$-conic map. (i) Suppose that $S \subset Y$ is $(\RP)^\ell$-connected (resp. $(\RP)^\ell$-conic). Then $\imin f(S)$ is $(\RP)^\ell$-connected (resp. $(\RP)^\ell$-conic). (ii) Suppose that $Z \subset X$ is $(\RP)^\ell$-conic. Then $f(Z)$ is $(\RP)^\ell$-conic.
\end{lem}

Let $X$ be a real analytic manifold and denote by $X_{sa}$ the associated subanalytic site. Assume that $X$ is endowed with $\ell$ subanalytic $\RP$-actions $\mu_1,\dots,\mu_\ell$ commuting with each other.

\begin{df}\label{def:conic sheaf} A sheaf of $k$-modules $F$ on $X_{sa}$ is $(\RP)^\ell$-conic if it is conic with respect to each $\mu_i$.
\begin{itemize}
\item[(i)]We denote by $\mod_{(\RP)^\ell}(k_{X_{sa}})$ the subcategory of
$\mod(k_{X_{sa}})$ consisting of $(\RP)^\ell$-conic sheaves.
\item[(ii)] We denote by $D^b_{(\RP)^\ell}(k_{X_{sa}})$,
the subcategory of $D^b(k_{X_{sa}})$ consisting of objects $F$
such that $H^j(F)$ belongs to $\mod_{(\RP)^\ell}(k_{X_{sa}})$ for all $j
\in \Z$.
\end{itemize}
\end{df}

Let us assume the following hypothesis
\begin{equation}\label{hypsam}
  \begin{cases}
 \text{(i) the pair $(X,\mu_i)$ satisfies \eqref{hypsa} for each $i=1,\dots,\ell$, }\\
 \text{(ii) every $U \in \op^c(X_{sa})$ has a finite covering consisting }\\
 \text{\ \ of $(\RP)^\ell$-connected subanalytic open subsets,}\\
 \text{(iii) we have $\RP_J U \in \op(X_{sa})$ for any $U \in \op^c(X_{sa})$}\\
 \text{\ \ and any $J \subset \{1,\ldots,\ell\}$.}\\
  \end{cases}
\end{equation}
In this situation the orbits of $\mu_i$, $i=1,\dots,\ell$ are either $\RP x \simeq \R$ or $\RP x = x$.\\




\begin{prop}\label{prop: RPUl} Assume \eqref{hypsam}. Let $U \in \op(X_{sa})$ be $(\RP)^\ell$-connected. Let $F \in D^b_{(\RP)^\ell}(k_{X_{sa}})$. Then
$$R\Gamma((\RP)^\ell U;F) \iso R\Gamma(U;F).$$
\end{prop}
\dim\ \ Suppose that $U$ satisfies \eqref{RPl connected}. By Proposition \ref{RPU} we have
$$R\Gamma(U;F) \simeq R\Gamma(\RP_{\sigma(1)}U;F) \simeq \cdots \simeq R\Gamma(\RP_{\sigma(1)}\cdots\RP_{\sigma(\ell)}U;F)$$
and $\RP_{\sigma(1)}\cdots\RP_{\sigma(\ell)}U=(\RP)^\ell U$. \qed


If $X$ satisfies \eqref{hypsam} (i)-(iii), then it follows from Proposition \ref{prop: RPUl} that for $(\RP)^\ell$-conic subanalytic sheaves it is enough to study the cohomology of the sections on $(\RP)^\ell$-conic open subsets.\\

\section{Multi-actions in $\R^n$}

Let $\mu$ be an action of $\RP$ on $\R^n$ defined by
$$
\mu:((x_1,\ldots,x_n),c) \mapsto (c_{1}x_1,\ldots,c_{n}x_n),
$$
where $c_{i}=1$ if $i \notin I$, $c_{i}=c \in \RP$ if $i \in I$, $I \subseteq \{1,\dots,n\}$.

\begin{lem}\label{lem: RP decomposition}  Each globally subanalytic $U \in \op(\R^n_{sa})$ has a finite covering consisting of globally subanalytic $\RP$-connected open subsets.
\end{lem}
\dim\ \ We may assume that $U$ is connected. In the rest of the proof we write for short subanalytic instead of globally subanalytic. Suppose that $I=\{1,\dots,m\}\subseteq\{1,\dots,n\}$ and consider an action
\begin{eqnarray*}
\mu: \R^n \times \RP & \to & \R^n \\
((x_1,\ldots,x_m,x_{m+1},\ldots,x_n),\lambda) & \mapsto & (\lambda x_1,\ldots,\lambda x_m,x_{m+1},\ldots,x_n).
\end{eqnarray*}
Let us consider the morphism of manifolds
\begin{eqnarray*}
\varphi : \mathbb{S}^{m-1} \times \R \times \R^{n-m} & \to & \R^n\\
(\vartheta,r,z) & \mapsto & (r i(\vartheta),z),
\end{eqnarray*}
where $i:\mathbb{S}^{m-1} \hookrightarrow \R^m$ denotes the
embedding. The map $\varphi$ is proper and subanalytic (even semi-algebraic). The subset
$\imin \varphi(U)$ is subanalytic in $\mathbb{S}^{m-1} \times \R \times \R^{n-m}$. Set $Z=\{0\} \times \R^{n-m}$.
\begin{itemize}
\item[(a)] Let us consider a cylindrical cell decomposition of $\imin\varphi(U \setminus Z)$ (for the definition see \cite{Co00}) with respect to the fibers of the projection $\pi:\mathbb{S}^{m-1} \times \R \times \R^{m-n} \to \mathbb{S}^{m-1} \times \R^{m-n}$. Each cell $D$ is defined by $(f,g):=\{(\vartheta,r,z);\;f(\vartheta,z)<r<g(\vartheta,z)\}$, where $f,g:D' \to \R \cup \{\pm \infty\}$ are subanalytic functions (we allow $f(\vartheta,z)=-\infty$ and $g(\vartheta,z)=+\infty$) and $D'$ is a cell of $\pi(\imin\varphi(U \setminus Z))$. Consider $h:D'\to \R$, $f<h<g$, for example $h={f+g \over 2}$ (or $h=f+1$ if $g=+\infty$, $h=g-1$ if $f=-\infty$, $h=0$ if $f=-\infty$ and $g=+\infty$). Extend the graph of $h$ to $\widetilde{h}$
on an open subanalytic neighborhood $U'$ of $D'$
such that the graph of $\widetilde{h}$ is contained in $\imin\varphi(U \setminus Z)$. For $p \in U'$ let $(m(p),M(p))$ be the connected component of $\imin {\pi}(p) \cap \imin\varphi(U \setminus Z)$ containing $\widetilde{h}(p)$. The functions $m,M:U'\to\R$ are subanalytic and the open set $U_D=\{(p,y);\; p \in U',\; m(p)<y<M(p)\}$ has the required properties.
In this way we obtain a finite covering $\{W_j\}_{j \in J}$ of $\imin \varphi(U \setminus Z)$ consisting of $\RP$-connected subanalytic open subsets.

\item[(b)] Let $p \in \pi(\imin\varphi(U \cap Z))$. Then $\imin \pi(p) \cap \imin\varphi(U)$ is a disjoint union of intervals. Let us consider the interval $(m(p),M(p))$, $m(p)<M(p) \in \R$ containing $0$. Set $W_Z=\{(p,r) \in U ;\; m(p)<r<M(p)\}$. The set $W_Z$ is open subanalytic, contains $\imin\varphi(U \cap Z)$ and its intersections with the fibers of $\pi$ are connected. Then $\varphi(W_Z)$ is an open $\RP$-connected subanalytic neighborhood of $U \cap Z$ and it is contained in $U$.
\end{itemize}
By (a) there exists a finite covering $\{\varphi(W_j)\}_{j \in J}$ of $U \setminus Z$ consisting of $\RP$-connected subanalytic open subsets, and $\varphi(W_Z) \cup \bigcup_{j \in J}\varphi(W_j)=U$. \qed

\begin{teo} \label{teo: hypsa} Let $\mu$ be an action of $\RP$ on $\R^n$ defined by
$$
\mu:((x_1,\ldots,x_n),c) \mapsto (c_{1}x_1,\ldots,c_{n}x_n),
$$
where $c_{i}=1$ or $c_{i}=c \in \RP$. Then $\R^n$ satisfies \eqref{hypsa} (i)-(iv).
\end{teo}
\dim\ \ \eqref{hypsa} (i) follows from Lemma \ref{lem: RPl decomposition}, \eqref{hypsa} (ii) follows from the fact that $\mu$ is a globally subanalytic map, hence if $U$ is globally subanalytic  $\mu(U,\RP)$ is still globally subanalytic,
(iii) and (iv) are trivial. \qed

Let $\{\mu_i\}_{i=1}^\ell$, $\ell \leq n$ be actions of $\RP$ on $\R^n$ defined by
$$
\mu_j:((x_1,\ldots,x_n),c_j) \mapsto (c_{1j}x_1,\ldots,c_{nj}x_n),
$$
where $c_{ij}=1$ if $i \notin I_j$, $c_{ij}=c_j \in \RP$ if $i \in I_j$. We assume the following hypothesis
\begin{equation}\label{H}
\text{Either $I_i \subsetneqq I_j$, $I_j \subsetneqq I_i$ or $I_i
\cap I_j = \emptyset$ holds
($i \ne j \in \{1,2,\dots,\ell\}$).}
\end{equation}
\begin{equation} \label{H3}
\text{ $I_i \varsupsetneqq \bigcup_{I_j \varsubsetneqq I_i}I_j$\; for $i=1,2,\ldots,\ell$.}
\end{equation}

\begin{df} Let $i,j = 1,\ldots,\ell$. When $I_i \cap I_j = \emptyset$ the actions $\mu_i$ and $\mu_j$ are said to be orthogonal and we write $\mu_i \bot \mu_j$.
\end{df}

\begin{lem}\label{lem: RPl bot decomposition}  Suppose that $\mu_i \bot \mu_j$ for each $i,j \in \{1,\ldots,\ell\}$. Each globally subanalytic $U \in \op(\R^n_{sa})$ has a finite covering consisting of globally subanalytic $(\RP)^\ell$-connected open subsets.
\end{lem}
\dim\ \ We may assume that $U$ is connected. In the rest of the proof we write for short subanalytic instead of globally subanalytic.\\

Suppose that $I_i \cap I_j = \emptyset$ for each $i,j \in \{1,\ldots,\ell\}$, i.e. all the actions are orthogonal. Let $U$ be $\RP_1$-connected and such that $\RP_1 \cdots \RP_{j-1} U$ is $\RP_j$-connected for each $j=2,\ldots,\ell-1$. We are going to find a finite cover $\{U_\alpha\}$ with $U_\alpha$ $\RP_1$-connected and such that $\RP_1 \cdots \RP_{j-1} U_\alpha$ is $\RP_j$-connected for each $j=2,\ldots,\ell$.

For any subset $J=\{j_1,\ldots,j_k\}$ of $\{1,\ldots,\ell-1\}$ let us consider
$$
S_J := \left\{x \in \R^n;\; \sum_{i \in I_k}x_i^2=1 \; (k \in J),\;x_j=0\;(j \in I_k,\; k \notin J)\right\}.
$$
By construction $S_J$ is invariant under $\mu_\ell$. Let us consider $\RP_1 \cdots \RP_{\ell-1}U \cap S_J$. By Lemma \ref{lem: RP decomposition} it admits a finite covering $\{V_{\beta_J}\}$ consisting of $\RP_\ell$-connected subanalytic open subsets of $S_J$. For each $\beta_J$ consider $U_{\beta_J}:=\RP_J(\imin {\pi_J}(\pi_J(V_{\beta_J})))$, where
$n_J=\sharp\left(\bigcup_{i \in J}I_i \cup I_\ell\right)$ (hence $\R^{n_J}$ is represented by the coordinates $x_k$, $k \in \bigcup_{i \in J}I_i \cup I_\ell$) and $\pi_J:\R^n \to \R^{n_J}$ is the projection. We want to show that $U \cap U_{\beta_J}$ has the desired properties. First remark that by Lemma \ref{lem: l-connectedness} (i) $\RP_1\cdots\RP_j(U \cap U_{\beta_J}) = \RP_1\cdots\RP_jU \cap U_{\beta_J}$ for each $j=1,\ldots,\ell-1$ since $U_{\beta_J}$ is conic with respect to $\mu_1,\ldots,\mu_{\ell-1}$ by construction.
In particular by Lemma \ref{lem: l-connectedness} it is $\RP_j$-connected for each $j=1,\ldots,\ell-1$.
Moreover $\RP_1\cdots\RP_{\ell-1}U \cap U_{\beta_J} = U_{\beta_J}$ which is $\RP_\ell$-connected by construction. Then $U \cap U_{\beta_J}$ has the desired properties and the result follows. \qed

\begin{oss}  Remark that the proof of Lemma \ref{lem: RPl bot decomposition} does not depend on the choice of the permutation $\sigma$ in \eqref{RPl connected}.
\end{oss}

\begin{lem}\label{lem: RPl decomposition}  Each globally subanalytic $U \in \op(\R^n_{sa})$ has a finite covering consisting of globally subanalytic $(\RP)^\ell$-connected open subsets.
\end{lem}
\dim\ \ We may assume that $U$ is connected. In the rest of the proof we write for short subanalytic instead of globally subanalytic.\\

Let $J$ be the subset of $\{1,\dots,\ell\}$ defined as follows:
$$
J=\left\{i \in \{1,\ldots,\ell\},\; I_i \cap \bigcup_{i\neq j}I_j=\emptyset\right\},
$$
i.e. $\sharp J$ denotes the number of $\mu_i$ such that $\mu_i \bot \mu_j$ for each $j\neq i$. Note that $0 \leq \sharp J \leq \ell$. We argue by double induction: by increasing induction for $\ell$ (from
1) and by decreasing induction for $\sharp J$ (from $\ell$ to 0).
If $\sharp J=\ell$ ($\ell=1,2,\dots$) the result follows by Lemma \ref{lem: RPl bot decomposition}. Suppose that the result holds for $\ell'$ smaller than $\ell$ and
any $0 \leq \sharp J \leq \ell$ and for $\ell=\ell'$ and $\sharp J = k,\dots,\ell$. We will
show it for $\ell=\ell'$ and $\sharp J = k-1$.
Thanks to \eqref{H} and \eqref{H3}, there exists $i \in \{1,\ldots,\ell\}$ such that $I_i \supset I_j$ for some $j \neq i$ and $I_i \setminus \bigcup_{j \neq i}I_j \neq \emptyset$. Up to take a permutation of $\{1,\ldots,\ell\}$ we may assume $i=\ell$.

Let $U$ be $\RP_1$-connected and such that $\RP_1 \cdots \RP_{j-1} U$ is $\RP_j$-connected for each $j=1,\ldots,\ell-1$. We are going to find a finite cover $\{U_\alpha\}$ with $U_\alpha$ $\RP_1$-connected and such that $\RP_1 \cdots \RP_{j-1} U_\alpha$ is $\RP_j$-connected for each $j=1,\ldots,\ell$.
\begin{itemize}

\item[(a)] Given $i_0 \in I_\ell \setminus \bigcup_{j\neq \ell}I_j$ let us consider $U':=U \cap \{x_{i_0} \neq 0\}$ and consider the subanalytic (even semi-algebraic) homeomorphism $\psi$ of $\{x_{i_0} \neq 0\}$ given by
    $$
    \begin{cases}
    \text{$\psi(x)_k=\displaystyle{x_k\over x_{i_0}}$ \ \ if $k \in I_j \subset I_\ell$, $j \neq \ell$,} \\
    \text{$\psi(x)_k=x_k$ \ \ otherwise.}
    \end{cases}
    $$
    Set $\mu'_k=\mu_k \circ (\psi \times \id_{(\RP)^\ell})$. By construction we have $\mu'_\ell \bot \mu'_j$ for each $j \neq \ell$, hence by the induction hypothesis we may find a finite cover $\{V_\alpha\}$ of $\psi(U')$ (which is subanalytic, being $\psi$ semi-algebraic) with the required properties. Then $\{\imin \psi(V_\alpha)\}$ is a finite cover of $U'$ consisting of $(\RP)^\ell$-connected open subanalytic subsets.

\item[(b)]  By \eqref{H} $\bigcup_{I_j \varsubsetneqq I_\ell}I_j$ is a disjoint union $\hat{I}_{i_1} \sqcup \cdots \sqcup \hat{I}_{i_k}$, $\{i_1,\dots,i_k\} \subseteq \{1,\dots,\ell-1\}$.
    Set $x=(x_i)$, $i \in I_\ell \setminus \bigcup_{j \neq \ell}I_j$, $y=(x_j)$, $j \in \hat{I}_{i_1} \sqcup \cdots \sqcup \hat{I}_{i_k}$ $z=(x_k)$, $k \notin I_\ell$ 
    Assume that $I_\ell \setminus \bigcup_{j \neq \ell}I_j=\{1,\ldots,m\}$.
    Set $S_0:=\{x=0\}$, and consider $V_0:=U \cap S_0$. It is an open subanalytic subset of $S_0$, it is
    $\RP_1$-connected and such that $\RP_1 \cdots \RP_{j-1} V_0$ is $\RP_j$-connected for each $j=1,\ldots,\ell$. Indeed on $S_0$ $\mu_\ell$ is generated by $\mu_1,\ldots,\mu_{\ell-1}$ and hence $\RP_1\cdots\RP_{\ell-1} V_0$ is conic with respect to $\mu_\ell$.
Let us consider the morphism of manifolds
\begin{eqnarray*}
\varphi : \mathbb{S}^{m-1} \times \R \times \R^{n-m} & \to & \R^n\\
(\vartheta,r,y,z) & \mapsto & (r i(\vartheta),y,z),
\end{eqnarray*}
where $i:\mathbb{S}^{m-1} \hookrightarrow \R^m$ denotes the
embedding. Endow $\mathbb{S}^{m-1} \times \R \times \R^{n-m}$ with the actions $\mu'_j=\mu_j$, $j=1,\dots,\ell-1$, $\mu'_\ell((\vartheta,r,y,z),\lambda)=(\vartheta,\lambda r,\lambda y,z)$. Then $\varphi$ is a conic map.

Denote by $\pi:\mathbb{S}^{m-1} \times \R \times \R^{n-m} \to \mathbb{S}^{m-1} \times \R^{n-m}$ the projection. Up to shrink $U$, we may suppose that $\imin \varphi(U)=\imin \varphi(U) \cap \imin \pi(\pi(\imin \varphi(V_0)))$. Let $p \in \imin \varphi(\RP_1 \cdots \RP_{\ell-1}V_0)$. Then $\imin \pi(p) \cap \imin \varphi(\RP_1 \cdots \RP_{\ell-1}U)$ is a disjoint union of intervals. Let us consider the interval $(m(p),M(p))$, $m(p)<M(p) \in \R$ containing $0$. Set
$$
U_0=\varphi\left(\{(p,r);\;p \in \imin \varphi(\RP_1 \cdots \RP_{\ell-1}V_0),\; m(p)<r<M(p)\}\right).
$$
The set $U_0$ is open subanalytic, contains $U \cap S_0$ and if $(x,y,z) \in U_0$, then $(0,y,z) \in U_0$. Moreover by Lemma \ref{lem: l-conic maps} (ii) $U_0$ is conic with respect to $\mu_1,\ldots,\mu_{\ell-1}$. We want to show that $U \cap U_0$ has the desired properties. First remark that by Lemma \ref{lem: l-connectedness} (i) $\RP_1\cdots\RP_j(U \cap U_0) = \RP_1\cdots\RP_jU \cap U_0$ for each $j=1,\ldots,\ell-1$ since $U_0$ is conic with respect to $\mu_1,\ldots,\mu_{\ell-1}$. In particular by Lemma \ref{lem: l-connectedness} (ii) it is $\RP_j$-connected for each $j=1,\ldots,\ell-1$. We prove that $\RP_1\cdots\RP_{\ell-1}U \cap U_{0} = U_{0}$ is $\RP_\ell$-connected. Suppose that $(x,y,z),(\lambda x,\lambda y,z) \in U_0$ for some $\lambda \in \RP$. By construction $(0,y,z),(0,\lambda y,z) \in U_0$. We may assume without loss of generality $\lambda >1$. We argue by contradiction. Suppose that there exists $1<\eta<\lambda$ such that $(\eta x,\eta y,z) \notin U_0$. Since $U_0$ is conic with respect to $\mu_{i_1},\dots,\mu_{i_k}$, for each $t \in \RP$ we have $(\eta x, t\eta y,z) \notin U_0$. Set $t=\lambda/\eta$. Then $(\eta x,\lambda y,z) \notin U_0$ which leads to a contradiction since by construction $(\eta x,\lambda y,z) \in U_0$ for each $\eta \in (0,\lambda)$.

\end{itemize}
Hence we have found coverings of $U \cap \{x_{i_k} \neq 0\}$ for each $i_k \in I_i \setminus \bigcup_{j \neq i}I_j$ and a neighborhood of $U \cap S_0$ with the required properties. Then the result follows.
\qed

\begin{teo} \label{teo: hypsam} Let $\{\mu_i\}_{i=1}^\ell$, $\ell \leq n$ be actions of $\RP$ on $\R^n$ defined by
$$
\mu_j:((x_1,\ldots,x_n),c_j) \mapsto (c_{1j}x_1,\ldots,c_{nj}x_n),
$$
where $c_{ij}=1$ or $c_{ij}=c_j \in \RP$ satisfying \eqref{H} and \eqref{H3}. Then $\R^n$ satisfies \eqref{hypsam} (i)-(iii).
\end{teo}
\dim\ \ \eqref{hypsam} (i) follows from Theorem \ref{teo: hypsa}, \eqref{hypsam} (ii) follows from Lemma \ref{lem: RPl decomposition} and \eqref{hypsam} (iii) follows from the fact that $\mu_i$ is a globally subanalytic map for $i=1,\ldots,\ell$, hence if $U$ is globally subanalytic  $\mu_i(U,\RP)$ is still globally subanalytic. \qed

%
\section{Proof of Proposition \ref{m4.1.4}}
To prove the proposition, we need to show the following stronger assertion.
We first clarify the geometrical situations.
Let $Y$ be a real analytic manifold which is countable at $\infty$, and
set $X=Y_{\vartheta} \times \mathbb{R}_x^n$ with coordinates $(\vartheta,\, x)$
where $\vartheta$ denotes a point in $Y$ and $x = (x_1,\dots,x_n)$
are coordinates of $\mathbb{R}^n$.
Let us consider a family $\{I_1,\,\dots, I_\ell\}$ of subsets in $\{1,\dots,n\}$ which
satisfies the conditions (\ref{eq:conditions-indices-set}).
By replacing $Y$ with $Y \times \mathbb{R}^d$ if $d := n - \#(\cup_j I_j) > 0$,
we may assume
$$
\bigcup_{j=1,\dots,\ell} I_j = \{1,\dots,n\}.
$$
Then as $\hat{I}_1 \sqcup \dots \sqcup \hat{I}_\ell = \{1,\dots,n\}$ holds
(see (\ref{eq:def-hat-I}) for the definition of $\hat{I}_j$),
we have
$$
X = Y \times \mathbb{R}_x^n =
Y \times \mathbb{R}^{n_1}_{x_{\hat{I}_1}} \times \dots \times
\mathbb{R}^{n_\ell}_{x_{\hat{I}_\ell}},
$$
where $n_j = \#\hat{I}_j$ and $x_{\hat{I}_j}$ denotes the coordinates
$(x_i)$ with $i \in \hat{I}_j$ of $\mathbb{R}^{n_j}_{x_{\hat{I}_j}}$.

Let $V$ ($j=1, \dots, \ell)$ be a subanalytic subset in $X$.
We denote by $Z_j$ the closed subset
$$
Z_j := \{(\vartheta,\, x) \in X;\, x_{\hat{I}_j} = 0\}
$$
and by $\pi_j$ the canonical projection from $X$ to $Z_j$
by forgetting the variables $x_{\hat{I}_j}$.
We introduce the following conditions Va.~and Vb.~of $V$ for each $j$.
\begin{enumerate}
	\item[{\bf{Va.}}] $V$ does not intersect $Z_j$.
\item[{\bf{Vb.}}]
	$\pi_j(V) \subset \pi_j(V \cap Z_j)$.
\end{enumerate}


Let $\widetilde{X} = X \times \mathbb{R}_t^\ell$ with coordinates
$(\vartheta, x, t) = (\vartheta, x_1, \dots, x_n, t_1, \dots, t_\ell)$ and let us consider
the actions $\mu_j$ on $\widetilde{X}$ defined by
$$
\mu_j((\vartheta,x,t), c)
=\left(\vartheta,\, c_{1j}x_1,\dots,c_{nj}x_n,\, t_1, \dots, \frac{t_j}{c},\dots, t_\ell\right) \quad \text{for $c > 0$},
$$
where $c_{ij} =c$ if $i \in I_j$ and $c_{ij} = 1$ otherwise.

In what follows, we identify
$X$ with the closed subspace $\{t = 0\}$ of $\widetilde{X}$, and
in particular, $V$ is identified with the subset
$\{(\vartheta,x,t) \in \widetilde{X};\,
t = 0,\, (\vartheta,x) \in V\}$.
Then we have the following proposition.
\begin{prop}
Let $V$ be a $(\RP)^\ell$-conic subanalytic subset in $X$ and
let $W$ be an open subanalytic neighborhood of $V$ in $\widetilde{X}$.
If $V$ satisfies the condition either Va.~or Vb.~for each $j$,
then there exists a subanalytic subset $W'$ in $\widetilde{X}$
satisfying the following conditions.
\begin{enumerate}
\item  $W'$ is an open neighborhood of $V$ and contained in $W$.
\item  $W'$ is $(\RP)^\ell$-connected.
\item  $\RP_1\dots\RP_k W'$
is also subanalytic in $\widetilde{X}$ for any $1 \le k \le \ell$.
\end{enumerate}
\end{prop}

\dim\ \
As the proof is long, we briefly explain an outline of the proof.
We first consider the sequence of
morphisms and the spaces
$$
\widetilde{X} \setminus Z \overset{\varphi}{\longrightarrow} \hat{X}
\overset{\hat{\pi}_{t_1}}{\longrightarrow}
\widetilde{X}^\flat
$$
where $Z$ is some closed submanifold in $\widetilde{X}$.
The purpose of $\varphi$ is to transform an orbit of the $\RP$-action $\mu_1$ into
a line in $\hat{X}$. That of $\hat{\pi}_{t_1}$ is to forget the
action in $\hat{X}$ associated with $\mu_1$. Then we will construct, in each space,
the pair $\hat{V}$ and $\hat{W}$ in $\hat{X}$
(resp. $V^\flat$ and $W^\flat$ in $\widetilde{X}^\flat$)	
which corresponds to $V$ and $W$ respectively in the original space $\widetilde{X}$.
As the number of actions in $\widetilde{X}^\flat$
is smaller than that in $\widetilde{X}$, we can apply the induction for the number
of the actions to the pair $V^\flat$ and $W^\flat$. The important problem here is
to confirm that the pair in $\widetilde{X}^\flat$ satisfies the geometrical
situation to which we can apply the induction, and the most part of the proof
is devoted to this confirmation. We emphasize that, as there are many spaces and
symbols, we introduce the following rule. The objects
defined on the space $\hat{X}$ (resp. $\widetilde{X}^\flat$)
are denoted by symbols with $\hat{\cdot}$ (resp. $\cdot^\flat$) like
$\hat{V}$, $\hat{W}$, etc., (resp $V^\flat$, $W^\flat$, etc.).

\

Set $\chi = \{I_1, \dots, I_\ell\}$.
By the reordering of the indices for elements in $\chi$, we always assume
$$
j_1 < j_2 \implies I_{j_1} \subset I_{j_2} \text{ or } I_{j_1} \cap I_{j_2} = \emptyset,
$$
which can be achieved by using the function $d(\cdot)$
defined in the proof of Proposition \ref{prop:canonical-coordinate}.
Note that, by this reordering, $\hat{I}_1 = I_1$ holds.

We show the proposition by induction on
$\#\chi$, i.e., the number of elements in $\chi$.
If $\#\chi = 0$, then the proposition clearly holds by taking $W' = W$.
Suppose that the proposition holds for $\#\chi = 0,1,\dots,\ell-1$
with  $\ell \ge 1$. Then
we will show the proposition for $\#\chi = \ell$.

By considering
$$
W \cap
\{(\vartheta,\, x,\,t) \in \widetilde{X};\, |t_j| < 1,\, j=1,\dots,\ell\},
$$
we may assume that $\pi_t: \widetilde{X} \to X$ is proper
on $\overline{W}$, where $\pi_t$ is the canonical projection defined by
forgetting the variables $t =(t_1,\dots,t_\ell)$.

Let $x'$ be the coordinates
$(x_i)$ with $i \in I_1$ and $x''$ the other coordinates $(x_i)$ with $i \in \{1,\dots,n\} \setminus I_1$.  Set $m := n_1 = \#I_1$ and
$$
Z = \{(\vartheta, x',x'', t) \in \widetilde{X};\, x' = 0\} \subset \widetilde{X}.
$$
If $V$ satisfies the condition Va.~for $j=1$, by replacing $W$ with $W \setminus Z$,
we may assume $W \cap Z = \emptyset$. If $V$ satisfies the condition Vb.~for $j=1$,
then we replace $W$ with $W \cap (\pi_{x'}^{-1}(W \cap Z))$ where
$\pi_{x'}: \widetilde{X} \to Z = Y \times \mathbb{R}^{n-m}_{x''} \times \mathbb{R}_t^{\ell}
$ is the canonical projection defined by forgetting the variables $x'$.
Hence, in this case, we may assume
\begin{equation}
	\pi_{x'}(W) = W \cap Z.
\end{equation}
Note that, by the conditions Va.~and Vb.,  $W$ is still an open neighborhood
of $V$ in both cases.

Let us consider the diagram:
$$
\begin{matrix}
\widetilde{X}
& \overset{\rho}{\leftarrow}
& Y_{\vartheta} \times \mathbb{S}_{\theta}^{m-1}
\times \mathbb{R}_r \times \mathbb{R}_{x''}^{n-m} \times \mathbb{R}_t^{\ell}
& \overset{\omega}{\to}
& \hat{X} \\
\end{matrix}
$$
Here
$$
\hat{X} :=Y_{\vartheta} \times \mathbb{S}_{\theta}^{m-1}
\times \mathbb{R}_{x''}^{n-m} \times \mathbb{R}_s \times
\mathbb{R}_{t}^{\ell}.
$$
We will explain the morphisms $\rho$ and $\omega$.
The morphism $\rho$ is nothing but the polar coordinates transformation with respect to $x'$, that is,
$$
(\vartheta,\, x',\, x'',\, t) =
\rho(\vartheta,\, \theta,\, r,\, x'',\, t) =
(\vartheta,\, r\kappa(\theta),\, x'',\, t)
$$
where $\kappa: \mathbb{S}_{\theta}^{m-1} \hookrightarrow \mathbb{R}^{m}_{x'}$ is
the canonical inclusion.
Set $\hat{Y}_{\hat{\vartheta}} = Y_\vartheta \times \mathbb{S}_\theta^{m-1}$
with coordinates $\hat{\vartheta} = (\vartheta,\, \theta)$ for simplicity.
The morphism $\omega$ is defined by
$$
(\hat{\vartheta},\, x'',\, s,\, t_1,\, t')
= \omega(\hat{\vartheta},\, r,\, x'',\, t) =
(\hat{\vartheta},\, x'',\, r t_1,\, r,\, t'),
$$
where $t'=(t_2, \dots, t_\ell)$. We emphasize that the $t_1$ coordinate
in $\hat{X}$ is given by $r$ through the map $\omega$.
Define the morphism
$$
\varphi := \omega|_{\{r > 0\}} \circ (\rho|_{\{r > 0\}})^{-1} : \widetilde{X} \setminus Z \to \hat{X}
$$
and define, for any subanalytic subset $U \subset \widetilde{X}$,
$$
\varphi_{\#}(U) := \omega(\rho^{-1}(U) \cap \{r \ge 0\}).
$$
We have the following lemma.
\begin{lem}{\label{lem:prop-of-varphi}}
\begin{enumerate}
\item $\varphi$ is isomorphism between $\widetilde{X} \setminus Z$ and $\{t_1 > 0\}$ of $\hat{X}$, in particular,
$\varphi$ is an open map from $\widetilde{X} \setminus Z$ to $\hat{X}$.
\item For any subanalytic subset $U$ in $\widetilde{X}$ such that $\pi_t$ is proper on $\overline{U}$,
	the subset $\varphi(U \setminus Z)$ and $\varphi_{\#}(U)$ are
	subanalytic in $\hat{X}$.
\item For any subanalytic subset $\Omega$ in $\hat{X}$,
	the subset $\varphi^{-1}(\Omega) \subset \widetilde{X} \setminus Z$ is subanalytic in $\widetilde{X}$.
\end{enumerate}
\end{lem}
\dim \
The claim 1.~is easy to see. The claim 2.~for $\varphi(U \setminus Z)$
follows from the facts that
$$
\varphi(U \setminus Z) = \omega(\rho^{-1}(U \setminus Z) \cap \{r > 0\})
$$
and $\omega$ is proper on $\overline{\rho^{-1}(U)}$ by  assumption on $U$.
The claim 2.~for $\varphi_{\#}(U)$ is shown by the same reasons.
As we have
$$
\varphi^{-1}(\Omega) = \rho(\omega^{-1}(\Omega) \cap \{r > 0\})
$$
and $\rho$ is always proper, the claim 3.~follows.
\qed

Set
$$
\hat{W} := \varphi(W \setminus Z),\qquad
\hat{V} := \varphi_{\#}(V).
$$
It follows from 2.~of the above lemma that
these subsets are subanalytic in $\hat{X}$.
Note that $\hat{V}$ is a
subset of $\{s = 0,\, t' = 0\} \subset \hat{X}$
where $t' = (t_2, \dots, t_\ell)$.
We also define the action $\hat{\mu}_1$ on $\hat{X}$ by
$$
\hat{\mu}_1(p,c) := (\hat{\vartheta},\, x'',\, s,\, ct_1,\, t')
\text{ for $p = (\hat{\vartheta},\, x'',\, s,\, t_1, t') \in \hat{X}$ and $c > 0$},
$$
which coincides with $\varphi(\mu_1(\varphi^{-1}(p), c))$ for $p \in \{t_1 > 0\}$ of $\hat{X}$.
It is easy to see that $\hat{V}$ is an
$\RP$-conic subset with respect to $\hat{\mu}_1$.

Now we define the continuous subanalytic function $\sigma$ defined on the whole
$Y_{\vartheta} \times {\mathbb{R}_{x''}^{n-m}}$
as follows.
If $V$ satisfies the condition Va.~for $j=1$, then
we set
$$
\sigma(q) = 1, \qquad
q \in
Y_{\vartheta} \times {\mathbb{R}_{x''}^{n-m}}.
$$
If $V$ satisfies the condition Vb.~for $j=1$,
we first define $\sigma'$ on $Y_{\vartheta} \times {\mathbb{R}_{x''}^{n-m}}$
by
$$
\sigma'(\vartheta,\, x'') = \dfrac{1}{2}
\left(\inf_{
(\vartheta,\,x',\, x'',\,0_t) \in K
}
| x' |
\right),
$$
where we set
$$
K := (\widetilde{X} \setminus W) \cup
\{(\vartheta,\, x',\,x'',\, t) \in \widetilde{X};\, |x'| \ge 1\}.
$$
Then $\sigma'$ is a lower semi-continuous subanalytic function.
As it is lower semi-continuous, we can find a continuous subanalytic
function $\sigma$ on $Y_{\vartheta} \times {\mathbb{R}_{x''}^{n-m}}$
which satisfies $0 \le \sigma(q) \le \sigma'(q)$ for any
$q \in Y_{\vartheta} \times {\mathbb{R}_{x''}^{n-m}}$
and
$\sigma(q) > 0$ for $q$ with $\sigma'(q) > 0$.
Hence we have obtained $\sigma$ for the both cases.
Now we set
$$
\widetilde{X}^{\flat} := \hat{Y}_{\hat{\vartheta}} \times \mathbb{R}_{x''}^{n-m} \times \mathbb{R}_s \times \mathbb{R}_{t'}^{\ell-1}
$$
where $t' = (t_2, \dots, t_\ell)$ and we denote by
$$
\hat{\pi}_{t_1}: \hat{X} = \hat{Y}_{\hat{\vartheta}} \times \mathbb{R}_{x''}^{n-m} \times \mathbb{R}_s \times \mathbb{R}_{t}^{\ell}
\to
\widetilde{X}^{\flat} = \hat{Y}_{\hat{\vartheta}} \times \mathbb{R}_{x''}^{n-m} \times \mathbb{R}_s \times \mathbb{R}_{t'}^{\ell-1}
$$
the canonical projection by forgetting the variable $t_1$. We define
the closed subanalytic hypersurface in $\hat{X}$ by
$$
\hat{T}_\sigma := \{(\hat{\vartheta} = (\vartheta,\, \theta),\, x'',\, s,\,t) \in \hat{X};\, t_1
= \sigma(\vartheta,\, x'')\}.
$$
Then we set
$$
W^{\flat} := \hat{\pi}_{t_1}(\hat{W} \cap \hat{T}_\sigma), \qquad
V^{\flat} := \hat{\pi}_{t_1}(\hat{V}). 
$$
Note that
$V^\flat$ 
is a subset of $\{s = 0,\, t' = 0\} \subset \widetilde{X}^\flat$.
The following lemma is crucial for the induction process.
\begin{lem}{\label{lem:lem-induction-subanalytic-connect}}
\begin{enumerate}
\item The subsets $W^{\flat}$ and $V^{\flat}$ are
	subanalytic in $\widetilde{X}^{\flat}$. Further $W^\flat$ is open in $\widetilde{X}^\flat$.
\item $W^{\flat}$ is an open neighborhood of $V^{\flat}$ in $\widetilde{X}^{\flat}$.
\end{enumerate}
\end{lem}
\dim \
The $W^\flat$ is a subanalytic open subset as $\hat{\pi}_{t_1}$ is finite on
$\hat{T}_\sigma$.
Since $\hat{V} \subset \{t_1 \ge 0\}$ and
it is $\RP$-conic with respect to $\hat{\mu}_1$,
we have
\begin{equation}\label{eq:flat_pi_equality}
	V^\flat = \hat{\pi}_{t_1}(\hat{V}) =
	\hat{\pi}_{t_1}(\hat{V} \cap \{t_1 = 1\}) \cup
	\hat{\pi}_{t_1}(\hat{V} \cap \{t_1 = 0\}).
\end{equation}
Hence $V^{\flat}$ is subanalytic in $\widetilde{X}^\flat$.


Now we prove the claim 2.
It suffices to show that $V^\flat \subset W^\flat$ as $W^\flat$ has
been shown to be open.
We first assume that $V$ satisfies the condition Va.~for $j=1$.
Then, in this case, we have $\sigma = 1$ and $V = V \setminus Z$.
Hence we obtain
$$
\begin{aligned}
V^\flat &=
\hat{\pi}_t(\varphi_{\#}(V)) =
\hat{\pi}_t(\varphi_{\#}(V \setminus Z)) =
\hat{\pi}_t(\varphi(V \setminus Z)) \\
& = \hat{\pi}_t(\varphi(V \setminus Z) \cap \hat{T}_\sigma)
\subset \hat{\pi}_t(\varphi(W \setminus Z) \cap \hat{T}_\sigma) = W^\flat.
\end{aligned}
$$

Suppose that $V$ satisfies the condition Vb.~for $j=1$. Then we have
\begin{equation}{\label{eq:zerosection-V-reduction}}
V^\flat = \hat{\pi}_{t_1}(\hat{V} \cap \{t_1 = 0\}).
\end{equation}
As a matter of fact, given a point
$$
p = (\hat{\vartheta} = (\vartheta,\, \theta),\, x'',\, 0_s,\,1_{t_1},\, 0_{t'})
\in \hat{V} \cap \{t_1 = 1\},
$$
we have
$$
\varphi^{-1}(p) = (\vartheta,\, \kappa(\theta)_{x'},\, x'',\, 0_{t}) \in V.
$$
By  condition Vb.~for $j=1$, the point
$
(\vartheta,\, 0_{x'},\, x'',\, 0_{t})
$
also belongs to $V$.
Hence we have
$$
(\hat{\vartheta} = (\vartheta,\, \theta),\, x'',\, 0_s,\, 0_{t})
\in \hat{V} \cap \{t_1 = 0\},
$$
which implies (\ref{eq:zerosection-V-reduction}).
Then, to show $V^\flat \subset W^\flat$, for any
$$
p = (\hat{\vartheta} = (\vartheta,\, \theta),\, x'',\, 0_s,\, 0_{t})
\in \hat{V} \cap \{t_1 = 0\},
$$
it suffices to prove
$$
\hat{p} = (\hat{\vartheta},\, x'',\, 0_s,\,
\sigma(\vartheta, x'')_{t_1},\, 0_{t'})
\in \hat{W} \cap \hat{T}_\sigma.
$$
As $q = (\vartheta,\, 0_{x'},\, x'',\, 0_{t}) \in V$ and $q$
belongs to $W$, we have $\sigma'(\vartheta,\, x'') > 0$.
Therefore, by the definitions of $\sigma$ and $\sigma'$,
we have $\sigma(\vartheta,\, x'') > 0$ and
$$
(\vartheta,\, \sigma(\vartheta,\, x'')\kappa(\theta)_{x'},\,
x'',\, 0_{t}) \in W \setminus Z.
$$
This implies $\hat{p} \in \hat{W} \cap \hat{T}_\sigma$.
The proof has been completed.
\qed
\begin{oss}{\label{rem:positive-sigma-V}}
The following fact follows from the above proof.
For any point $p = (\hat{\vartheta} = (\vartheta,\, \theta),\, x'',\, 0_s,\, t_1, 0_{t'})
\in \hat{V}$, we have $\sigma(\vartheta,\, x'') > 0$.
\end{oss}

Now we apply the induction on the number of actions
to $W^{\flat}$ under the following geometrical situation.
Let $j_0$ be the minimal $j > 1$ such that $I_1 \subset I_j$.
If such a $j$ does not exist, set $j_0=+\infty$ by convention.
We consider
$$
Y^{\flat} := \hat{Y}_{\hat{\vartheta}} = Y \times \mathbb{S}_\theta^{m-1},\,\,
X^{\flat} := Y^{\flat} \times (\mathbb{R}_{x''}^{n-m} \times \mathbb{R}_s) \text{ and }
\widetilde{X}^{\flat} := X^{\flat} \times \mathbb{R}_{t'}^{\ell-1}
$$
if $j_0 < +\infty$
and
$$
Y^{\flat} := \hat{Y}_{\hat{\vartheta}} \times \mathbb{R}_s,\,\,
X^{\flat} := Y^{\flat} \times \mathbb{R}_{x''}^{n-m} \text{ and }
\widetilde{X}^{\flat} := X^{\flat} \times \mathbb{R}_{t'}^{\ell-1},
$$
if $j_0 = +\infty$, i.e., when there is no index $j > 1$ with $I_1 \subset I_j$.
We regard these $Y^\flat$, $X^\flat$ and $\widetilde{X}^\flat$
as $Y$, $X$ and $\widetilde{X}$ in the proposition respectively.

Then, to define $\ell-1$ actions on $\widetilde{X}^{\flat}$, we introduce
a family $I^{\flat}_2, \dots I^{\flat}_\ell$ of subsets in
$$
\left\{
\begin{array}{ll}
\{0,\, 1,\dots,n\} \setminus I_1,\qquad
&j_0 < +\infty, \\
\{1,\dots,n\} \setminus I_1, \qquad
&j_0 = +\infty.
\end{array}
\right.
$$
Here we note that the coordinates $(x'')$ of $\mathbb{R}_{x''}^{n-m}$
consists of $x_i$'s with $i \in \{1,\dots,n\} \setminus I_1$ and
the coordinate $s$ of $\mathbb{R}_s$ is also denoted by $x_0$ for convenience.
Set
$$
I_j^{\flat} := \left\{
\begin{array}{ll}
\{0\} \cup (I_j \setminus I_1), \qquad & \text{for $j > 1$ with $I_1 \subset I_j$}, \\
I_j, \qquad &\text{for $j > 1$ with $I_1 \cap I_j = \emptyset$}.
\end{array}
\right.
$$
Note that
$
I_j \subset I_1
$
never occurs for $j > 1$ by the ordering of $\chi$.
These $I^{\flat}_j$'s also satisfy the conditions (\ref{eq:conditions-indices-set}).
By using $I^{\flat}_j$'s,
we can define  $\ell-1$ actions
$\mu^{\flat}_2$, $\dots$, $\mu^{\flat}_\ell$ on $\widetilde{X}^{\flat}$
in the same way as that for $\mu_j$'s.
\begin{lem}
Assume $\ell > 1$.
Then $V^\flat$ is $(\RP)^{\ell-1}$-conic in $X^\flat$ and
it satisfies the condition either Va.~or Vb.~for each $j > 1$.
\end{lem}
\dim \
We first note
$$
X^\flat =
Y^\flat \times \mathbb{R}^{n^\flat_2}_{x_{\hat{I}_2^\flat}}
\times \dots
\times \mathbb{R}^{n^\flat_\ell}_{x_{\hat{I}_\ell^\flat}}
$$
where $x_{{\hat{I}^\flat_j}}$ denotes the coordinates $(x_i)$ with
$i \in \hat{I}^\flat_j$ and the coordinate $x_0$ is an alias of the coordinate $s$
as already noted.
In this situation,
the subset $Z^\flat_j$ is defined by $\{x_{\hat{I}^\flat_j} = 0\}$ ($j > 1$) and
$\pi^\flat_j$ is the canonical projection from $X^\flat$ to $Z^\flat_j$.
It follows from the definition that we have, for $j > 1$,
$$
\hat{I}^\flat_j = \left\{
\begin{array}{ll}
\hat{I}_j, \qquad & j \ne j_0, \\
\{0\} \cup \hat{I}_j, \qquad & j = j_0.
\end{array}
\right.
$$
Therefore we identify the space $\mathbb{R}^{n_j}_{x_{\hat{I}_j}}$
with $\mathbb{R}^{n^\flat_j}_{x_{\hat{I}^\flat_j}}$ for $j \ne j_0 > 1$
and
$\mathbb{R}^{n_j}_{x_{\hat{I}_j}}$
with the subspace $\{x_0 = 0\}$ (i.e., $\{s=0\}$) in
$\mathbb{R}^{n^\flat_j}_{x_{\hat{I}^\flat_j}}$ for $j = j_0$.
Then $V^\flat$ is given by
$$
\left\{
\begin{aligned}
(\hat{\vartheta} = (\vartheta,\, \theta),\, x'',\, s,\, t') \in \widetilde{X}^\flat;\,
& s = 0,\, t' = 0,\, \\
& \left(
\begin{aligned}
&(\vartheta,\, \kappa(\theta),\, x_{\hat{I}_2},\, \dots,\, x_{\hat{I}_j}) \in V
\text{ or }\\
&(\vartheta,\, 0,\, x_{\hat{I}_2},\, \dots,\, x_{\hat{I}_j}) \in V
\end{aligned}
\right)
\end{aligned}
\right\}
$$
where we use the identification
$x_{{\hat{I}^\flat}_j} = x_{\hat{I}_j}$ if $j \ne j_0 > 1$ and
$x_{{\hat{I}^\flat}_j} = (s,\, x_{\hat{I}_j})$ if $j = j_0$.
The claim of the lemma easily follows from this.
\qed

Therefore, by taking the above lemma and
Lemma \ref{lem:lem-induction-subanalytic-connect} into account,
we can apply the induction on $\#\chi$ to an open neighborhood
$W^\flat$ of $V^\flat$ in $\widetilde{X}^\flat$
under the geometrical situation already described.

Then we obtain
a subanalytic open neighborhood ${W^\flat}' \subset W^\flat$ of $V^\flat$
which satisfies the conditions 1., 2.~and 3.~of the proposition
for $I^{\flat}_j$ ($j=2,\dots,\ell$).
Set
$$
\hat{T}_{\sigma,\,{W^\flat}'} :=
\hat{\pi}^{-1}_{t_1}\left({W^\flat}'\right) \cap \hat{T}_\sigma\,
\subset \hat{W}
$$
which is a subanalytic subset in $\hat{X}$ and open in $\hat{T}_\sigma$. We also define
$$
\hat{W}'_{V \setminus Z} := \underset{p \in \hat{T}_{\sigma,\,{W^\flat}'}}{\cup} L(p)
\subset \hat{W}
$$
where $L(p)$ is a connected component of
$\hat{W} \cap \hat{\pi}_{t_1}^{-1}(\hat{\pi}_{t_1}(p))$ which contains
the point $p$.
\begin{oss}{\label{rem:property-of-T_sigma}}
	By definition, we have
$
\hat{T}_{\sigma,\,{W^\flat}'} \subset (\hat{W} \cap \hat{T}_\sigma)
\subset \hat{W} \subset \{t_1 > 0\}
$
in $\hat{X}$.  Hence, for any point
$
p = (\hat{\vartheta} = (\vartheta,\, \theta),\, x'',\, s,\, t_1, t')
\in \hat{T}_{\sigma,\,{W^\flat}'},
$
the subset $L(p)$ is not empty and $t_1 = \sigma(\vartheta,\, x'') > 0$ there.
\end{oss}
By the following lemma, $\hat{W}'_{V \setminus Z}$ is still an open subanalytic subset
in $\hat{X}$.

\begin{lem}\label{lem:projection-connected-fiber-subanalytic}
Let $X$ be a real analytic manifold and $U$ an open subanalytic subset in $X$. Let
$\pi: X \times \mathbb{R}_t \to X$ be the
canonical projection. Let $f(x)$ be a continuous subanalytic function on $X$.
Then, for any subanalytic open subset $W$ in $X \times \mathbb{R}_t$,
the set
$$
W' := \underset{p \in U}{\cup} L(p)
$$
is an open subanalytic subset in $X \times \mathbb{R}_t$
where $L(p)$
is a connected component of
$W \cap \pi^{-1}(p)$ which contains the point $(p,\, f(p))$.
\end{lem}
\dim \
The fact that $W'$ is open is clear by definition.
Since the problem is local with respect to $X$, we may assume
$X = \mathbb{R}^n$ and $U$ is a relatively compact open subanalytic subset
of $X$. Set, for $n=1,2,\dots$,
$$
W_n := W \cap \{(x,\,t) \in X \times \mathbb{R}_t;\, |t| < n\}.
$$
We denote by $W'_n$ the corresponding result of the operation described
in the lemma. Since $W_n$ is globally subanalytic in $X \times \mathbb{R}_t$,
$W'_n$ is subanalytic in $X \times \mathbb{R}_t$ by the same reasoning as that
in the proof of Lemma \ref{lem: RP decomposition}.
As $f(x)$ is a continuous function, the increasing sequence $\{W'_n\}$
is locally stationary in $X \times \mathbb{R}_t$. Hence
$
W' = \cup_n W'_n
$
is subanalytic in $X \times \mathbb{R}_t$.
\qed

Then we set
$$
W'_{V \setminus Z} := \varphi^{-1}(\hat{W}'_{V \setminus Z}) \subset W.
$$
It follows from 3.~of Lemma \ref{lem:prop-of-varphi} that
$W'_{V \setminus Z}$ is an open subanalytic subset in $\widetilde{X}$.
To show the properties of $W'_{V \setminus Z}$, we introduce
the following notation. Let $U$ be an open subanalytic subset of $\widetilde{X}$.
We define
$$
I_Z(U) := \left(\overline{U^C \setminus Z}\right)^C
$$
where $A^C$ denotes the complement of a subset $A$ in $\widetilde{X}$.
The following facts are easy to see.
\begin{enumerate}
	\item A point $p \in \widetilde{X}$ belongs to $I_Z(U)$
if and only if there exists an open neighborhood
$\Omega$ of $p$ such that $(\Omega \setminus Z) \subset U$.
\item $I_Z(U)$ is an open subanalytic subset in $\widetilde{X}$.
\item $I_Z(U) \setminus Z = U \setminus Z$.
\item If $\pi_{x'}(U) = U \cap Z$ holds, then $I_Z(U) = U$.
\end{enumerate}

\begin{lem}{\label{lem:W-is-open-neighbor-of-V}}
\begin{enumerate}
\item $W'_{V \setminus Z}$ is an open neighborhood of $V \setminus Z$.
\item If $V$ satisfies the condition Vb.~for $j=1$, then $I_Z(W'_{V \setminus Z})$
is contained in $W$ and it
is a subanalytic open neighborhood of $V \cap Z$.
\end{enumerate}
\end{lem}
\dim \
We first show the claim 1.
By Lemma \ref{lem:prop-of-varphi}, the map $\varphi$ is an isomorphism between
$\widetilde{X} \setminus Z$ and $\{t_1 > 0\} \subset \hat{X}$. We have
$$
\hat{W} = \varphi(W \setminus Z) \supset
\varphi(V \setminus Z) = (\varphi_{\#}(V) \setminus \{t_1 = 0\})
= (\hat{V} \setminus \{t_1 = 0\}).
$$
Hence it suffices to consider the corresponding problem in $\hat{X}$, that is,
we will show that
$\hat{V} \setminus \{t_1 = 0\} \subset \hat{W}'_{V \setminus Z}$.
Let
$$
p^* = (\hat{\vartheta}^{*} = (\vartheta^{*},\, \theta^{*}),\,
x''^{*},\, 0_s,\, t_1^{*},\, 0_{t'})
\in \hat{V} \setminus \{t_1 = 0\}.
$$
By Remark \ref{rem:positive-sigma-V}, we have
$\sigma(\vartheta^*,\, x''^*) > 0$.
We may assume $t_1^{*} \le \sigma(\vartheta^*,\, x''^*)$ as the proof
for the case $t_1^{*} \ge \sigma(\vartheta^*,\, x''^*)$ is the same.
%
As $\hat{V}$ is $\RP$-conic with respect to $\hat{\mu}_1$, the segment
$$
G := \{p \in \hat{X};\, p \in \hat{\pi}^{-1}_{t_1}(\hat{\pi}_{t_1}(p^*)),\,
t_1^* \le t_1 \le \sigma(\vartheta^*,\, x''^*)\}
$$
is contained in $\hat{V}$,
and hence, it is contained in $\hat{W}$ also.
Since $\hat{\pi}_{t_1}(p^*)$ belongs to $V^\flat \subset {W^\flat}'$,
we obtain $G \subset \hat{W}'_{V \setminus Z}$ and, in particular, $p^* \in \hat{W}'_{V \setminus Z}$.

\

Now let us show the claim 2. As $V$ satisfies
 condition Vb.~for $j=1$, we take $W$ so that $\pi_{x'}(W) = W \cap Z$ holds. Hence we
obtain
$$
W = I_Z(W) \supset I_Z(W'_{V \setminus Z}).
$$
This shows the first assertion of the claim 2.

We will show $V \cap Z \subset I_Z(W'_{V \setminus Z})$.
Suppose that a point
$$
q^* = (\vartheta^*,\, 0_{x'},\, x''^*,\, 0_t) \in V \cap Z
$$
does not belong to $I_Z(W'_{V \setminus Z})$.
Then there exist points $q^{(k)} \in \widetilde{X} \setminus Z$ $(k=1,2,\dots)$
such that $q^{(k)} \notin W'_{V \setminus Z}$ and
$q^{(k)} \to q^*$. These points are represented by their polar coordinates
$$
\tilde{q}^{(k)} = (\vartheta^{(k)},\, \theta^{(k)},\, r^{(k)},\, x''^{(k)},\, t^{(k)}) \in
Y_{\vartheta} \times \mathbb{S}_{\theta}^{m-1}
\times \mathbb{R}_r \times \mathbb{R}_{x''}^{n-m} \times \mathbb{R}_t^{\ell}
$$
satisfying $\rho(\tilde{q}^{(k)}) = q^{(k)}$ and
$$
\vartheta^{(k)} \to \vartheta^*,\quad
r^{(k)} \to 0\;(r^{(k)} > 0),\quad x''^{(k)} \to x''^*,\quad t^{(k)} \to 0.
$$
As $\mathbb{S}_{\theta}^{m-1}$ is compact, by taking a subsequence,
we may assume $\theta^{(k)} \to \theta^*$ for some
$\theta^*  \in \mathbb{S}_{\theta}^{m-1}$. Note that
$$
p^* := \omega(\vartheta^*,\, \theta^*,\, 0_r,\, x''^*,\, 0_t)
=(\vartheta^*,\, \theta^*,\, x''^*,\, 0_s,\, 0_{t})
$$
belongs to $\hat{V}$, and hence we have $\hat{\pi}_{t_1}(p^*) \in V^\flat$.

On the other hand, as $q^* \in V \subset W$,
it follows from the definition of $\sigma$ that
there exist $\epsilon > 0$ and
open neighborhoods  $U_\vartheta$, $U_{x''}$
of $\vartheta^*$, $x''^*$ in $Y$, $\mathbb{R}^{n-m}_{x''}$ respectively
such that the open subset
$$
U := \{\vartheta \in U_\vartheta,\, x'' \in U_{x''},\,
|x'| < \sigma(\vartheta^*,\, x''^*) + \epsilon,\,
|t_j| < \epsilon\, (j=1,\dots,\ell)\}
$$
of $\widetilde{X}$ is contained in $W$.
Then the image $\varphi(U \setminus Z)$ in $\hat{X}$ is given by
$$
\left\{
\begin{aligned}
&\vartheta \in U_\vartheta,\,
\theta \in \mathbb{S}^{m-1},\,
x'' \in U_{x''},\,
|t_j| < \epsilon\, (j=2,\dots,\ell) \\
& 0 < t_1 < \sigma(\vartheta^*,\, x''^*) + \epsilon,\,
|s| < \epsilon t_1
\end{aligned}
\right\}
$$
which is contained in $\hat{W}$.
Now we set
$$
\hat{U} := \left\{
\begin{aligned}
&\vartheta \in U_\vartheta,\,
\theta \in U_\theta,\,
x'' \in U_{x''},\,
|t_j| < \delta\, (j=2,\dots,\ell) \\
& 0 < t_1 < \sigma(\vartheta^*,\, x''^*) + \epsilon,\,
|s| < \epsilon t_1,\, |s| < \delta
\end{aligned}
\right\}
$$
for an open neighborhood $U_\theta$ of $\theta^*$ in $\mathbb{S}^{m-1}$
and $0 < \delta < \epsilon$ which is still contained in $\hat{W}$.
Then, as $\hat{\pi}_{t_1}(p^*) \in V^\flat \subset {W^\flat}'$,
if we take $U_\vartheta$, $U_\theta$, $U_{x''}$ and $\delta > 0$
sufficiently small, we have $\hat{\pi}_{t_1}(\hat{U}) \subset {W^\flat}'$.
Therefore, as $\sigma$ is continuous,
taking $U_\vartheta$, $U_{x''}$ and
$\delta > 0$ sufficiently small again,
we conclude that the subset $\hat{U}$ is contained in $\hat{W}'_{V \setminus Z}$.

Then the sequence
$$
\varphi(q^{(k)}) = \omega(\tilde{q}^{(k)}) =
\left(\vartheta^{(k)},\, \theta^{(k)},\, x''^{(k)},\,
(t_1^{(k)}r^{(k)})_{s},\, (r^{(k)})_{t_1},\, t'^{(k)}\right) \in \hat{X}
$$
clearly belongs to $\hat{U}$ for sufficiently large $k$, which
implies $\varphi(q^{(k)}) \in \hat{W}'_{V \setminus Z}$.
Hence we obtain $q^{(k)} \in W'_{V \setminus Z}$ for such a $k$.
This contradicts the assumption $q^{(k)} \notin W'_{V \setminus Z}$.
Therefore we have obtained the claim 2.  The proof has been completed.
\qed

It has been shown that $W'_{V \setminus Z}$ is
an open neighborhood of $V \setminus Z$.
We need, however, an additional subset to make
$W'_{V \setminus Z}$ an open neighborhood of $V$ if $V \cap Z$ is non-empty,
i.e., $V$ satisfies the condition Vb.~for $j=1$.
Note that it follows from the definition of $I_Z(\cdot)$ that,
for any open subanalytic subset $W_Z'$ in $Z$ with
$W_Z' \subset (I_Z(W'_{V \setminus Z}) \cap Z)$,
the subset $W_Z' \cup W'_{V \setminus Z}$ becomes an open subanalytic
subset in $\widetilde{X}$.
Therefore the problem can be reduced to find an appropriate open subanalytic subset
$W'_Z \subset (I_Z(W'_{V \setminus Z}) \cap Z)$
in  $Z$ which satisfies the conditions of the proposition.

\

%

Suppose that $V$ satisfies the condition Vb.~for $j=1$.
Let us consider the canonical projection
$$
\pi_{t_1}: Z =
Y \times \mathbb{R}_{x''}^{n-m} \times \mathbb{R}_{t}^{\ell}
\to
\widetilde{X}^\flat :=
Y \times \mathbb{R}_{x''}^{n-m} \times \mathbb{R}_{t'}^{\ell -1}
$$
defined by forgetting the variable $t_1$.
We set
$$
Y^\flat = Y,\quad X^\flat := Y^\flat \times \mathbb{R}_{x''}^{n-m},\quad
\widetilde{X}^\flat := X^\flat \times \mathbb{R}_{t'}^{\ell-1},
$$
and $T := \{t_1 = 0\}$ in $Z$.
We define, for $W_Z := I_Z(W'_{V \setminus Z}) \cap Z$,
$$
\begin{aligned}
	W^\flat &:= \pi_{t_1}( W_Z \cap T),\quad
	V^\flat := \pi_{t_1}(V \cap Z).
\end{aligned}
$$
These are subanalytic subsets in $\widetilde{X}^\flat$ and $W^\flat$ is
an open neighborhood of $V^\flat$
by 2.~of Lemma \ref{lem:W-is-open-neighbor-of-V}.
Now we define the
subsets $I^\flat_2$, \dots, $I^\flat_\ell$ of $\{1,\dots,n\} \setminus I_1$
by $I^\flat_j := I_j \setminus I_1$, $j=2,\dots,\ell$. These satisfy the
conditions (\ref{eq:conditions-indices-set}).
We also define
actions $\mu^\flat_j$ ($j=2,\dots,\ell$) by using $I^\flat_j$'s.
Clearly $V^\flat$ is $(\RP)^{\ell-1}$-conic with respect to these actions
and it satisfies either Va.~or Vb.~for each $j > 1$.

Hence we have constructed the geometrical situation to which the induction
can be applied.
Then, by applying the induction to $W^\flat$, we obtain a subanalytic
subset ${W^\flat}'$ which is an open neighborhood of $V^\flat$ and satisfies the
conditions 1., 2.~and 3.~of the proposition.

We set
$$
W'_Z := \underset{p \in T \cap \pi_t^{-1}({W^\flat}')}{\cup} L(p)
$$
where $L(p)$ is the connected component of
$W_Z \cap \pi^{-1}_{t_1}(\pi_{t_1}(p))$ which
contains the point $p$.
It follows from Lemma \ref{lem:projection-connected-fiber-subanalytic} that
$W'_Z$ is an open subanalytic subset of $Z$. Furthermore it
is an open neighborhood of $V \cap Z$ in $Z$ which is contained in
$W_Z$.

\

We have obtained the desired $W'$ of the proposition. As a matter of fact,
we set
$$
W' :=
\left\{
\begin{aligned}
	W'_{V \setminus Z},\qquad &\text{if $V$ satisfies the condition Va.~for $j=1$}, \\
	W'_{V \setminus Z} \cup W'_Z,\qquad &\text{if $V$ satisfies the condition Vb.~for $j=1$}.
\end{aligned}
\right.
$$
Then $W'$ is an open subanalytic neighborhood of $V$ and it is contained in $W$.
Hence, to complete the proof, it suffices to show
that $W'$ satisfies the conditions 2.~and 3.~of the proposition, which
comes from the following lemma.
The long proof of the proposition has been completed.
\qed

\begin{lem}
The subset $W'$ constructed above satisfies the
conditions 2.~and 3.~of the proposition.
\end{lem}
\dim \
Each $\RP$ action with respect to $\mu_j$ ($j=1,\dots.\ell$) is stable on
the regions $\widetilde{X} \setminus Z$ and $Z$ respectively.
Hence it suffices to show the claim for $W'_{V \setminus Z}$ and $W'_Z$
separately. Here we prove the claim only for $W'_{V \setminus Z}$.
The proof for $W'_Z$ is the same as that for
$W'_{V \setminus Z}$ and much easier.

We denote by $\hat{\pi}^{-1}_{t_1 > 0}(q)$ ($q \in \widetilde{X}^\flat$)
the subset $\hat{\pi}^{-1}_{t_1}(q) \cap \{t_1 > 0\}$ in $\hat{X}$ and
by $\hat{\mu}_j$ the action on $\{t_1 > 0\}$ of $\hat{X}$
$$
\hat{\mu}_j(p, c) := \varphi(\mu_j(\varphi^{-1}(p), c)),\qquad (c > 0).
$$
By noticing Remark \ref{rem:property-of-T_sigma},
we can easily see,
$$
\RP_{\mu_1} W'_{V \setminus Z} =
\varphi^{-1}(\hat{\pi}^{-1}_{t_1 > 0}({W^\flat}'))
$$
and, for $j > 1$ and any subset $A$ in  $\widetilde{X}^\flat$,
$$
\RP_{\hat{\mu}_j} \left(\hat{\pi}^{-1}_{t_1 > 0}(A)\right) =
\hat{\pi}^{-1}_{t_1 > 0}\left(\RP_{{\mu^\flat}_j}A\right).
$$
Therefore we get, for $1 \le k \le \ell$,
\begin{equation}{\label{eq:actions-commutative-relations}}
\RP_{\mu_1}\dots\RP_{\mu_k} W'_{V \setminus Z}
= \varphi^{-1}\left(\hat{\pi}^{-1}_{t_1 > 0}
\left(\RP_{{\mu^\flat}_2} \dots
\RP_{{\mu^\flat}_k} {W^\flat}'\right)\right)
\end{equation}
Then, by the induction hypothesis,
the subset $\hat{\pi}^{-1}_{t_1 > 0}\left(
\RP_{{\mu^\flat}_2} \dots
\RP_{{\mu^\flat}_k} {W^\flat}'\right)$ is subanalytic in $\hat{X}$.
Then, by 3.~of Lemma \ref{lem:prop-of-varphi}, the subset
$\RP_{\mu_1}\dots\RP_{\mu_k} W'_{V \setminus Z}$ is also subanalytic
in $\widetilde{X}$. This shows the condition 3.~of the proposition.

Now we show the condition 2.~of the proposition.
The subset $W_{V \setminus Z}'$ is clearly $\RP_{\mu_1}$-connected.
We can check easily that, for any $\RP_{\mu^\flat_j}$-connected
subset $A$ in $\widetilde{X}^\flat$ ($j > 1$), the subset
$\hat{\pi}^{-1}_{t_1 > 0}(A)$ is also
$\RP_{\hat{\mu}_j}$-connected in $\hat{X}$.
By the induction hypothesis, the subset
$\RP_{{\mu^\flat}_2} \dots
\RP_{{\mu^\flat}_k} {W^\flat}'$
is $\RP_{\mu^\flat_{k+1}}$-connected.
Hence, by (\ref{eq:actions-commutative-relations}),
the subset $\RP_{\mu_1}\dots\RP_{\mu_k} W'_{V \setminus Z}$
is $\RP_{\mu_{k+1}}$-connected.
This completes the proof.
\qed

\end{section}
\addcontentsline{toc}{section}{

\noindent
\parbox[t]{.48\textwidth}{
Naofumi HONDA \\
Department of Mathematics, \\
Faculty of Science, \\
Hokkaido University, \\
060-0810 Sapporo, Japan. \\
honda@math.sci.hokudai.ac.jp
} \hfill
\parbox[t]{.48\textwidth}{
Luca PRELLI\\
Dipartimento di Matematica \\
  Pura ed Applicata,\\
Universit\`{a} degli studi di Padova,\\
Via Trieste 63,\\
35121 Padova, Italia. \\
lprelli@math.unipd.it }

\end{document}